\documentclass[12pt]{amsart}
  \usepackage{latexsym}
   \usepackage{amsmath}
   \usepackage{amssymb}
   \usepackage{amscd}
\usepackage{epsfig, graphics}
\usepackage{hyperref}

\input amssym.def
\input amssym.tex

\newcommand{\nc}{\newcommand}

\newcommand{\End}{\,{\rm End}\,}

\newcommand{\Sym}{\,{\rm Sym}\,}

\newcommand{\vol}{ \,{\rm vol}\, }

\newcommand{\beq}{\begin{equation}}
\newcommand{\eeq}{\end{equation}}
\newcommand{\beqst}{\begin{equation*}}
\newcommand{\eeqst}{\end{equation*}}
\newcommand{\barr}{\begin{array}}
\newcommand{\earr}{\end{array}}
\newcommand{\beqar}{\begin{eqnarray}}
\newcommand{\eeqar}{\end{eqnarray}}
\newtheorem{theorem}{Theorem}
\newtheorem{corollary}[theorem]{Corollary}
\newtheorem{lemma}[theorem]{Lemma}
\newtheorem{prop}[theorem]{Proposition}
\newtheorem{proposition}[theorem]{Proposition}
\newtheorem{definition}[theorem]{Definition}
\newtheorem{remit}[theorem]{Remark}

\newtheorem{examplit}[theorem]{Example}

\newcommand{\cal}{\mathcal}

\newenvironment{rem}{\begin{remit}\rm}{\end{remit}}

\newcommand{\RR}{{\Bbb R }}
\newcommand{\CC}{{\Bbb C }}
\newcommand{\ZZ}{{\Bbb Z }}
\newcommand{\PP}{ {\Bbb P } }

\newcommand{\UU}{{\Bbb U }}

\newcommand{\cald}{{\mbox{$\cal D$}}}

\newcommand{\calf}{{\mbox{$\cal F$}}}

\newcommand{\cali}{{\mbox{$\cal I$}}}

\newcommand{\calm}{{\mbox{$\cal M$}}}

\def\a{\alpha}
\def\b{\beta}
\def\g{\gamma}

\def\e{\epsilon}

\def\k{\kappa}
\def\l{\lambda}
\def\m{\mu}

\def\G{\Gamma}

\def\L{\Lambda}

\nc{\Proof}{\noindent{\em Proof:} }
\nc{\Lie}{ {\rm Lie} }

\nc{\liek}{{ \bf k} }
\nc{\lieks}{\liek^*}
\nc{\hk}{H^*_K}
\nc{\xvec}{X}
\nc{\hfeta}{h_F^{\eta}}
\nc{\mnd}{\calm(n,d)}

\nc{\lieksa}{\lieks_{{\rm ad}} }
\nc{\lieka}{\liek_{{\rm ad}} }
\nc{\liekas}{(\liek_s)_{{\rm ad}} }
\nc{\lieksas}{(\lieks_s)_{{\rm ad}} }
\nc{\liekaq}{(\liek/\liek_s)_{{\rm ad}} }
\nc{\lieksaq}{(\lieks/\lieks_s)_{{\rm ad}} }
\nc{\conv}{{\rm Conv}}

\nc{\proj}{\pi}
\nc{\zloc}{\mu^{-1}(0)}
\nc{\nf}{\nu_F}
\nc{\lx}{\l(X) }
\nc{\rlx}{\frac{d \l(X)}{\l(X) }  }
\nc{\srlx}{ {d \l(X)} / {\l(X) }  }
\nc{\intlat} {\L^I}

\nc{\sumpl}{ { \sum_{F \in \calf_+} } }

\nc{\liet}{{\bf t}}
\nc{\nusym}{\cald}
\nc{\hfe}{h_F^{\eta} }
\nc{\liets}{\liet^*}
\nc{\expmf}{e^{\isq \mu_T(F)(X) } }
\nc{\expbfj}{e^{\isq \bfj (X) } }
\nc{\expmfs}{e^{\isq \mu_T(F)X } }
\nc{\sumf}{\sum_{F \in \calf_+} }
\nc{\mf}{\mu_T(F)}
\nc{\tck}{{C^K_1}}
\nc{\cko}{C^K_2}
\nc{\bom}{{\bar{\omega}}}
\nc{\fk}{F_K}
\nc{\ft}{F_T}

\nc{\sss}{s}
\nc{\wf} { { { \cal W}_F } }

\nc{\bfj}{ {\beta_{F,j} } }
\nc{\nfj}{\nu_{F,j} }
\nc{\res}{{\rm res} }

\nc{\stabo}{n_0}
\nc{\Jac}{\Delta}

\nc{\xg}{M/\!/G}
\nc{\xtg}{\tilde{M}/\!/G}
\nc{\txg}{\tilde{M}/\!/G}
\nc{\mred}{\xg}
\nc{\xred}{\xg}

\newcommand{\stab}{{\rm Stab}}

\newcommand{\git}{/\!/}
\newcommand{\quott}{\git}
\newcommand{\hT}{H_T}

\newcommand{\mlnd}{{\mathcal M}_\Lambda(n,d) }

\newcommand{\renorm}{{ \setcounter{equation}{0} }}

\begin{document}

\title[Intersection pairings on singular moduli spaces]{Intersection pairings on
singular moduli spaces of bundles over a Riemann surface and their
partial desingularisations}
\date{\today}

\begin{abstract}
This paper studies intersection theory on the compactified moduli
space $\calm (n,d)$ of holomorphic bundles of rank $n$ and degree
$d$ over a fixed compact Riemann surface $\Sigma$ of genus $g \geq
2$ where $n$ and $d$ may have common factors. Because of the
presence of singularities we work with the intersection cohomology
groups $I\!H^*(\calm (n,d))$ defined by Goresky and MacPherson and
the ordinary cohomology groups of a certain partial resolution of
singularities $\tilde{\calm} (n,d)$ of $\calm (n,d)$. Based on our
earlier work \cite{jkkw}, we give a precise formula for the
intersection cohomology pairings and provide a method to calculate
pairings on $\tilde{\calm}(n,d)$. The case when $n=2$ is discussed
in detail. Finally Witten's integral is considered for this
singular case.
\end{abstract}

\address{{\rm Lisa C. Jeffrey}, Department of Mathematics, University of Toronto,
Toronto ON, M5S 2E4, Canada} \email{jeffrey@math.toronto.edu}
\address{{\rm Young-Hoon Kiem}, Department of Mathematics, Seoul
National University, Seoul 151-747,
Korea}\email{kiem@math.snu.ac.kr}
\address{{\rm Frances C. Kirwan}, Mathematical Institute, Oxford University
 OX1 3LB, UK}\email{frances.kirwan@balliol.oxford.ac.uk}
\address{{\rm Jonathan M. Woolf}, Department of Mathematical
Sciences, Liverpool, L69 7ZL, UK}\email{jonathan.woolf@liverpool.ac.uk}

\author [Jeffrey, Kiem, Kirwan and Woolf]
{L.C. Jeffrey, Y.-H. Kiem, F.C. Kirwan, J. Woolf}

\thanks{LCJ was partially supported by a grant from NSERC.}
\thanks{YHK was partially supported by KOSEF
R01-2003-000-11634-0.}

\maketitle

\renorm
\section{Introduction}

This paper studies intersection theory on the compactified moduli
spaces $\calm (n,d)$ and $\calm_\Lambda (n,d)$ of holomorphic
bundles of rank $n$ and degree $d$ over a fixed compact Riemann
surface $\Sigma$ of genus $g \geq 2$ (and with fixed determinant
line bundle $\Lambda$ in the case of $\calm_\Lambda (n,d)$), and
their partial desingularisations $\tilde{\calm}(n,d)$ and
$\tilde{\calm}_\Lambda (n,d)$ in the sense of \cite{K4,K5}. Here
$n$ and $d$ may have common factors so that $\calm (n,d)$ and
$\calm_\Lambda (n,d)$ may be singular; when $n$ and $d$ are
coprime then $\tilde{\calm}(n,d) = \calm (n,d)$ and
$\tilde{\calm}_\Lambda (n,d) = \calm_\Lambda (n,d)$ and the
results of this paper have already been obtained in \cite{JK2}.
Indeed, intersection theory on these moduli spaces when $n$ and
$d$ are coprime has been studied intensively for several decades
 \cite{AB,Be,DR,D,G,HN,HS,JK2,KN,LP,NR,NS,
TelemanWoodward,T,tdgr,Z}; more recently work has also been done
on the singular\footnote{The moduli spaces $\calm(n,d)$ and
$\calm_\Lambda (n,d)$ are singular if and only if $n$ and $d$ have
a common factor, except in one special case: when $g=2$ then
$\calm(2,d)$ and $\calm_\Lambda (2,d)$ are nonsingular when $d$ is
even as well as when $d$ is odd.} moduli spaces, in particular
$\calm (2,d)$ when $d$ is even \cite{Kiem,Kiem2,Kiem3}. Our aim
here is to extend the results of \cite{JK2} on intersection
pairings in the cohomology of $\mlnd$ to the case when $n$ and $d$
are not coprime, by using the methods of \cite{jkkw}. Because of
the presence of singularities we work with the intersection
cohomology groups $I\!H^*(\calm (n,d))$ and $I\!H^*(\calm_\Lambda
(n,d))$ defined by Goresky and MacPherson \cite{GM1,GM2} and the
ordinary cohomology groups of the partial resolutions of
singularities $\tilde{\calm} (n,d)$ and $\tilde{\calm}_\Lambda
(n,d)$ of $\calm (n,d)$ and $\calm_\Lambda (n,d)$ \cite{K4,K5}.
The formulas we obtain for the intersection pairings in
$I\!H^*(\calm (n,d))$ and $I\!H^*(\calm_\Lambda (n,d))$ (see Theorem
\ref{mthmsect6}) follow easily from the results of
\cite{JK2,Kiem,jkkw} and are essentially the same as in the
coprime case studied in \cite{JK2}. As was shown in \cite{Kiem}
these pairings can be regarded as providing some of the
intersection pairings in $H^*(\tilde{\calm} (n,d))$ and
$H^*(\tilde{\calm}_\Lambda (n,d))$. The complete picture of the
pairings in $H^*(\tilde{\calm} (n,d))$ and
$H^*(\tilde{\calm}_\Lambda (n,d))$ is much more complicated to
describe; we give a method for calculating these (see Theorem
\ref{mthmsect6}) and carry out the case when $n=2$ in detail (see
Theorem \ref{mthmsect8}).

The original motivation for both \cite{JK2} and this article (as
well as much other work) was \cite{tdgr}, where Witten studied the
moduli spaces $\calm (n,d)$ as symplectic reductions of infinite
dimensional affine spaces by infinite dimensional groups following
\cite{AB}. He found formulas for intersection pairings on these
moduli spaces for coprime $n$ and $d$ from asymptotic expansions
of certain infinite dimensional integrals as a parameter $\e$
tends to $0$. He did this by showing that each integral is a sum
of terms tending to $0$ exponentially fast with $\e$, together
with a polynomial in $\e$ whose coefficients are intersection
pairings on the moduli space. Witten's formulas were later proved
using finite-dimensional methods in \cite{JK2} for $\calm(n,d)$
and $\calm_\Lambda(n,d)$ with coprime $n$ and $d$, and for moduli
spaces of principal bundles for more general compact groups in
\cite{Meinrenken}. Witten also gave formulas for the asymptotic
expansions of his integrals in the case of bundles of rank two and
even degree, and noted that powers of $\e^{1/2}$ appeared. In fact
(cf. \cite{jkkw,paradan,Woodward,BeasleyWitten}) when $n$ and $d$
are not coprime Witten's integrals are always sums of polynomials
in $\e^{1/2}$ rather than $\e$, together with terms tending to $0$
exponentially fast with $\e$. Some of the coefficients of integral
powers of $\e$ in these expressions can be interpreted as
intersection pairings; however it is not clear whether there is a
geometrical interpretation of the coefficients of the
half-integral powers of $\e$.

The layout of this paper is as follows. In $\S$2 we recall the
results we shall need from \cite{jkkw} on intersection pairings on
geometric invariant theoretic quotients and symplectic reductions.
$\S$3 recalls facts about moduli spaces of vector bundles over
curves and their partial desingularisations. In $\S$4 we review
the finite-dimensional methods used in \cite{JK2} to rederive
Witten's formulas in the case when $n$ and $d$ are coprime, and in
$\S$5 we combine these methods with the results of \cite{jkkw} to
calculate the pairings in $I\!H^*(\calm(n,d))$ in Theorem
\ref{mthmsect6}. The pairings in the intersection homology $I\!H^*(\calm_\Lambda(n,d))$ of the moduli space of bundles with fixed determinant can be obtained from these calculations. Indeed, the latter part of the paper deals only with $\calm(n,d)$ and $\tilde \calm(n,d)$ since computations for $\calm_\Lambda(n,d)$ and $\tilde \calm_\Lambda(n,d)$ can be easily derived from these, see Remark \ref{r:finitecover}. $\S$6 begins the much more laborious task of
extending these calculations to cover all the pairings in
$H^*(\tilde{\calm}(n,d))$, by expressing such pairings as sums of
formulas like those already seen together with certain \lq
wall-crossing terms' given by integrals over symplectic quotients
of projective bundles (Theorem \ref{mthmsect6}). $\S$7 completes
the calculation of pairings on $\tilde{\calm}(n,d)$ by explaining
how to use the techniques of \cite{jkkw} to compute the
wall-crossing terms inductively via integration over the fibres;
unfortunately this becomes very cumbersome in practice for large
$n$. In $\S$8 we consider the case when $n=2$ in detail; here it
is not hard to give explicit formulas for the pairings in
$H^*(\tilde{\calm}(2,d))$ as well as in $I\!H^*(\calm(2,d))$ (see
Theorem \ref{mthmsect8}). Finally in $\S$9 we look at Witten's
integrals.

\renorm

\section{Pairings on singular quotients}

Let $\xg$ be the quotient in the sense of Mumford's geometric
invariant theory \cite{MFK} of a nonsingular connected complex
projective variety $M$ by a linear action of a connected complex
reductive group $G$. In \cite{jkkw} we gave formulas, under
certain conditions on the group action, for the pairings of
intersection cohomology classes of complementary degrees in the
intersection cohomology $I\!H^*(\xg)$ of $\xg$. (Intersection
cohomology is defined with respect to the middle perversity
throughout this paper, and all cohomology and homology groups have
complex coefficients). We also gave formulas for intersection
pairings on resolutions $\txg$ (or more precisely partial
resolutions, since orbifold singularities are allowed) of the
quotients $\xg$.

Recall that if every semistable point of $M$ is stable then $0$ is
a regular value of the moment map and the stabiliser in the
maximal compact subgroup $K$ of  $G$ of every point of $\zloc$ is
finite.
 This
implies that the cohomology $H^*(\xg)$ of the quotient $\xg$ is
naturally isomorphic to the equivariant cohomology $\hk(\zloc)
\cong \hk (M^{ss})$ (recall that we are working with cohomology
with complex coefficients). The restriction map $\hk(M) \to
\hk(M^{ss})$ is surjective (\cite{K2} 5.4), and so the composition
of the restriction map $\hk(M) \to \hk(M^{ss})$ and the
isomorphism $\hk(M^{ss}) \to H^*(\xg)$ gives us a natural
surjective ring homomorphism \beq \label{1.1} \k_M :\hk(M) \to
H^*(\xg). \eeq Since $H^*(\xg)$ satisfies Poincar\'{e} duality,
the kernel of this surjection is then determined by the formula
obtained in \cite{JK1} (see (\ref{res}) below) for pairings of
cohomology classes of complementary dimensions in $\xg$ in terms
of equivariant cohomology classes in $M$ which represent them.

If there are semistable points of $M$ which are not stable (we
assume only that there do exist some stable points, or
equivalently that there exist some points in $\zloc$ where the
derivative of $\mu$ is surjective) then there is no longer a
natural surjection from $\hk(M)$ to $H^*(\xg)$, and $\xg$ is in
general singular so its cohomology $H^*(\xg)$ may not satisfy
Poincar\'{e} duality. However its intersection cohomology groups
satisfy
Poincar\'{e} duality, and there is a surjection from $\hk(M)$ to
the intersection cohomology $I\!H^*(\xg)$, which we will call
$\k_M$ since it coincides with (\ref{1.1}) when semistability is
the same as stability. This surjection $\k_M: \hk(M) \to
I\!H^*(\xg)$ arises as follows.

We can construct a canonical partial resolution of singularities
$\txg$ of the quotient $\xg$ (see \cite{K4}), by blowing $M$ up
along a sequence of nonsingular $G$-invariant subvarieties, all
contained in the complement $M-M^s$ of the set $M^s$ of stable
points of $M$. This eventually gives us a nonsingular projective
variety $\tilde{M}$ with a linear $G$-action, lifting the action
on $M$, for which every semistable point of $\tilde{M}$ is stable.
The quotient $\txg$ has only orbifold singularities, and the
blowdown map $\pi:\tilde{M} \to M$ induces a birational morphism
$$\pi_G:
\txg \to \xg$$ which is an isomorphism over the dense open subset
$M^s/G$ of $\xg$. The intersection cohomology $I\!H^*(\xg)$ of
$\xg$ is a direct summand of the cohomology of this partial
resolution of singularities $\txg$, and the composition \beq
\label{eq2} \k_M: \hk(M) \to \hk(\tilde{M}) \to H^*(\txg) \to
I\!H^*(\xg) \eeq is surjective (\cite{K6,Woolf2}). In fact this
surjection is the composition of the restriction map from $\hk(M)$
to $\hk(M^{ss})$ and a surjection $\k_M^{ss}: \hk(M^{ss}) \to
I\!H^*(\xg)$.

The work of the second author \cite{Kiem} allows us to understand
pairings in $I\!H^*(\xg)$ of intersection cohomology classes on
the singular quotient $\xg$ in terms of this surjection
$\k_M^{ss}: \hk(M^{ss}) \to I\!H^*(\xg)$. It is shown in
\cite{Kiem} that if the action of $G$ on $M$ is weakly balanced (a
condition satisfied in the case of moduli spaces of bundles over
Riemann surfaces: see \cite{Kiem}), then there is a naturally
defined subset $V_M$ of $\hk(M^{ss})$ such that $\k_M^{ss}:
\hk(M^{ss}) \to I\!H^*(\xg)$ restricts to an isomorphism \beq
\label{kmss} \k_M^{ss}: V_M \to I\!H^*(\xg). \eeq It is also shown
in \cite{Kiem} that the intersection pairing of two elements
$\k_M(\alpha)$ and $\k_M(\beta)$ of complementary degrees in
$I\!H^*(\xg)$ is equal to the evaluation of the image in
$H^*(\tilde{M}/\!/G)$ of the product $\alpha\beta \in \hk(M^{ss})$
on the fundamental class $[\tilde{M}/\!/G]$, provided that
$\alpha$ and $\beta$ lie in $V_M$. This means that we can study
intersection pairings in $I\!H^*(\xg)$ via intersection pairings
in the ordinary cohomology of the quotient $\txg$ for which
semistability is the same as stability.

The residue formula of \cite{JK1}\footnote{See Theorem 3.1 of
\cite{JK2} for a corrected version.}
 is a formula,
in the case when semistability equals stability, for pairings
 of cohomology classes $\k_M(\a)$ and
$\k_M(\beta)$ of complementary dimensions in $\xg$ in terms of
equivariant cohomology classes $\a$ and $\beta$ in $M$ which
represent them. Let $T$ be a maximal torus in $K$; its
complexification $T_c$ is then a maximal complex torus of $G$. Let
$\Gamma$ be the set of roots of $K$ regarded as elements of the
dual $\liets$ of the Lie algebra $\liet$ of $T$, and let $\G_+$
and $\G_-$ be the subsets of $\G$ consisting of the positive and
negative roots of $K$. Let $\bom = \omega + \mu$ be the
standard\footnote{Here we follow the conventions of
\cite{JK2,jkkw} which differ slightly from those used in
\cite{JK1} and by Witten in \cite{tdgr}; in particular we have no
factors of $i$.} extension of the symplectic form $\omega$ to an
equivariantly closed differential form on $M$. We shall assume for
simplicity throughout that the stabiliser in $K$ of a generic
point of $\zloc$ is trivial. Then if $\calf$ is the set of
components of the fixed point set $M^T$ of $T$ on $M$, the residue
formula is \beq \label{res} \k_M(\alpha) \k_M(\beta)[\xg] = \frac{
(-1)^{s+n_+}}{|W| \vol(T)} \res \bigl( \cald(\xvec)^2 \sum_{F \in
\calf} \int_F \frac{i_F^*(\alpha \beta
e^{\bom})(\xvec)}{e_F(\xvec)} [d\xvec] \bigr ), \eeq where
$\vol(T)$ and $[d\xvec]$ are the volume of $T$ and the measure on
its Lie algebra $\liet$ induced by the restriction to $\liet$ of
the fixed inner product on $\liek$, while $W$ is the Weyl group of
$K$, the polynomial function $\cald(\xvec)= \prod_{\g \in \G_+}
\g(X)$ of $\xvec \in \liet$ is the product of the positive
roots\footnote{In this paper, as in \cite{JK2}, we adopt the
convention that weights $\beta \in \liets$ send the integer
lattice $\L^I = \ker(\exp:\liet \to T)$ to $\ZZ$ rather than to
$2\pi \ZZ$, and that the roots of  $K$ are the nonzero weights of
its complexified adjoint action. This is one reason why the
constant in the residue formula above differs from that of
\cite{JK1} Theorem 8.1.} of $K$ and $n_+ = (s-l)/2$ is the number
of those positive roots\footnote{Notice that
$(-1)^{n_+}(\cald(X))^2= \prod_{\g \in \G}\g(X)$.}; $s$ is the
dimension  of $K$ and $l$ is the dimension of $T$. Also if $F \in
\calf$ is a component of the fixed point set $M^T$ then $i_F:F \to
M$ is the inclusion and $e_F$ is the equivariant Euler class of
the normal bundle to $F$ in $M$. The multivariable residue map
$\res$ which appears in the formula is a linear map, but in order
to apply it to the individual terms in the formula some choices
have to be made which do not affect the residue of the whole sum.
Once the choices have been made, many of the terms in the sum have
residue zero and the formula can be rewritten as a sum over a
subset $\calf_+$ of $\calf$. When $T$ has dimension one (which is
the only case we shall need explicitly in this paper, because we
shall be using an inductive argument modelled on that of
\cite{JK2}) we can take $\calf_+$ to consist of those $F \in
\calf$ on which the constant value taken by the $T$-moment map
$\mu_T:M\to\liets \cong \RR$ is positive, and then $\res$ applied
to those terms in the sum labelled by $F\in \calf_+$ is the usual
one-variable residue $\res_{\xvec =0}$ at $0$. Indeed for $K =
U(1)$ we have \beq \label{eq3.4}  \k_M(\alpha) \k_M(\beta) [\xg] =
 \res_{\xvec=0} \Bigl (  \sum_{F \in \calf_+}
 \int_F \frac{i^*_F\alpha (\xvec) \beta (\xvec)}{e_F(\xvec)}
\Bigr )  \eeq (see \cite{JK1,Kalkman,Wu}) where $\res_{\xvec=0} $
denotes the coefficient of $1/\xvec$ when $\xvec \in \RR$ has been
identified with $2 \pi i \xvec \in \liek$.

Thus, in principle, when $M$ has semistable points which are not
stable we can apply the residue formula (\ref{res}) above to
$\tilde{M}$ to obtain pairings on the partial desingularisation
$\xtg$ and in the intersection cohomology $I\!H^*(\xg)$ of the
singular quotient $\xg$. When the action of $G$ on $M$ is weakly
balanced  (which will be the case for the actions to be considered
in this paper), so that $\k_M^{ss} : \hk(M^{ss}) \to I\!H^*(\xg)$
restricts to an isomorphism $\k_M^{ss}: V_M \to I\!H^*(\xg)$ as at
(\ref{kmss}), and if  $\alpha|_{M^{ss}}$ and $\beta|_{M^{ss}}$ lie
in $V_M$ and $\k_M(\alpha)$ and $\k_M(\beta)$ have complementary
degrees, then their intersection pairing  in $I\!H^*(\xg)$ is
 given
by \beq \label{((3.10))} \langle
\k_{{M}}(\a),\k_{{M}}(\b)\rangle_{I\!H^*(\xg)} =\frac{
 (-1)^{n_+}}{
 |W| } \k_M^{(\delta)} (\alpha \beta \cald^2)
[\mu_T^{-1}(\delta)/T] \eeq and thus by \beq \label{thm8.4}
\langle \k_{{M}}(\a),\k_{{M}}(\b)\rangle_{I\!H^*(\xg)} =\frac{
 (-1)^{s+n_+}}{|W| \vol(T)} \res \bigl( \cald(\xvec)^2
\sum_{{F} \in {\calf}} \int_{{F }} \frac{i_{{F}}^*(\alpha \beta
e^{\bom-\delta})(\xvec)}{e_{F}(\xvec)} [d\xvec] \bigr ),\eeq for
any sufficiently small $\delta \in \liets$ which is a regular
value of the moment map $\mu_T$ (see \cite{jkkw} Theorem 18). Here
$\k_M^{(\delta)}$ is defined as at (\ref{1.1}) but with $K$
replaced by $T$ and the moment map $\mu$ replaced by $\mu_T -
\delta$. Moreover in fact (\ref{((3.10))}) is valid even when $M$
is not compact and has singularities away from $\zloc$, provided
that $\zloc$ is compact and $M$ is smooth near $\zloc$.

\newcommand{\tmg}{\tilde{M}/\!/G}

Thus the intersection pairings in $I\!H^*(\xg)$ are given by a
very simple modification of the residue formula (\ref{res}).
However when we try to apply (\ref{res}) to $\tilde{M}$ to obtain
pairings on the partial desingularisation $\tmg$ then
complications arise. The main difficulty is that, although the
construction of $\tilde{M}^{ss}$ and of $\txg$ from the linear
$G$-action on $M$ is canonical and explicit, the construction of
$\tilde{M}$ is not. In \cite{K4} the set $M^{ss}$ of semistable
points of $M$ is blown up along a sequence of nonsingular
$G$-invariant closed subvarieties $V$, and after each blow-up any
points which are not semistable are thrown out, so that eventually
we arrive at $\tilde{M}^{ss}$ and thus obtain $\xtg =
\tilde{M}^{ss}/G$. If necessary $\tilde{M}$ itself could be
constructed by resolving the singularities of the closures
$\bar{V}$ of these subvarieties $V$ and blowing up along their
proper transforms, but in practice this is not usually simple.
Unfortunately the residue formula of \cite{JK1} involves the set
of components of the fixed point set of the action of the maximal
torus $T$ of $K$, so applying it directly to $\tilde{M}$ would be
very complicated, and knowledge of the set of semistable points
$\tilde{M}^{ss}$ alone would not suffice.  However there is an
alternative way to calculate the pairings which only requires
information about $\tilde{M}^{ss}$. This makes use of the method
of reduction to the maximal torus \cite{GK,martin, martin2}.

When $M^{ss} = M^s$ one can reduce to the maximal torus as
follows. Let $\mu_T:M \to \liets$ be the $T$-moment map given by
composing $\mu:M \to \lieks$ with the natural map $\lieks \to
\liets$. As Guillemin and Kalkman observe in \cite{GK}, it follows
immediately from the residue formula (\ref{res}) that if $0$ is a
regular value\footnote{We are assuming that $0$ is a regular value
of $\mu$, or equivalently that $K$ acts with finite stabilisers on
$\mu^{-1}(0)$.}  of $\mu_T$ then we have a surjection
$\k_M^{T}:H^*_T(M) \to H^*(M/\!/T_c) =H^*(\mu_T^{-1}(0)/T)$
defined as at (\ref{1.1}), and if $\alpha, \beta \in \hk(M)$ then
\beq \label{rtmt} \k_M(\alpha \beta) [\xg ] = \frac{1} {|W|}
\k_M^{T}(\alpha \beta \prod_{\g\in \G} \g) [\mu_T^{-1}(0)/T].\eeq
This formula requires some interpretation since  $\alpha, \beta\in
\hk(M)$ and $\prod_{\g \in \Gamma} \g \in H^*_T=H^*(BT)$, which we
think of as the equivariant cohomology of a point. We regard
$\alpha$ and $\beta$  as elements of  $H^*_T(M)$ via the natural
identification of $H^*_K(M)$ with the Weyl invariant part
$(H^*_T(M))^W$ of $H^*_T(M)$ and $\gamma$  as an element of
$H^*_T(M)$ via the natural inclusion of the equivariant cohomology
of a point in $H^*_T(M)$. In fact Martin \cite{martin, martin2}
gives a direct proof of (\ref{rtmt}) without appealing to the
residue formula. His proof shows also that, provided
$\mu^{-1}(0)/T$ is oriented appropriately, \beq \label{rtmt2}
\k_M(\alpha \beta) [\xg ] = \frac{1} {|W|}\k_M^{T} (\alpha \beta
\prod_{\g\in \G_+} \g) [\mu^{-1}(0)/T]\eeq where the product is
now over only the positive roots of $K$, and his argument shows in
addition that it is possible to represent the cohomology classes
$\k_M(\prod_{\g \in \G_+} \g)$ and $\k_M(\prod_{\g \in \G_-} \g)$
(which of course only differ by a sign $(-1)^{n_+}$) by closed
differential forms on $\mu_T^{-1}(0)/T$ with support in an
arbitrarily small neighbourhood of $\mu^{-1}(0)/T$. Thus there is
in fact no need to assume in (\ref{rtmt}) and (\ref{rtmt2}) that
$0$ is a regular value of $\mu_T$; it is enough to have $0$ a
regular value of $\mu$, and $M$ itself may have singularities away
from $\zloc$. This is important in \cite{JK2} when these ideas are
applied to the moduli spaces $\mnd$ when $n$ and $d$ are coprime,
and it will be similarly important in this paper.

If $M^{ss} \neq M^s$ then we can  apply (\ref{rtmt}) to the
blow-up $\tilde{M}$ of $M$ to get \beq \label{tor}
\k_{\tilde{M}}(\alpha \beta )[\txg]=\frac{
 (-1)^{n_+}}{
 |W|}
\k_{\tilde{M}}^T(\alpha \beta \cald^2)[\tilde{M}/\!/T_c]. \eeq It
is shown in \cite{jkkw} $\S$8 that it is possible to choose a
value $\xi \in \liets$  which is regular for both $\mu_T$ and
$\tilde{\mu}_T$, such that the difference between
$\k_{\tilde{M}}^T(\alpha \beta\cald^2)[\tilde{M}/\!/T_c]$ and the
evaluation on the fundamental class of
$\tilde{\mu}_T^{-1}(\xi)/T=\tilde{M}/\!/_{\xi} T_c$ of the
 cohomology class  induced by $\alpha \beta \cald^2 \in \hT(M)$
can be calculated in terms of data determined purely by the
construction of $\tilde{M}^{ss}$ from $M^{ss}$, which is canonical
and explicit, instead of the construction of $\tilde{M}$ from $M$,
and moreover this evaluation on $[\tilde{M}/\!/_{\xi} T_c]$ equals
the evaluation on the fundamental class of  $\mu_T^{-1}(\xi)/T =
M/\!/_\xi T_c$ of the cohomology class  induced by $\alpha \beta
\cald^2 $, which can be calculated by using the residue formula
(\ref{res}) applied to the action of $T$ on $M$ with the moment
map $\mu_T - \xi$. Combining all these calculations
 enables us to
calculate pairings in the cohomology of the partial
desingularisation $\xtg$ of $M/\!/G$.

Once we have reduced to calculating pairings on
$[\tilde{M}/\!/_{\xi} T_c]$, an alternative strategy to the use of
the residue formula (\ref{res}) is to
 follow the approach taken by Guillemin and
Kalkman in \cite{GK} and Martin in \cite{martin, martin2}, which
was applied to the moduli spaces $\mnd$ when $n$ and $d$ are
coprime in \cite{JK2}. This is to consider the change in
$$\k_M^{(\zeta)}(\alpha \beta) [\mu_T^{-1}(\zeta)/T],$$
for fixed $\alpha, \beta \in H^*_T(M)$, as $\zeta$ varies through
the regular values of $\mu_T$. This is sufficient, if $M$ is a
compact symplectic manifold, because the image $\mu_T(M)$ is
bounded, so if $\zeta$ is far enough from $0$ then
$\mu_T^{-1}(\zeta)/T$ is empty and thus $\k_M^{(\zeta)}(\alpha
\beta) [\mu_T^{-1}(\zeta)/T]=0.$ Now
 the image $\mu_T(M)$ is a convex
polytope \cite{Aconv,GSconv}; it is the convex hull in $\liets$ of
the set
$$\{ \mu_T(F) : F\in\calf\}$$
of the images $\mu_T(F)$ (each a single point of $\liets$) of the
connected components $F$ of the fixed point set $M^T$. This convex
polytope is divided by codimension-one walls into subpolytopes,
themselves convex hulls of subsets of $\{ \mu_T(F) : F\in\calf\}$,
whose interiors consist entirely of regular values of $\mu_T$.
When $\zeta$ varies in the interior of one of these subpolytopes
there is no change in $\k_M^{(\zeta)}(\alpha \beta)
[\mu_T^{-1}(\zeta)/T],$ so it suffices to understand what happens
as $\zeta$ crosses a codimension-one wall.

Any such wall is the image $\mu_T(M_1)$ of a connected component
$M_1$ of the fixed point set of a circle subgroup $T_1$ of $T$.
The quotient group $T/T_1$ acts on $M_1$, which is a symplectic
submanifold of $M$, and the restriction of the moment map $\mu_T$
to $M_1$ has an orthogonal decomposition
$$\mu_T|_{M_1} = \mu_{T/T_1} \oplus \mu_{T_1}$$
where $\mu_{T/T_1}: M_1 \to (\liet/\liet_1)^*$ is a moment map for
the action of $T/T_1$ on $M_1$ and $\mu_{T_1}:M_1 \to \liets_1$ is
constant (because $T_1$ acts trivially on $M_1$). If $\zeta_1$ is
a regular value of $\mu_{T/T_1}$ then there is a symplectic
quotient
$$ \mu_{T/T_1}^{-1}(\zeta_1 )/ (T/T_1),$$
and it is shown in \cite{GK} that the change in
$\k_M^{(\zeta)}(\alpha \beta) [\mu_T^{-1}(\zeta)/T]$ as $\zeta$
crosses the wall $\mu_T(M_1)$ is
$$({\rm res}_{M_1} (\alpha \beta))_{\zeta_1} [\mu_{T/T_1}^{-1}(\zeta_1 )/ (T/T_1) ]$$
for a suitable residue operation
$${\rm res}_{M_1} : H^*_T(M) \to H^{*-d_1}_{T/T_1}(M_1)$$
where $d_1 = {\rm codim} M_1 -2$. When $T$ is itself a circle,
this residue operation is given by restricting to $M_1$, dividing
by the equivariant Euler class of the normal bundle to $M_1$ in
$M$, and taking the ordinary residue  $\res_{\xvec=0}$ at $0$ on
$\CC$. This gives an inductive method for calculating the change
in $\k_M^{(\zeta)}(\alpha \beta) [\mu_T^{-1}(\zeta)/T]$ as the
wall is crossed, in terms of data on $M$ localised near $M^T$; it
is essentially equivalent to the residue formula (\ref{res}) when
$\dim(T)=1$, but differs from it for groups of higher rank.

\begin{rem} \label{endsection2} The advantage of this method over
the residue formula (\ref{res}) for our purposes is that it can be
applied in situations when $M$ is not compact, and indeed when the
fixed point set of the action of $T$ on $M$ has infinitely many
components so that the residue formula cannot be applied directly.
 This method was used in \cite{JK2} with $M$ as the extended moduli space
of \cite{J} to obtain formulas for the pairings on $\mnd$ when $n$
and $d$ are coprime (see $\S$4 below), and exactly the same
arguments will provide us with formulas for pairings on
${M}/\!/_{\xi} T_c$ when $\xi$ is a regular value of $\mu$
sufficiently close to $0$.
\end{rem}

\begin{rem}
When it is unlikely to cause confusion we will simplify the
notation a little and write $\kappa$ instead of $\kappa_M$,
$\kappa_M^{ss}$ etc.
\end{rem}

\renorm
\section{Moduli spaces of bundles and their partial desingularisations}

Recall that a holomorphic vector bundle $E$ of rank $n$ and degree
$d$ over the compact Riemann surface $\Sigma$ of genus $g \geq 2$
is called semistable (respectively stable) if every proper
subbundle $E'$ of $E$ satisfies ${\rm deg}(E')/{\rm rank}(E') \leq
d/n$ (respectively ${\rm deg}(E')/{\rm rank}(E') < d/n$).  There
is a moduli space $\calm^s(n,d)$ of isomorphism classes of stable
bundles of rank $n$ and degree $d$ over $\Sigma$, which is a
nonsingular quasi-projective variety with a natural
compactification $\mnd$ whose points are represented by semistable
bundles of rank $n$ and degree $d$ over $\Sigma$. The compactified
moduli space $\mnd$ is a projective variety which is in general
singular, although if $d$ and $n$ have no common factors then
$\mnd$ coincides with $ \calm^s(n,d)$ and is a nonsingular
projective variety.

The spaces $\mnd$ can be represented in several different ways as
quotients in the sense of Mumford's geometric invariant theory
\cite{MFK} or as  quotients in an analogous sense for infinite
dimensional group actions, leading to  constructions of  partial
desingularisations $\tilde{{\cal M}}(n,d)$ of $\mnd$ \cite{K5},
which are projective varieties with only orbifold singularities.
In this section we shall follow the infinite-dimensional point of
view taken by Atiyah and Bott in \cite{AB}.

Let ${\cal E}$ be a fixed $C^{\infty}$ complex hermitian vector
bundle of rank $n$ and degree $d$ over $\Sigma$. Let ${\cal C}$ be
the space of all holomorphic structures on ${\cal E}$, let ${\cal
C}^s$ (respectively ${\cal C}^{ss}$) be the open subset of ${\cal
C}$ consisting of all stable (respectively semistable) holomorphic
structures on ${\cal E}$, let ${\cal G}$ denote the gauge group of
${\cal E}$ (the group of all $C^{\infty}$ unitary automorphisms of
${\cal E}$) and let ${\cal G}_c$ denote its complexification which
is the group of all $C^{\infty}$ complex automorphisms of ${\cal
E}$. When it is necessary for clarification we shall write ${\cal
C}(n,d)$ and ${\cal G}(n,d)$ instead of ${\cal C}$ and ${\cal G}$.
 The
moduli space ${\cal M}^s(n,d)$ can be identified naturally with
${\cal C}^s/{\cal G}_c$ and $\mnd$ can be identified naturally
with the quotient of ${\cal C}^{ss}$ by the equivalence relation
for which semistable structures are equivalent if and only if the
closures of their ${\cal G}_c$-orbits meet in ${\cal C}^{ss}$.
Thus we can think of $\mnd$ as a quotient ${\cal C}/\!/{\cal G}_c$
in a sense analogous to geometric invariant theoretic quotients,
and so construct a partial desingularisation $\tilde{{\cal
M}}(n,d)$ of $\mnd$ \cite{K5}. In fact in \cite{K5} $\tilde{{\cal
M}}(n,d)$ is not constructed using the representation of $\mnd$ as
the geometric invariant theoretic quotient of ${\cal C}$ by ${\cal
G}_c$, although it is noted at \cite[p.246]{K5} that this
representation of $\mnd$ would lead to the same partial
desingularisation. Instead in \cite{K5} $\mnd$ is represented as a
geometric invariant theoretic quotient of a finite-dimensional
nonsingular quasi-projective variety $R(n,d)$ by a linear action
of $SL(p;\CC)$ where $p=d+n(1-g)$ with $d \gg 0$.

Alternatively $\mnd$ can be thought of as a symplectic quotient of
${\cal C}$ by the gauge group ${\cal G}$, where the r\^{o}le of
the normsquare of the moment map is played by the Yang-Mills
functional.

The construction of $\tilde{{\cal M}}(n,d)$ involves a set ${\cal
R}$ of representatives $R$ of the conjugacy classes of reductive
subgroups of ${\cal G}_c$ which occur as the connected components
of stabilisers in ${\cal G}_c$ of semistable points of
$\mathcal{C}$, together with their fixed point sets $Z_R$ in
$\mathcal{C}$. Equivalently we look for automorphism groups of
semistable bundles over $\Sigma$. Such conjugacy classes
correspond to unordered sequences $(m_1,n_1),...,(m_q,n_q)$ of
pairs of positive integers such that $m_1 n_1 + ... + m_q n_q = n$
and $n$ divides $n_i d$ for each $i$ (cf. \cite[pp. 248-9]{K5}).
An element $R$ of the corresponding conjugacy class is given by
\begin{equation} \label{labelR} R= GL(m_1;\CC) \times ... \times GL(m_q;\CC). \end{equation}
In the notation of \cite{K4} and \cite{K5} we construct
$\tilde{{\cal C}}^{ss}$ from ${\cal C}^{ss}$ by blowing up along
the subvarieties ${\cal G}_c Z_R^{ss}$ (or rather their proper
transforms),  in decreasing order of $\dim R$, and removing the
points which are not semistable at each stage, where $Z_R^{ss}$ is
the set of semistable holomorphic structures fixed by $R$.
Equivalently we construct $\tilde{{\cal M}}(n,d)$ from $\mnd$ by
blowing up along the images $Z_R/\!/N^R$ in $\mnd$ of the
subvarieties ${\cal G}_c Z_R^{ss}$, where $N=N^R$ is the
normaliser of $R$ in ${\cal G}_c$. A holomorphic structure fixed
by $R$ is semistable (equivalently belongs to $Z_R^{ss}$) if and
only if it is semistable for the induced action of $N/R$ on $Z_R$;
we denote by $Z_R^s$ the open subset of $Z_R^{ss}$ consisting of
those holomorphic structures fixed by $R$ which are stable for the
induced action of $N/R$ on $Z_R$. When $R$ is as at (\ref{labelR})
above then ${\cal G}_c Z_R^{ss}$ consists of all those holomorphic
structures $E$ with
\begin{equation}
\label{newo*} E \cong (\CC^{m_1} \otimes D_1) \oplus \cdots \oplus
(\CC^{m_q} \otimes D_q)
\end{equation}
with $D_1, \ldots ,D_q$ all semistable and $D_i$ of rank $n_i$ and
degree $d_i = n_i d/n$, while ${\cal G}_c Z_R^{s}$ consists of all
those holomorphic structures $E$ as above where $D_1, \ldots ,D_q$
are all stable and not isomorphic to one another.
 Moreover the normaliser $N$ of $R$ in
${\cal G}_c$ has connected component
\begin{equation} \label{N_0} N_0 \cong \prod_{1 \leq i \leq q} (GL(m_i;\CC) \times {\cal G}_c(n_i,d_i))/\CC^*
\end{equation}
where $\CC^*$ is the diagonal central one-parameter subgroup of
$GL(m_i;\CC) \times {\cal G}_c(n_i,d_i)$. The group $\pi_0(N) =
N/N_0$ is the product
\begin{equation} \label{N/N_0} \pi_0(N) = \prod_{j\geq 0, k \geq 0}
\Sym({\#\{i:m_i=j \mbox{ and }n_i=k\}})
\end{equation}
where $\Sym(b)$ denotes the symmetric group of permutations of a
set with $b$ elements. Furthermore the semistable holomorphic
structures which become unstable after the blow-up corresponding
to the conjugacy class of $R$ are those with a filtration $0=E_0
\subset E_1 \subset ... \subset E_s= E$ such that $E$ is not
isomorphic to $\bigoplus_{1 \leq k \leq s} E_k/E_{k-1} $ but
$$\bigoplus_{1 \leq k \leq s} E_k/E_{k-1}
\cong (\CC^{m_1} \otimes D_1) \oplus \cdots \oplus (\CC^{m_q}
\otimes D_q)
$$
where $D_1, \ldots ,D_q$ are all stable and not isomorphic to one
another, and $D_i$ has rank $n_i$ and degree $d_i$ \cite[p.
248]{K5}.

In \cite{K8} the action of $R$ on the normal ${\cal N}_{R}$ to
${\cal G}_cZ^{ss}_{R}$ at a point represented by a holomorphic
structure $E$ of the form (\ref{newo*}) and the induced action on
$\PP({\cal N}_{R})$ are studied. If a $C^{\infty}$ isomorphism of
our fixed $C^{\infty}$ bundle ${\cal E}$ with $ (\CC^{m_1} \otimes
D_1) \oplus \cdots \oplus (\CC^{m_q} \otimes D_q)$ is chosen, then
we can identify ${\cal C}$ with the infinite-dimensional vector
space
$$\Omega^{0,1}(\mbox{\End}( (\CC^{m_1} \otimes D_1) \oplus \cdots \oplus (\CC^{m_q} \otimes D_q)))$$
and the normal to the ${\cal G}_c$-orbit at $E$ is given by
$H^1(\Sigma, \mbox{End} E)$, where $\mbox{End} E$ is the bundle of
holomorphic endomorphisms of $E$ \cite[$\S$7]{AB}. The normal to
${\cal G}_cZ^{ss}_{R}$ can then be identified with
\[
H^1(\Sigma,\mbox{End}'_{\oplus} E) \cong \bigoplus_{i_1,i_2=1}^{q}
\CC^{m_{i_1}m_{i_2}-\delta_{i_1}^{ i_2}} \otimes
H^1(\Sigma,D^*_{i_1} \otimes D_{i_2})
\]
where $\delta_i^j$ denotes the Kronecker delta and
$\mbox{End}'_{\oplus} E$ is the quotient of the bundle
$\mbox{End}E$ of holomorphic endomorphisms of $E$ by the subbundle
$\mbox{End}_{\oplus} E$ consisting of those endomorphisms which
preserve the decomposition (\ref{newo*}). The action of $ R =
\prod_{i=1}^{q} GL(m_i;\CC) $ on this is given by the natural
action on $\CC^{m_{i_1} m_{i_2} - \delta_{i_1}^{ i_2}}$ identified
with the set of $m_{i_1} \times m_{i_2}$ matrices if $i_1 \neq
i_2$ and the set of trace-free matrices if $i_1 = i_2$; its
weights are of the form $ \xi - \xi'$ where $\xi$ and  $\xi'$ are
weights of the standard representation of $R$ on $\oplus_{i=1}^{q}
\CC^{m_i}$.

The linear action of $R$ on the normal ${\cal N}_R$ induces a
stratification of $\PP({\cal N}_R)$ with the semistable set
$\PP({\cal N}_R)^{ss}$ as its open stratum \cite{K2}. An element
$\beta$ of the indexing set ${\cal B}_R$ of this stratification
 is represented by the closest point to $0$ of the convex hull of some nonempty set of the weights
of the action of $R$ on ${\cal N}_R$, and two such closest points
can be taken to represent the same element of ${\cal B}_R$ if and
only if they lie in the same $Ad(N)$-orbit, where $N$ is the
normaliser of $R$ (see \cite{K2} or \cite{K7}). This indexing set
${\cal B}_R$ is described more explicitly in \cite{K8} as follows.
Let us take our maximal compact torus $T_R$ in $R$ to be the
product of the standard maximal tori of the unitary groups
$U(m_1)$,..., $U(m_q)$ consisting of the diagonal matrices, and
let $\liet_R$ be its Lie algebra. Let
$$M = m_1 + ... + m_q$$
and let $e_1,...,e_M$ be the weights of the standard
representation of $T_R$ on $\CC^{m_1}\oplus ... \oplus \CC^{m_q}$.
We use the usual invariant inner product on the Lie algebra ${\bf
u}(m_i)$ of $U(m_i)$ for $1 \leq i \leq q$ given by $\langle A,
B\rangle = -{\rm tr} A\bar{B}^t$ and multiply by a positive scalar
factor $\rho_i$ (to be chosen later) to induce an inner product on
the Lie algebra of $T_R$ such that $e_1,...,e_M$ are mutually
orthogonal and $|\!|e_j|\!|^2 = \rho_i$ if $m_1 +...+ m_{i-1} < j
\leq m_1 +...+m_i$. Note that in \cite{K8} $\rho_i$ is chosen to
be $(n_i + d_i(1-g))^{-1}$, but this does not affect the proof of
the following result which is \cite{K8} Proposition 5.1.

\begin{prop} \label{5.1} Let $\beta$ be any nonzero element of the Lie algebra
$\liet_R$ of the maximal compact torus $T_R$ of $R$. Then $\beta$
represents an element of ${\cal B}_R \backslash \{0 \}$ if and
only if there is a partition
$$\{ \Delta_{h,m}:(h,m) \in J \}$$
of $\{1,...,M \}$, indexed by a rectangle $J$ in $\ZZ \times \ZZ$,
with the following properties. If
$$r_{h,m} = \sum_{j \in \Delta_{h,m}} |\!|e_j|\!|^{-2}$$
and
$$\epsilon(h) = \left( \sum_{m} m r_{h,m}\right)
\left(\sum_{m} r_{h,m}\right)^{-1},$$ then $-1/2 \leq \epsilon(h)
< 1/2$ and $\epsilon(1) > \epsilon(2) >...$, and
$$\frac{\beta}{|\!| \beta |\!|^2} = \sum_{(h,m) \in J}
\sum_{j \in \Delta_{h,m}} (\epsilon(h) - m) \frac{e_j}{|\!| e_j
|\!|^2}.$$ Moreover the conditions on the function $\epsilon$
ensure that the partition $\{\Delta_{h,m}:(h,m) \in J \}$ and its
indexing can be recovered from the coefficients of $\beta$ with
respect to the basis
$$e_1/|\!|e_1 |\!|^2, ..., e_M/|\!|e_M|\!|^2$$
of $\liet_R$.
\end{prop}

\newcommand{\lij}{\lambda_{ij}}
\newcommand{\lji}{\lambda_{ji}}
\newcommand{\lik}{\lambda_{ik}}
\newcommand{\ljk}{\lambda_{jk}}
\newcommand{\lki}{\lambda_{ki}}
\newcommand{\lkj}{\lambda_{kj}}
\newcommand{\lijbh}{\lambda_{ij}^{\beta h}}
\newcommand{\lijb}{\lambda_{ij}^{\beta}}

The proof of this proposition involves studying the convex hull of
$$\{ e_i - e_j: (i,j) \in S \}$$
for a nonempty subset $S$ of $\{ (i,j) \in \ZZ \times \ZZ: 1 \leq
i,j \leq M \}$. From $S$ we can construct a directed graph $G(S)$
with vertices $1,...,M$ and directed edges from $i$ to $j$
whenever $(i,j)\in S$. Let $\Delta_1$,..., $\Delta_s$ be the
connected components of this graph. Then $\{ e_i - e_j: (i,j) \in
S \}$ is the disjoint union of its subsets $\{ e_i - e_j: (i,j)
\in S \ \textrm{and}\  i,j \in \Delta_h
 \}$ for $1 \leq h \leq t$, and $\{ e_i - e_j: (i,j) \in S\ \textrm{and}\   i,j \in \Delta_h
 \}$ is contained in the vector subspace of $\liet_R$ spanned by the basis
vectors $\{e_k: k \in \Delta_h\}$. Since these subspaces are
mutually orthogonal for $1 \leq h \leq s$, the closest point to
$0$ in the convex hull of $\{e_i - e_j:(i,j) \in S \}$ is
\begin{equation} \label{beta} \beta =\left(   \sum_{h=1}^s \frac{1}{|\!| \beta_h |\!|^2}
\right )^{-1} \sum_{h=1}^s \frac{\beta_h}{|\!| \beta_h |\!|^2}
\end{equation}
where $\beta_h$ is the closest point to $0$ of the convex hull of
$\{ e_i - e_j: (i,j) \in S\ \textrm{and}\   i,j \in \Delta_h  \}$
for $1 \leq h \leq s$. It is shown in the proof of \cite{K8}
Proposition 5.1 that we can express each $\Delta_h$ as a disjoint
union
$$\Delta_h = \bigcup_m {\Delta}_{h,m} $$
such that
$$\frac{\beta_h}{|\!| \beta_h |\!|^2} = \sum_{m}
\sum_{j \in \Delta_{h,m}} (\epsilon(h) - m) \frac{e_j}{|\!| e_j
|\!|^2}$$ where $\epsilon(h)$ has the required properties, and
from this the result follows.

Recall from Remark \ref{endsection2} that in order to compute
intersection pairings on $\tilde{{\cal M}}(n,d)$ we will consider
the change in pairings as walls are crossed between convex
subpolytopes whose interiors consist of regular values of a torus
moment map. To deal with the blow-up corresponding to the subgroup
$$ R= GL(m_1;\CC) \times ... \times GL(m_q;\CC) $$
as above at (\ref{labelR}), we will fix a ray
\begin{equation} \label{ray} \RR_+(e_1 + e_2 + \cdots + e_{M-1} - (M-1)e_M) \end{equation}
in $\liet_R$ and consider the wall crossings needed to approach
$0$ along this ray. Here the walls are convex hulls of subsets of
the set of weights $\{e_i-e_j : 1 \leq i,j\leq M\}$ for the action
of $R$ on ${\cal N}_R$, and hence all lie within the codimension
$1$ subspace in $\liet_R$ given by $\{\sum_i \lambda_ie_i : \sum_i
\lambda_i =0\}$.

Any wall which needs to be crossed lies in a hyperplane in this
subspace obtained by intersecting the subspace with one of the
affine hyperplanes $\beta + \beta^{\perp}$ determined by some
$\beta$ representing an element of ${\cal B}_R \setminus \{0\}$.
Such a $\beta$ corresponds to a partition $\{\Delta_{h,m}:(h,m)
\in J \}$ satisfying the conditions in Proposition  \ref{5.1}, or
equivalently is the closest point to zero in the convex hull of a
nonempty subset $\{e_i-e_j : (i,j) \in S\}$ of weights of the $R$
action on ${\cal N}_R$. The subset $S$ determines a directed graph
$G(S)$ as above.

\begin{lemma}
\label{no loops} The directed graph $G(S)$ corresponding to a
non-zero $\beta$ contains no directed loops.
\end{lemma}
\noindent{\bf Proof:} Suppose we have a directed loop in $G(S)$,
that is a sequence of edges corresponding to weights
$e_{i_1}-e_{i_2}, e_{i_2}-e_{i_3},\ldots, e_{i_r}-e_{i_1}$. Then
$\frac{1}{r} \left((e_{i_1}-e_{i_2})+\ldots
+(e_{i_r}-e_{i_1})\right)=0$ lies in the convex hull of the
weights. This contradicts the assumption that $\beta$ is the
closest point to the origin of this convex hull.\hfill $\Box$

\bigskip

\begin{lemma}
\label{lemma2} If the ray $ \RR_+(e_1 + e_2 + \cdots + e_{M-1} -
(M-1)e_M)$ meets the wall determined by $\beta$ then $M$ is the
only vertex of the directed graph $G(S)$ with no outgoing edges.
\end{lemma}

\noindent{\bf Proof:} Since the ray (\ref{ray}) meets the wall
determined by $\beta$ we have
$$\lambda (e_1 + e_2 + \cdots + e_{M-1} - (M-1)e_M) = \sum_{(i,j) \in S} \lambda_{ij} (e_i - e_j)$$
for some $\lambda >0$ and $\lambda_{ij} \geq 0$ with $\sum_{(i,j)
\in S} \lij = 1$. Equating coefficients of $e_i$ for $i \leq M-1$
gives
$$\lambda = \sum_{j:(i,j) \in S} \lij - \sum_{j:(j,i) \in S} \lji$$
and hence
$$\sum_{j:(i,j) \in S} \lij = \lambda + \sum_{j:(j,i) \in S} \lji >0$$
so there is some $j$ for which $(i,j) \in S$ and hence $i$ has an
outgoing edge.

Suppose $M$ has an outgoing edge. Then all the vertices would have
outgoing edges and $G(S)$ would contain a directed loop,
contradicting Lemma \ref{no loops}.\hfill $\Box$

\bigskip

\begin{lemma}
\label{graph connected} If the ray $ \RR_+(e_1 + e_2 + \cdots +
e_{M-1} - (M-1)e_M)$ meets the wall determined by $\beta$ then
the directed graph $G(S)$ is connected.
\end{lemma}
\noindent{\bf Proof:} If $G(S)$ is not connected there is a
component not containing $M$. Every vertex in this component has
an outgoing edge by Lemma \ref{lemma2} and so the component
contains a directed loop contradicting Lemma \ref{no loops}.\hfill
$\Box$

\bigskip

\begin{rem} \label{rem1} Since $G(S)$ is connected the index $h$ for
the partition $\{\Delta_{h,m}:(h,m) \in J \}$ can be omitted. We
will also relabel the partition $\{\Delta_{m} \}$ by adding a
constant to $m$ so that it is indexed by $m \in \{1, \ldots, t\}$;
the only difference this makes is that we can no longer assume
that $\epsilon$ lies in $ [-1/2, 1/2)$. We also know from Lemma
\ref{lemma2} that $M$ is the only \lq top' element of the graph.
It follows from the definition of the $\Delta_{h,m}$ in the proof
of \cite{K8} Proposition 5.1 that $M$ is then the only element in
$\Delta_t$.

\end{rem}

Subpolytopes of the wall determined by $\beta$ are the
intersections of convex hulls of subsets of $\{ e_i - e_j:(i,j)
\in S\}$ with exactly $M-1$ elements which are linearly
independent. These correspond to subgraphs of $G(S)$ with
precisely $M-1$ edges. If the ray (\ref{ray}) meets the
subpolytope it must correspond to a connected subgraph with $M$ as
the only vertex with no outgoing edge (by the same arguments as
for Lemmas \ref{no loops}, \ref{lemma2} and \ref{graph
connected}). Thus these subgraphs are trees; they are the minimal
connected subgraphs of $G(S)$ with the same vertices as $G(S)$ and
having $M$ as the only vertex with no outgoing edges.

\begin{lemma} \label{lem3} If the connected graph $G(S)$ has the property that its only vertex with no
outgoing edges is $M$, then $G(S)$ has a minimal connected
subgraph with the same property.
\end{lemma}

\noindent{\bf Proof:} If $1 \leq i \leq M-1$ then we can pick some
$j(i) \in \{1, \ldots, M\}$ such that $(i,j(i)) \in S$. Since
$G(S)$ contains no directed loops we can check that $\{(i,j(i)):1
\leq i \leq M-1 \}$ determines a minimal connected subgraph with
the required property.\hfill $\Box$
\bigskip

\begin{lemma} \label{lem4}
Suppose that $G(S)$ has the property that its only vertex with no
outgoing edges is $M$, and let $G_0=G(S_0)$ be a minimal connected
subgraph with the same property. Then the ray $ \RR_+(e_1 + e_2 +
\cdots + e_{M-1} - (M-1)e_M)$ meets the interior of the convex
hull of $\{ e_i - e_j:(i,j) \in S_0\}$.
\end{lemma}

\noindent{\bf Proof:} Without loss of generality we can assume
that the vertices $i$ such that $(i,M) \in S_0$ are precisely
$M-1,M-2, \ldots,M-k.$ If we remove the vertex $M$ and the edges
joining $M-1,M-2, \ldots,M-k$ to $M$ from the graph $G(S_0)$, then
the resulting graph has $k$ connected components $G_1 = G(S_1),
\ldots, G_k = G(S_k)$, say, where for $1 \leq j_0 \leq k$ the
vertex $M - {j_0}$ is the only vertex in $G_{j_0}$ with no
outgoing edges. Let  $M_{j_0}$ be the number of vertices in the
connected component $G_{j_0}$ and let the other vertices of
$G_{j_0}$ apart from $M-j_0$ be $i_1^{j_0}, \ldots,
i^{j_0}_{M_{j_0}}$. Then $G_{j_0}$ has $M_{j_0} - 1$ edges and is
a minimal connected graph on its vertices with the property that
$M-j_0$ is the only vertex with no outgoing edges. Hence by
induction on $M$ we can assume that there exist
$\lambda^{(j_0)}>0$ and $\lij^{(j_0)}>0$ for $(i,j) \in S_{j_0}$
with
$$\sum_{(i,j) \in S_{j_0}} \lij^{(j_0)} = 1 $$
and
$$\lambda^{(j_0)} \left(e_{i^{j_0}_1} + e_{i^{j_0}_2} + \cdots +e_{i^{j_0}_{M_{j_0}-1}} -
(M_{j_0} -1)e_{M-j_0} \right) = \sum_{(i,j) \in S_{j_0}}
\lij^{(j_0)} (e_i - e_j).$$ Dividing by $\lambda^{(j_0)}$ (which
is strictly positive) gives
$$e_{i^{j_0}_1} + e_{i^{j_0}_2} + \cdots +e_{i^{j_0}_{M_{j_0}-1}} -
(M_{j_0} -1)e_{M-j_0}  = \sum_{(i,j) \in S_{j_0}}
\frac{\lij^{(j_0)}}{\lambda^{(j_0)}} (e_i - e_j)$$ and summing
over $j_0=1,\ldots,k$ gives
$$e_{1} + e_{2} + \cdots +e_{M-k-1} -
(M_{k} -1)e_{M_k} - \cdots - (M_1-1)e_{M-1}  = \sum_{j_0=1}^k
\sum_{(i,j) \in S_{j_0}} \frac{\lij^{(j_0)}}{\lambda^{(j_0)}} (e_i
- e_j).$$ Adding $M_k(e_{M-k} - e_M) + \cdots + M_1 (e_{M-1} -
e_M)$ to each side and using the equality $M_1 + \cdots + M_k =
M-1$ gives
$$e_{1} + e_{2} + \cdots +e_{M-1} -
 (M-1)e_{M}  = \sum_{j_0=1}^k \sum_{(i,j) \in S_{j_0}}
\frac{\lij^{(j_0)}}{\lambda^{(j_0)}} (e_i - e_j) + \sum_{j_0=1}^k
M_{j_0} (e_{M-j_0} - e_M).$$ Since $S_{j_0} \subseteq S_0$ and
$(M-j_0,M) \in S_0$ for $1 \leq j_0 \leq k$, and in addition
$\lij^{(j_0)} / \lambda^{(j_0)} >0$ and $M_{j_0} >0$ for $1 \leq
j_0 \leq k$ and $(i,j) \in S_{j_0}$, we can rewrite this as
$$e_{1} + e_{2} + \cdots +e_{M-1} -
 (M-1)e_{M}  =  \sum_{(i,j) \in S_{0}}
\lij' (e_i - e_j) $$ where $\lij' \geq 0$ for all $(i,j) \in S_0$.
Indeed since
$$S_0 = \{ (M-j_0,M): 1 \leq j_0 \leq k \} \cup \bigcup_{j_0=1}^k S_{j_0}$$
we have $\lij' >0$ for all $(i,j) \in S_0$. Finally dividing each
side by
$$\sum_{(i,j) \in S_0} \lij'$$
gives
$$\lambda(e_{1} + e_{2} + \cdots +e_{M-1} -
 (M-1)e_{M})  =  \sum_{(i,j) \in S_{0}}
\lij (e_i - e_j) $$ where $\lambda>0$ and $\lij >0$ for all $(i,j)
\in S_0$ and
$$\sum_{(i,j) \in S_0} \lij = 1$$
as required.  \hfill $\Box$

\bigskip

Combining Lemmas \ref{lemma2}, \ref{graph connected}, \ref{lem3}
and \ref{lem4} gives us

\begin{prop} \label{propn5} The ray $\RR_+(e_{1} + e_{2} + \cdots +e_{M-1} -
 (M-1)e_{M}) $ meets the wall determined by $\beta$ if and only if $M$ is
the only vertex of $G(S)$ with no outgoing edges.
\end{prop}

\begin{rem} We now know which minimal connected subgraphs of $G(S)$
determine subpolytopes of this wall through which the ray
$\RR_+(e_{1} + e_{2} + \cdots +e_{M-1} -
 (M-1)e_{M}) $ passes. They are the minimal connected subgraphs with the property
that $M$ is the only vertex with no outgoing edges. The case $M=4$
is illustrated in Figure 1.
\end{rem}

\medskip

\begin{figure}[h]
 \centerline {
\includegraphics[width=4in]{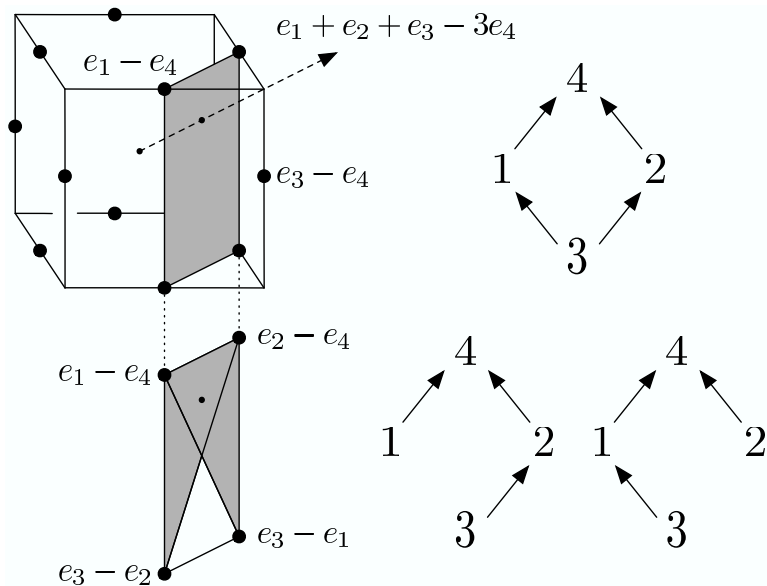}
} \caption{\label{figure} The case $M=4$. 
The $12$ weights $\{e_i-e_j : 1 \leq i,
j \leq 4, i\neq j\}$ are shown as dots forming the midpoints of a
hexahedron. A wall, corresponding to the subset
$S=\{(1,4),(2,4),(3,1),(3,2)\}$, which meets the ray
$\mathbb{R}^+(e_1+e_2+e_3-3e_4)$ is shown shaded. The connected
graph $G(S)$, in which $4$ is the only vertex with no outgoing
edges, is shown on the right. The subpolytopes of the wall which
meet the ray are shown below, and on their right are the
corresponding minimal connected subgraphs of $G(S)$, again in
which $4$ is the only vertex with no outgoing edges.}
\end{figure}

The procedure described in $\S$8 of \cite{jkkw} for calculating
intersection pairings works most efficiently if the wall crossing
takes place at a point $\beta^*$ of the affine hyperplane $\beta +
\beta^{\perp}$ in the Lie algebra of the maximal torus of $R$
which satisfies $\stab \beta \subseteq \stab \beta^* $ where
$\stab \beta$ denotes the stabiliser of $\beta$ under the adjoint
action of ${\cal G}_c$ (cf. \cite{jkkw} Remark 27). Recall that
$$\frac{\beta}{|\!| \beta |\!|^2} =  \sum_{m = 1}^{t}
\sum_{j \in \Delta_{m}} (\epsilon - m) \frac{e_j}{|\!| e_j
|\!|^2}$$ where $\epsilon$ is a constant. Recall also that we
chose an invariant inner product on the Lie algebra of $R$ such
that $e_1,\ldots,e_M$ are mutually orthogonal and $|\!|e_j|\!|^2 =
\rho_i$ if $m_1 +...+ m_{i-1} < j \leq m_1 +...+m_i$ where
$\rho_i$ can be any strictly positive scalar for $1 \leq i \leq
q$. Thus for generic choices of $\rho_1, \ldots, \rho_q$ we will
have the required condition
\begin{equation} \label{stabsubset} \stab \beta \subseteq \stab \beta^* \end{equation}
provided that $\beta^*$ is of the form
\begin{equation} \label{betastarform} \beta^* =   \sum_{m = 1}^{t}
\sum_{j \in \Delta_{m}} f(m) \frac{e_j}{|\!| e_j |\!|^2}
\end{equation} for some function $f(m)$ of $m \in \{ 1, \ldots, t
\}$. But by Proposition \ref{propn5} the condition for the ray
$\RR_+(e_{1} + e_{2} + \cdots +e_{M-1} -
 (M-1)e_{M}) $ to meet the wall determined by $\beta$ is that $M$ should be
the only vertex of $G(S)$ with no outgoing edges, and hence the
only element of $\Delta_{t}$ (see Remark \ref{rem1}). If this
condition is satisfied, then the ray meets the wall in a point
$\beta^*$ of the form (\ref{betastarform}) above, since
$$\lambda(e_{1} + e_{2} + \cdots +e_{M-1} -
 (M-1)e_{M})  = \sum_{m = 1}^{t}
\sum_{j \in \Delta_{m}} f(m) \frac{e_j}{|\!| e_j |\!|^2}$$ where
$f(m) = -\lambda (M-1) |\!| e_M|\!|^2$ if $m=t$ and $f(m) =
\lambda |\!| e_m|\!|^2$ if $1 \leq m < t$.

\begin{rem}
\label{generic rhos remark} With such generic choices of
$\rho_1,\ldots,\rho_q$ we get $R \cap \mathrm{Stab} \beta = \prod
GL(m_i^k)$ where $m_i^k = |\Delta_i^k|$ and
$$
\Delta_i^k = \Delta_k \cap \{m_1+\ldots+m_{i-1}+1, \ldots ,
m_1+\ldots+m_{i}\}
$$
cf. Definition 5.5 in \cite{K8}. We will make use of this in \S6.
\end{rem}

\renorm
\section{Intersection theory on nonsingular moduli spaces of bundles}

\newcommand{\diag}{{\rm diag}}
\newcommand{\inpr}[1]{{\langle #1 \rangle}}
\newcommand{\ktg}{{K^{2g}} }
\newcommand{\mext}{ { M_K^{\rm ext} } }

In this section we review the results of \cite{JK2} (see also
\cite{JS}). Throughout we will use a fixed invariant inner product
on the Lie algebra $\liek$ of a compact Lie group $K$ to identify
$\liek$ with $\liek^*$.

\subsection{Generators of the cohomology ring}
A set of generators for the cohomology\footnote{In this paper, all
cohomology groups are assumed to be with complex coefficients.}
  $H^*(\mnd)$
of the moduli space $\mnd$ of stable holomorphic vector bundles of
coprime rank $n$ and degree $d$  on a compact Riemann surface
$\Sigma$ of genus $g \ge 2 $
 is given in \cite{AB} by Atiyah and Bott.
It may be described as follows. There is a universal rank $n$
vector  bundle
$$ \UU \to \Sigma \times \mnd $$
which is unique up to tensor product with the pullback of any
holomorphic line bundle on $\mnd$; for definiteness Atiyah and
Bott impose an extra normalizing condition which determines the
universal bundle up to isomorphism, but this is not crucial to
their argument (see \cite{AB}, p582).
 Then by \cite{AB} Proposition 2.20 the following elements
of  $H^*(\mnd)$ for $2\leq r\leq n$ make up a set of generators:
\begin{equation} \label{e:fdef}
f_r =  c_r(\UU)/[\Sigma], ~~2 \le r \le n \end{equation}
 \begin{equation} \label{e:bdef}  b_r^j =  c_r(\UU)/\a_j, ~~1 \le r \le n \end{equation}
\begin{equation} \label{e:adef} a_r =  c_r(\UU)/ 1,~~2 \le r \le n \end{equation}
Here, $[\Sigma]$
 $ \in H_2(\Sigma)$ and $\a_j \in H_1(\Sigma)$ $(j = 1, \dots, 2g)$
form standard bases  of $H_2(\Sigma, \ZZ)$ and  $H_1(\Sigma,
\ZZ)$, and $/$  represents the slant product $H^N(\Sigma \times
\mnd) \otimes H_j(\Sigma) \to H^{N-j}(\mnd). $

\begin{rem}
\label{abf lift to equivariant coh} These generators lift to
equivariant cohomology: there is a $\mathcal{G}(n,d)$-equivariant
universal bundle over $\Sigma \times \mathcal{C}(n,d)$ and the
slant product of its equivariant Chern classes with $1 \in
H_0(\Sigma)$, $[\Sigma] \in H_2(\Sigma)$ and the $\a_j \in
H_1(\Sigma)$ give classes which (abusing notation) we denote
respectively $a_1,\ldots , a_r$,  $f_2,\ldots, f_r$ and $b^j_r$
for $1 < r \leq n$ and $1 \leq j \leq 2g$. These classes generate
$H_{\mathcal{G}(n,d)}^*(\mathcal{C}(n,d))$.

The quotient $\mathcal{G}(n,d) \to \overline{\mathcal{G}}(n,d)$ by
the central subgroup $S^1$ induces an inclusion
$$H_{\overline{\mathcal{G}}(n,d)}^*(\mathcal{C}(n,d)) \hookrightarrow H_{\mathcal{G}(n,d)}^*(\mathcal{C}(n,d))$$ such that identifying $H^*(BS^1)$ with the polynomial subalgebra generated by $a_1$ gives an isomorphism
$$
H_{\mathcal{G}(n,d)}^*(\mathcal{C}(n,d))\cong
H_{\overline{\mathcal{G}}(n,d)}^*(\mathcal{C}(n,d)) \otimes
H^*(BS^1).
$$
Hence $a_2,\ldots , a_r$, $ f_2,\ldots, f_r$ and $b^j_r$ for $1< r
\leq n$ and $1 \leq j \leq 2g$ determine generators of
$$H_{\overline{\mathcal{G}}(n,d)}^*(\mathcal{C}(n,d)).$$ When $n$
and $d$ are coprime the latter is isomorphic to $H^*(\mnd)$ and
these equivariant classes correspond to those defined in
(\ref{e:fdef}), (\ref{e:bdef}) and (\ref{e:adef}). When $n$ and
$d$ are not coprime we think of $a_2,\ldots , a_r$, $ f_2,\ldots,
f_r$ and $b^j_r$ for $1< r \leq n$ and $1 \leq j \leq 2g$ as
equivariant cohomology classes.
\end{rem}

If we replace $\mnd$ by the moduli space $\mlnd$ of stable
holomorphic vector bundles  of coprime rank $n$ and degree $d$ and
fixed determinant line bundle, the set of generators is given by
(\ref{e:fdef}) and (\ref{e:adef}), while (\ref{e:bdef}) is
replaced by

 \begin{equation}
\label{e:bdef'}  b_r^j =  c_r(\UU)/\a_j, ~~2 \le r \le n.
\end{equation} The paper \cite{JK2} treats intersection numbers in
the cohomology of $\mlnd$ rather than that of $\mnd$, but the two
sets of pairings  are very closely related since
$$
H^*(\mnd) \cong H^*(\mlnd) \otimes H^*({\rm Jac})
$$
where ${\rm Jac} \cong U(1)^{2g}$ is the Jacobian (see Remark
\ref{r:finitecover} below).

\subsection{Extended moduli space}

Let $K=SU(n)$. We recall the construction of the extended moduli space from \cite{J}. This is a finite-dimensional (but non-compact) Hamiltonian $K$-space $\mext$ with moment map $\mu$ such that the symplectic quotient $\mu^{-1}(0)/K$ is homeomorphic to the moduli space $\mlnd$.

First consider the space $K^{2g}$. We define the map
 $\Phi: \ktg \to K$ by
$$\Phi(h_1,  \dots, h_{2g}) = \prod_{j = 1}^g h_{2j-1} h_{2j}
h_{2j-1}^{-1} h_{2j}^{-1} $$  and $\omega \in \Omega^2(\ktg)^K$ is
the 2-form defined in \cite[\S 8]{J2} (8.28)-(8.30) (cf. \cite[\S
9.1]{AMM}, Theorem 9.1). The space $\ktg$ is a quasi-Hamiltonian
$K$-space in the terminology of \cite{AMM}. We recall the
definition:

\begin{definition} \label{d:qham} (\cite{AMM})
A manifold  $M$ with a $K$-action and a 2-form $\omega$ and a map
$\Phi: M \to K$
 is a quasi-Hamiltonian $K$-space if  it satisfies
the following three axioms.

   (1) The differential of $\omega$ is given by:

    \[ d\omega = - {\Phi}^* \chi_K \]
where $\chi_K = (\theta_K, [\theta_K, \theta_K])$ is the
differential form on $K$ which represents the generator of $H^3(K,
\ZZ)$, in terms of the left invariant Maurer-Cartan form $\theta_K
\in \Omega^1(K) \otimes {\bf k}. $

   (2) The map $\Phi$ satisfies

    \[ \iota (\nu _\xi)\omega = \frac{1}{2} \Phi^* (\theta_K + \bar{\theta}_K, \xi)\]
where $\theta_K$ is the left invariant Maurer-Cartan form and
$\bar{\theta}_K$ is the corresponding right invariant
Maurer-Cartan form. Here, for $\xi \in \bf{k}$ we denote by
$\nu_\xi$ the  vector field on $M$ arising from the  action of
$K$.

   (3) At each $x \in M$, the kernel of $\omega _x$  is given by
       \[ \ker \omega _x = \{ \nu _{\xi}(x)~| ~\xi \in \ker(Ad_{\Phi(x)} + 1) \} \]

\end{definition}

We can construct a corresponding Hamiltonian $K$-space $\mext$ as
follows. We  choose an element $c \in Z(K)$, by setting
$$ c  = e^{2 \pi i d/n} {\bf I}$$
where $ {\bf I} $ is the identity matrix. Let
$${\mext} = \{ (m, \Lambda) \in \ktg \times {\bf{k}}\ :\ \Phi(m) =
c \exp (\Lambda) \} $$ and let $\mu: \mext \to {\bf{k}}$ be
defined by
\[
    \mu(m,\Lambda) = -\Lambda.
\]
Then we get the following commutative diagram.
 \begin{equation} \label{e:commdiag}
   \begin{CD}
   \mext     @>{-\mu}>> {\bf{k}}  \\
   \pi_1@VVV             @VVV c\exp   \\
   \ktg         @>{\Phi}>>  K
   \end{CD}
\end{equation}

Since $d \exp^*\chi_K = 0$ (where $\chi_K \in \Omega^3(K)$
represents the generator of $H^3(K, \ZZ)$; see Definition
\ref{d:qham} above)
 we can find some $\sigma \in \Omega^2({\bf{k}})$ such that $ d \sigma =  \exp^*\chi_K $.
We see that $d(\pi_1^* \omega - \mu^* \sigma) =0$. The space
$\mext$ defined by (\ref{e:commdiag}) becomes a Hamiltonian
$K$-space with moment map $\mu$ and the invariant 2-form
\begin{equation} \label{e:twoform}
{\tilde{\omega}} = \pi_1^* \omega - \mu^* \sigma \in
\Omega^2({\mext})^K .
\end{equation}
The space $\mext$ is the extended moduli space defined in
\cite{J}.

There are classes $\tilde{a}_r \in H^{2r}_K(\mext) $ (for $r = 2,
\dots, n$) which pass to $a_r$ $ \in H^{2r}(\mlnd)$ under the
composition of the restriction map $H^*_K(\mext) \to H^*_K(
\mu^{-1}(0))$ and the isomorphism $H^*_K(\mu^{-1}(0)) \cong
H^*(\mlnd)$ (this composition is sometimes referred to as the
Kirwan map). Likewise (for $r = 2, \dots, n$ and $j = 1, \dots,
2g$) there are
 classes
$\tilde{b}_r^j \in H^{2r-1}_K(\mext )$ and $\tilde{f}_r \in
H^{2r-2}_K(\mext )$ which pass to $b_r^j$ and $f_r$ under the
Kirwan map (see \cite{JK2}). The classes $\tilde{a}_r $ and
$\tilde{b}_r^j$ are invariant under translation in ${\bf k}$.

\begin{rem} We can modify the infinite dimensional description
used in $\S$3 so it  applies to  the space $\mlnd$. It is also
possible to  treat the space $\mnd$ using the finite dimensional
methods of the present section, with $U(n)$ replacing $SU(n)$.
The extended moduli space $\mext$ may be constructed as in
\cite{J} by a partial reduction of the infinite-dimensional space
${\mathcal C} $ by the
 based gauge group.
Hence, for our purposes working with the infinite dimensional
description using the gauge group and the space of all complex
structures is equivalent to working with the  finite dimensional
description via the extended moduli space for both $\mnd$ and
$\mlnd$.
\end{rem}

 \subsection{Equivariant Poincar\'e Dual} \label{ss:eqpd}
Since we know that $\ktg \times {\bf{k}}$ is always smooth, we
will work with integration over $\ktg \times {\bf{k}}$, instead of
working with integration over its subset $\mext$.

We work with the Cartan model of equivariant cohomology, for which
if the space $Y$ is equipped with an action of $K$ \beq
\label{e:eq1}  H^*_K(Y) = H^*(\Omega^*_K(Y), d_K) \end{equation}
Here, \beq \label{e:eq2} \Omega^*_K(Y) = \Bigl ( \Omega^*(Y)
 \otimes S(\liek^*) \Bigr )^K \eeq
and the equivariant differential $d_K  $ is given by \beq
\label{e:eq3} (d_K \alpha)(\xi) = d(\alpha(\xi)) -
\iota_{\nu_\xi}\alpha \eeq for $\xi \in \liek$ and $\alpha \in
\Omega^*_K(Y)$ where $\nu_\xi$ is the vector field on $Y$
generated by $\xi$.

\medskip

 \begin{lemma}  \label{l3.1} (\cite{JK2} Corollary 5.6)
  Let $T$ be the maximal torus of $K=SU(n)$ acting on $K$ by conjugation. If $c \in T$ then we
  can find a $T$-equivariantly closed differential form
$\hat{\alpha} \in \Omega_T^*(K)$ on $K$ with
  support arbitrarily close to $c$ such that
$$
    \int_K \eta \hat{\alpha} = \eta |_c \in H_T^*
$$
  for all $T$-equivariantly closed differential forms $\eta \in \Omega_T^*(K)$.
 \end{lemma}

 \begin{proposition} \label{p3.2}  (\cite{JK2}, Proposition 5.7)
 Let $P : \ktg \times {\bf k} \rightarrow K$ be
defined by
 \begin{equation} \label{e:pdef}
P: (m,\Lambda) \mapsto \Phi(m) \exp(-\Lambda)
\end{equation}
so that ${\mext} = P^{-1}(c)$. Let $\alpha = P^* \hat{\alpha}. $
Then
 $$
   \int_{\ktg \times {\bf{k}}} \eta \alpha = \int_{{\mext}} \eta
 $$
 for all $T$-equivariantly closed differential forms $\eta \in \Omega_T^*(\ktg \times {\bf k}).$

\end{proposition}

From now on we shall use $\mu$ to denote the map
 $$
    \mu : \ktg \times {\bf k} \rightarrow {\bf k}, ~ \mu: (m, \Lambda) \mapsto -\Lambda
  $$
as well as its restriction to $\mext$. Then $\mext = P^{-1}(c)$
$\subset \ktg \times \liek$  and $\mext\quott K = (P^{-1}(c) \cap
\mu^{-1}(0))/K$.

  Let $V$ be a small neighbourhood of $c$ in $K$.
  In fact, if $V^{'}$ is any neighbourhood of $c$ in $K$ containing
the closure of ${ V}$ then
  $$
    \int_{{P^{-1}(V^{'})}} \eta \alpha = \int_\mext \eta \in H_T^*
  $$

 \begin{lemma} \label{l3.3}
  If $0$ is a regular value of  $\mu$,
then $(P^{-1}(V) \cap \mu^{-1}(0))/T$ is
  an orbifold.
\end{lemma}

\noindent \textbf{Proof:}   Our first observation is that near
$\Phi^{-1}(c)$, the manifold $\ktg$ is endowed with a  symplectic
structure, as there is a $K$-invariant neighbourhood $V \subset K$
containing $c$ such that the restriction of the closed 2-form
$\tilde{\omega}$ to $(\Phi\circ \pi_1)^{-1}(V)$ is nondegenerate
(where $\tilde{\omega}$ was defined at (\ref{e:twoform})). This is
true for the following reason. Consider the diagram
(\ref{e:commdiag}). The space $\ktg$ is a smooth manifold, so the
space $\mext$ is smooth whenever $d(c\exp)$ is surjective (a
condition satisfied on a neighbourhood $V\subset K$ of $c$). The
two-form $\tilde{\omega}$ defined by (\ref{e:twoform}) is closed
on $\mext$. The map $\mu$ satisfies the moment map condition $
d\mu_\xi  = \tilde{\omega} (\nu_\xi, \cdot) $ on $\mext$, where
(for $\xi \in \bf{k}$) we denote by $\nu_\xi$ the  vector field on
$\ktg$ arising from the  action of $K$ (see \cite{AMM} and
\cite{J}). Furthermore the 2-form $\tilde{\omega}$ descends under
symplectic reduction from $(\Phi \circ \pi_1)^{-1}(V)$
 to the standard symplectic form on
$\calm(n,d)$. It follows that $\tilde{\omega}$ is nondegenerate on
$(\Phi\circ \pi_1)^{-1}(V), $ which is an open neighbourhood of
$\mu^{-1}(0)$ in $\mext$.

We know that $c$ is a regular value for $P : \ktg \times\mathbf{k}
\to K$, and therefore we can choose the neighbourhood $ V$  so
that all points of $V$ are also regular values of $P$ (by standard
properties of  the rank of a differentiable map). Because $
\Phi^{-1}(V)$ is symplectic with moment map $\mu$ related to
$\Phi$ as in diagram (\ref{e:commdiag}), the action of $K$ has
finite stabilizers at all points of $ \Phi^{-1}(V)$. This implies
that $T$ also acts with finite stabilizers at all points of $
(\Phi\circ \pi_1)^{-1}(V)$. Hence $P^{-1}(V) \cap \mu^{-1}(0)/T$
is an orbifold. \hfill $\Box$

\begin{rem}
In fact in the case when $K = SU(n)$ and $c = e^{2 \pi i d/n} {\bf
I}$ generates $Z(K)$, we may choose $V$ small enough to guarantee
that
  $T/Z(G)$ acts freely on
$P^{-1}(V) \cap \mu^{-1}(0)$ (since the action of $Z(G)$ is
trivial), so the quotient $P^{-1}(V) \cap \mu^{-1}(0)/T$  is a
smooth manifold (see \cite{JK2}, Lemma 5.10).
\end{rem}

We extend  the definition of the  composition
 $$
    \kappa : H_T^*(P^{-1}(c)) \rightarrow H_T^*(P^{-1}(c) \cap \mu^{-1}(0))
           \cong H^*(P^{-1}(c) \cap \mu^{-1}(0)/T)
 $$
  to
 $$
    \kappa :  H_T^*(P^{-1}(V)) \rightarrow H_T^*(P^{-1}(V) \cap \mu^{-1}(0))
           \cong H^*(P^{-1}(V) \cap \mu^{-1}(0)/T)
  $$

By Proposition \ref{p3.2} if $\eta \in H_T^*(\ktg \times {\bf k})$
we have
$$
  \int_{P^{-1}(c) \cap \mu^{-1}(0)/T} \kappa(\eta)
         = \int_{P^{-1}(V) \cap \mu^{-1}(0)/T} \kappa(\eta \alpha)
 $$
where the class $\alpha$ is the equivariant {Poincar\'e  dual} of
$P^{-1}(c)$ in $P^{-1}(V)$. The quantity $\int_{P^{-1}(V) \cap
\mu^{-1}(0)/T} \kappa(\eta \alpha)$ is
 given by a
 formula involving iterated residues, as we will see below.

\subsection{Periodicity} \label{ss:periodicity}

 We define a one dimensional torus ${\hat{T}}_1 \cong S^1$ in $K$
generated by the element ${\hat{e}}_1= (1, -1,0,\dots,0) $ in the
Lie algebra of the standard maximal torus $T$. Then $\hat{T_1}$
 is identified with $S^1$ via
 $$
   e^{2 \pi i \theta}  \in S^1 \mapsto (\exp{\theta \hat{e}_1 })
\in {\hat{T}}_1.
 $$
 The one dimensional Lie algebra
${\hat{\bf{t}}}_1$ is spanned by ${\hat{e}}_1$. Its
 orthocomplement in ${\bf t}$ is denoted by  ${\bf t}_{n-1}$.
 Define $T_{n-1}$ to be the torus given by $\exp ( {\bf t}_{n-1} )$.
Explicitly  this is
 $$
   T_{n-1} = \{(t_1, t_1, t_3,...,t_{n-1},t_n) \in U(1)^n  : (t_1)^2 ( \prod_{j=3}^n t_j ) = 1 \}.
$$
 Then
$T_{n-1}$ is isomorphic to the maximal torus of $SU(n-1)$ (i.e.
$T_{n-1} \cong
  (S^1)^{n-2}$).
We have
 $$
   {\bf t}_{n-1} = \{(X_1,...,X_n) \in {\Bbb R}^n : X_1 = X_2, \ \sum_{j=1}^n X_j = 0\}.
 $$
The torus $T_{n-1}$ has the same Lie algebra as the torus
$T/\hat{T_1}$.

\newcommand{\xc}{{ M(c)} }

\begin{lemma} \label{l3} (\cite{JK2}, Lemma 6.1)
Let $W$ be the Weyl group of $K=SU(n)$ so that the order of $W$ is
$|W|=n!$ and let $c=\diag(e^{2\pi i d/n},\ldots,e^{2\pi i d/n})$
where $d$ is coprime to $n$. If $V$ is a sufficiently small
neighbourhood of $c$ in $K$ that the quotient $T/Z_n$ of $T$ by
the centre $Z_n$ of $K=SU(n)$ acts freely on $P^{-1}(V) \cap
\zloc$, then for any $\eta\in\hk(\mext)$ we have
$$ \int_{\mlnd} \kappa(\eta e^{ \bom} ) =
 \frac{1}{|W|} \int_{N(c)}
\kappa(\nusym\eta e^{ \bom} ) = \frac{1}{|W|}\int_{N(V)}
\kappa(\nusym \eta e^{ \bom}\alpha)$$ where
\begin{equation} \label{e:ndef} N(c)= (\mext \cap \mu^{-1}(0))/T \end{equation}
for $\mu:K^{2g} \times \liek \to \liek$ given by minus the
projection onto $\liek$ and
$$N(V) =  (P^{-1}(V)\cap \zloc) /T.$$
Also $\alpha$ is a $T$-equivariantly closed form on $K^{2g}\times
\liek$ representing the $T$-equivariant Poincar\'{e} dual to
$\mext$, which is  chosen as in Proposition \ref{p3.2} so that the
support of $\alpha$ is contained in $P^{-1}(V)$ and has compact
intersection with $\zloc$.
\end{lemma}

 \begin{prop}  \label{p4.1} (\cite{JK2}, Proposition 6.3;\cite{Goldin})
  For any symplectic manifold $M$ equipped with a
Hamiltonian action of $T = T_n$ such that $T_{n-1}$
  acts locally freely on $\mu_{T_{n-1}}^{-1}(0)$,
the symplectic quotient $\mu_{T_n}^{-1}(0)/T_n$
  may be identified with the symplectic quotient of
$\mu_{T_{n-1}}^{-1}(0)/T_{n-1}$ by the
  induced Hamiltonian action of ${\hat{T_1}}$. Moreover if in addition
$T_n$ acts locally freely
  on $\mu_{T_n}^{-1}(0)$ then the ring homomorphism
$\kappa : H_{T_n}^*(M) \rightarrow
  H^*(\mu_{T_n}^{-1}(0)/T_n)$ which is the composition of restriction with the
natural isomorphism $\kappa : H_{T_n}^*( \mu_{T_n}^{-1}(0) ) \cong
  H^*(\mu_{T_n}^{-1}(0)/T_n)$
factors as
 \begin{equation} \label{e:factor}
    \kappa = {\hat{\kappa}}_1 \circ \kappa_{n-1}
\end{equation}
  where
  \[
    \kappa_{n-1} : H_{T_n}^*(M) \rightarrow
H_{T_n}^*(\mu_{T_{n-1}}^{-1}(0))
             \cong H_{{\hat{T}}_1 \times T_{n-1}}^*(\mu_{T_{n-1}}^{-1}(0))
             \cong H_{{\hat{T}}_1}^*(\mu_{T_{n-1}}^{-1}(0)/T_{n-1})
  \]
  and
  \[
   {\hat{\kappa}}_1 :  H_{{\hat{T}}_1}^*(\mu_{T_{n-1}}^{-1}(0)/T_{n-1})
\rightarrow
                           H^* \Bigl ( (\mu_{T_{n-1}}^{-1}(0) \cap
\mu_{{\hat{T_1}}}^{-1}(0))/(T_{n-1} \times {\hat{T}}_1 ) \Bigr )
                       \cong H^*(\mu_{T_n}^{-1}(0)/T_n)
  \]
are the corresponding compositions of restriction maps with
similar isomorphisms.
 \end{prop}

\begin{rem} \label{r4.2}
Let $T=T_n$ and $T_{n-1}$ be as in Proposition \ref{p4.1}. Note
that it is also true that
\begin{equation} \label{e:factor2}
\kappa = {\hat{\kappa}}_{n-1} \circ {\kappa_1} \end{equation}
where
$${\kappa_1}: H^*_T(M) \to H^*_{(T/\hat{T_1})} (\mu_{\hat{T_1}}^{-1}(0)/\hat{T_1})$$
and
$$\hat{\kappa}_{n-1}: H^*_{(T/\hat{T_1})} (\mu_{\hat{T_1}}^{-1}(0)/\hat{T_1} )
\to H^*_T(\mu_T^{-1}(0))\cong H^*(\mu_T^{-1}(0)/T)$$ are defined
in a similar way to $\hat{\kappa}_1$ and $\kappa_{n-1}$.
\end{rem}

 \begin{prop}
\label{p:4.2}
 {\bf(Dependence of symplectic quotients on parameters)} [Guillemin-Kalkman \cite{GK}; S. Martin \cite{martin}]
Let $M$ be a Hamiltonian $T$-space (where
 $T = U(1)$). If $n_0^T$ is the order of the stabilizer in $T$ of a generic point of
 $M$ then
 \[
   \int_{\mu_T^{-1}(\xi_1)/T} (\eta e^{\bar{\omega}})_{\xi_1}
   - \int_{\mu_T^{-1}(\xi_0)/T} (\eta e^{\bar{\omega}})_{\xi_0}
   = n_0^T \sum_{E \in {\mathcal E} : \xi_0 < \mu_T(E) < \xi_1} \mathrm{res}_{X=0} e^{\mu_T(E)X} \int_E \frac{\eta(X) e^\omega }{ e_E(X) }
 \]
 where $X \in {\Bbb C}$ has been identified with $2 \pi i X \in {\bf t} \otimes {\Bbb C}$
 and $\xi_0 < \xi_1$ are two regular values of the moment map.
Here ${\mathcal E}$ is the set of components of the fixed point
set of $\hat{T_1}$ on $M $.
 \end{prop}

Using Remark \ref{r4.2} the following lemma is proved exactly as
Lemma 6.7 of \cite{JK2}.

 \begin{lemma} \label{l:4.1}
  Suppose that $0$ is a regular value of  $\mu_{T}$ and that
 $\eta$ is an equivariant  cohomology class on $\mext$, which
is a polynomial in the classes $\tilde{a}_r$ and $\tilde{b}_r^j$.
Suppose also that $0$ is a regular value of $\mu_{T/\hat{T_1}}: E
\to \liet/\hat{\liet}_1 $ for all components $E$ of the fixed
point set of $\hat{T}_1$.
   If $V$ is a sufficiently small $T$-invariant neighbourhood of
  $c$ in $K$ so that $P^{-1}(V) \cap \mu^{-1}({\hat{\bf t}}_1)/T_{n-1}$
 is an orbifold and we define
  $N(V) = P^{-1}(V) \cap \mu^{-1}(0)/T$, then
  \[
    \int_{N(V)} \kappa(\eta e^{\bar{\omega}} e^{-Y_1} \alpha)
        = \int_{ P^{-1}(V) \cap \mu^{-1}({\hat{e}}_1)/T_n}
\kappa(\eta e^{\bar{\omega}} \alpha) \\
  \]
  \[
      = \int_{N(V)} \kappa(\eta e^{\bar{\omega}} \alpha)
            -n_0 \sum_{ E \in {\mathcal E} : -||{\hat{e}}_1||^2 < \langle {\hat{e}}_1,\mu(E) \rangle <0 }
            \int_{E\quott T_{n-1}}
   \kappa_{T/\hat{T_1} }  \mathrm{res}_{Y_1 =0} \frac{
\eta e^{\bar{\omega}} \alpha } {e_E}
  \]
  where ${\mathcal E}$ is the set of
components of the fixed point set of the action of ${\hat{T}}_1$
on $\ktg \times \liet$, and $e_E$ is the $\hat{T_1}$-equivariant
Euler class of the normal to $E$ in $\ktg$,
 while $n_0$
is the order of the subgroup of ${\hat{T}}_1/({\hat{T}}_1 \cap
T_{n-1})$
  that acts trivially on
$\ktg \times \liet$. Here $Y_1$ is a complex variable defined by
$<\hat{e_1},X> = Y_1$, where $X \in {\bf t} \otimes \CC$, and
$\alpha$ is
  the $T$-equivariantly
closed differential form on $\ktg \times {\bf t}$ given by
Proposition \ref{p3.2}
  which represents the equivariant Poincar\'e  dual  of
$\mext$, chosen so that the
  support of $\alpha$ is contained in $P^{-1}(V)$.
 \end{lemma}

\begin{rem} The proof of this lemma uses the fact that
the restriction of $P : \ktg \times {\bf k} \rightarrow K$ to
$\mu^{-1}({\bf t}) = \ktg \times {\bf t}$ is invariant
 under the translation $s_{\Lambda_0} : \ktg \times {\bf k} \rightarrow
\ktg \times {\bf k}$ defined by
 \[
     s_{\Lambda_0} : (m, \Lambda) \mapsto (m, \Lambda + \Lambda_0)
 \]
 for $\Lambda_0 \in \Lambda^I = \ker(\exp)$ in ${\bf t}$, and so is the polynomial $\eta$ in
the classes $\tilde{a}_r$ and $\tilde{b}_r^j$ (see \cite{JK2}).
\end{rem}

\begin{rem}
If $\mu_T: M \to \liet^*$ is a moment map for the action of a
torus $T$ on a symplectic manifold $M$, then we can add any
constant $\xi$ in $\liet^*$ to $\mu_T$ and get another moment map.
Then by applying Proposition \ref{p4.1} to $\mu_T - \xi$ we can
generalise the proposition to apply to
$\mu_{T_{n-1}}^{-1}(\xi_{n-1})$ and $\m_{\hat{T}_1}^{-1}(\xi_1)$
where $\xi_{n-1}$ and $\xi_1$ are the projections of $\xi$ into
$\liet_{n-1}$ and $\hat{\liet}_1$. A similar generalisation is
available for Lemma \ref{l:4.1}.
\end{rem}

From Lemma \ref{l:4.1}, we get
 \begin{eqnarray}
   \int_{N(V)} \kappa(\eta e^{\bar{\omega}} \alpha)
    &-& \int_{N(V)} \kappa(\eta e^{\bar{\omega}} e^{Y_1} \alpha)
= \int_{N(V)} \kappa(\eta e^{\bar{\omega}} (1-e^{Y_1}) \alpha) \nonumber \\
    &=& n_0 \sum_{ E \in {\mathcal E} : -||{\hat{e}}_1||^2 < \langle {\hat{e}}_1,\mu(E) \rangle <0 }
             \int_{E\quott T_{n-1}} \kappa_{n-1}
 \mathrm{res}_{Y_1 =0}\frac{ \eta e^{ \bar {\omega } } \alpha }{ { e_E }}
 \end{eqnarray}
 Therefore
 \begin{equation} \label{e:a7.14}
   \int_{N(V)} \kappa(\eta e^{\bar{\omega}} \alpha)
   =  n_0 \sum_{ E \in {\mathcal E} :
 -||{\hat{e}}_1||^2 < \langle {\hat{e}}_1,\mu(E) \rangle <0 }
          \int_{E\quott T_{n-1}}
  \kappa_{n-1} \mathrm{res}_{Y_1 =0}
\frac{ ( \eta e^{\bar{\omega} } \alpha) }{e_E(1-e^{Y_1} ) }
 \end{equation}

\subsection{Formulas for intersection pairings}

From \cite{JS} we have

 \begin{lemma} \label{l:4.2}
    Let $(M, \omega, \mu)$ be a quasi-Hamiltonian $K$-space and $
\hat{T_1}  \cong S^1$ be a circle subgroup of $K$. If $H$ is the
fixed point set of the
    adjoint action of $\hat{T_1}$ on $K$, then
$H$ is a Lie subgroup of $K$ and the fixed point set
$M^{\hat{T_1}}$ is a quasi-Hamiltonian $H$-space.
\end{lemma}

\begin{rem} In our particular case we take
$M = K^{2g}$, and $M^{\hat{T_1}}$ then has the form $H^{2g}$ where
$T \subset H$ and $\hat{T_1} \subset Z(H)$. Thus $T_{n-1} =
T/\hat{T_1}$ is a group of rank $n - 2$ with a quasi-Hamiltonian
action
 on $H^{2g}$. This enables us to perform an inductive argument.
\end{rem}
 \bigskip

\newcommand{\bracearg}[1]{ { [[ #1 ]] } }
\newcommand{\tildarg}[1]{ { [[ #1 ]]  } }
\newcommand{\tfq}{{ \tilde{f}_{ (q) } }}
\newcommand{\tf}{\tilde{f}}
\newcommand{\abk}{{\kappa}}
\newcommand{\tc}{\tilde{c}}

\newcommand{\tar}{\tilde{a}_r}
\newcommand{\tbr}{\tilde{b}_r^j}

The main result of \cite{JK2} is the following.
\begin{theorem} \label{mainab} (\cite{JK2}, Theorem 8.1)
Let $c=\diag\, (e^{2\pi i d/n},\ldots, e^{2\pi i d/n})$ where $$d
\in\{1,\ldots,n-1\}$$ is coprime to $n$, and suppose that $\eta\in
H^*_{K}(\mext)$ is a polynomial $Q(\tilde{a}_2,\ldots,\tilde{a}_n,
\tilde{b}_2^1,\ldots,\tilde{b}_n^{2g})$ in the equivariant
cohomology classes $\tilde{a}_r$ and $\tilde{b}_r^j$ for $2\leq
r\leq n$ and $1\leq j\leq 2g$ introduced in $\S$4.2. Then the
pairing
$$Q(a_2,\ldots,a_n,b_2^1,\ldots,b_n^{2g})\exp (f_2) [\mlnd]$$
is given by
$$\int_{\mlnd} \kappa (\eta e^{\bom} )
 =  \frac{(-1)^{n_+(g-1)}}{n!} \res_{Y_{1} =0} \dots \res_{Y_{n-1} =0}
\Bigl ( \frac{\sum_{w \in W_{n-1}} e^{\inpr{ \tildarg{w \tc},X}  }
 \int_{T^{2g}} \eta
e^{ \omega} } { \nusym^{2g-2} \prod_{1\leq j \leq n-1} ( \exp
(Y_j)-1 ) } \Bigr ), $$ where $n_+ = \frac{1}{2} n(n-1)$ is the
number of positive roots of $K=SU(n)$ and $X\in {\bf t}$ has
coordinates $Y_1=X_1-X_2,\ldots,Y_{n-1}=X_{n-1}-X_n$ defined by
the simple roots, while $W_{n-1} \cong S_{n-1}$ is the Weyl group
of $SU(n-1)$ embedded in $SU(n)$ in the standard way using the
first $n-1$ coordinates. The element $\tilde{c}$ is the unique
element of $\liet_n$  which satisfies $e^{2\pi i\tilde{c}} = c$
and belongs to the fundamental domain defined by the simple roots
for the translation action on $\liet_n$ of the integer lattice
$\L^I$. Also,  the notation $\bracearg{\gamma}$ means the unique
element which is in the fundamental domain defined by the simple
roots for the translation action  on $\liet_n$ of the integer
lattice and for
 which $\bracearg{\gamma}$ is equal to $\gamma$ plus some element of the
integer lattice; the notation $\nusym$  refers to the product of
all the positive roots of $SU(n)$.
\end{theorem}

Theorem \ref{mainab} can be proved inductively using
 Lemma \ref{l3} and (\ref{e:a7.14}), together with
Lemma \ref{l:4.2}. The proof reduces to proving the following
Proposition \ref{beg}, which is in fact more general because we
are no longer assuming that $c=\diag\, (e^{2\pi i d/n},\ldots,
e^{2\pi i d/n})$ (and therefore this proposition provides formulas
for  pairings in moduli spaces of parabolic bundles \cite{EK}).

\begin{prop} \label{beg} Let $c=\diag(c_1,\ldots,c_n) \in T$ be such that
the product of no proper subset of $c_1,\ldots,c_n$ is 1. If
$\eta(X)$ is a polynomial in the $\tar$ and $\tbr$ then
$$\int_{N(c)} \kappa ({\mathcal D}_n \eta e^{ \bom} ) = (-1)^{n_+(g-1)}
\res_{Y_{1} =0} \dots \res_{Y_{n-1} =0} \Bigl ( \frac{\sum_{w \in
W_{n-1}} e^{ \langle [[w\tc]] ,X  \rangle } \int_{T^{2g}} \eta e^{
\omega} } {{\mathcal D}^{2g-2} \prod_{1\leq j \leq n-1} ( \exp
(Y_j)-1 ) } \Bigr ), $$ where  $\tc = (\tc_1,\dots,\tc_n) \in
\liet_n$ satisfies $e^{2\pi i \tc}=c$ and belongs to the
fundamental domain defined by the simple roots for the translation
action on $\liet_n$ of the integer lattice $\Lambda^I$ and
$$N(c) = (\mu^{-1}(0) \cap \mext)/T.$$
The other notation is as in Theorem \ref{mainab}.
\end{prop}

\subsection{Extension to general pairings}

So far we have considered pairings of powers of the classes $a_r$,
$b_r^j$  and the K\"ahler class $f_2$. We now explain the
 general case.

\newcommand{\he}[1]{{\hat{e_{#1}} } }
\newcommand{\tarnox}{{\tilde{a}_r} }
\newcommand{\tbrnox}{{\tilde{b}_r^j}}

We  define $q \in S(\lieks)^K$ to be an invariant polynomial,
which is given in terms of the elementary symmetric polynomials
$\tau_j$ by \beq \label{9.0077} q(X) = \tau_2(X) + \sum_{r = 3}^n
\delta_r \tau_r(X). \eeq The associated element $\tfq $ of
$\hk(\mext)$ is defined by \beq \label{9.001} \tfq = \tf_2 +
 \sum_{r = 3}^n \delta_r \tf_r. \eeq
Here, the
 $\delta_r$ are formal
  nilpotent parameters: we
expand
 $\exp \tfq$ as a formal power series in the
$\delta_r$.

\begin{theorem} \label{t9.5}  (\cite{JK2}, Theorem 9.11(a))
For $q$ and $\tfq$ defined as above and
$$\eta = Q(\tilde{a}_2,\ldots,\tilde{a}_n,
\tilde{b}_2^1,\ldots,\tilde{b}_n^{2g})$$
 as
in Theorem \ref{mainab}, the pairing
$$Q(a_2,\ldots,a_n,b_2^1,\ldots,b_n^{2g})\exp (f_2 + \sum_{r=3}^n
\delta_r f_r) [\mlnd]$$ is given by
\begin{equation} \label{coff}
\frac{1}{n!}\int_{N(c) } \abk
 ( e^{\tfq} \nusym \eta )
=
 \frac{(-1)^{n_+(g-1)} }{n!}\sum_{w \in W_{n-1} }
   \res_{Y_{1} = 0} \dots \res_{Y_{n-1} = 0}
\frac { \int_{T^{2g}\times  \{ -\tildarg{w\tc}  \}  }\Bigl ( e^{
\tfq(X) }   \eta(X) \Bigr )  }{\nusym(X)^{2g-2} \prod_{j =
1}^{n-1} \exp(-B(-X)_{j} - 1)
 }, \end{equation}
when $c=\diag\, (e^{2\pi i d/n},\ldots, e^{2\pi i d/n})$ and $$d
\in\{1,\ldots,n-1\}$$ is coprime to $n$. In addition (\ref{coff})
is true more generally for any $c=\diag(c_1,\ldots,c_n) \in T$
such that the product of no proper subset of $c_1,\ldots,c_n$ is
1. Here $B(X)_j = -(dq)_X (\he{j}) $; we have used the fixed
invariant inner product on $\liek$ to identify $dq_X: \liet \to
\RR$ with an element of $\liet$ and thus define the  map $B: \liet
\to \liet.$ The other notation $\tildarg{\gamma}$ is as in Theorem
\ref{mainab} and Proposition \ref{beg}.
\end{theorem}

\begin{rem} \label{r:finitecover}
There are exact sequences
\begin{equation} \label{e:seqone}
1\to SU(n) \to U(n) \stackrel{det}{\longrightarrow} U(1) \to 1
\end{equation}
\begin{equation} \label{e:seqtwo}
1 \to {\bf  Z}_n \to SU(n) \times U(1) \to U(n) \to 1
\end{equation} and thus a finite covering map
\begin{equation} \label{e:covmap}
\mlnd \times {\rm Jac} \to \mnd
\end{equation}
with fiber ${\bf  Z}_n^{2g}$ which induces an isomorphism
$$H^*(\mnd) \cong H^*(\mlnd) \otimes H^*({\rm Jac})$$
(see \cite{AB}). As a result the cohomology of the space $\mnd$ is
related to that of $\mlnd$ by introducing the additional
generators $b_1^j \in H^1(\mnd)$ (corresponding to the generators
of the cohomology of the Jacobian) where $j = 1, \dots, 2g$.
Intersection pairings on $\mnd$ are related to the corresponding
pairings on $\mlnd$ by a factor $n^{2g}$ corresponding to the
order of the fiber in (\ref{e:covmap}).

Similarly we have a covering map
\begin{equation} \label{e:covpar}
\mu_{SU(n)}^{-1}({\tilde c})/T_{SU(n)} \times {\rm Jac} \to
\mu_{U(n)}^{-1}({\tilde c})/T_{SU(n)}
\end{equation}
with fiber ${\bf Z}_n^{2g}$.

The results of the remainder of this paper could be phrased
equally well in terms of the  moduli spaces $\mnd$ or the moduli
spaces $\mlnd$ of holomorphic vector bundles with fixed
determinant, though some care is needed when comparing the partial
desingularisations $\tilde{{\calm}}(n,d)$ and
$\tilde{{\calm}}_\Lambda(n,d)$: the covering map
$$\tilde{{\calm}}_\Lambda(n,d) \times {\rm Jac} \to \tilde{{\calm}}(n,d)$$
does  {\em not} induce an isomorphism from
$H^*(\tilde{{\calm}}(n,d))$ to $ H^*(\tilde{{\calm}}_\Lambda(n,d))
\otimes H^*({\rm Jac})$.

For simplicity, from now on we have chosen to restrict our
treatment to the moduli spaces $\mnd$.
\end{rem}

\renorm
\section{Pairings in intersection cohomology}

In this section, we show that the pairings in the intersection
cohomology $I\!H^*(\mnd)$ are given by essentially the same formulas
as in the nonsingular case reviewed in $\S$4, but with a small
shift.

In \cite{Kiem} it is shown that the weakly balanced condition
(\cite{Kiem} $\S$7 or \cite{jkkw} $\S$5) is satisfied for the
geometric invariant theoretic construction of $\mnd$. This means that the
subspace $V(n,d)=V_{M^{\rm ext}_{U(n)}}$ of $H^*_K(\mu^{-1}(0))\cong
H^*_{\overline{\mathcal{G}}(n,d)}(\mathcal{C}(n,d)^{ss})$ defined
in \cite{Kiem} is isomorphic to $I\!H^*(\mnd)$ where $K=SU(n)$ and
$\overline{\mathcal{G}}(n,d)={\mathcal{G}}(n,d)/U(1)$. To
define $V(n,d)$ we consider as in $\S$3 a set of representatives
\begin{equation}\label{5.1yh} R=\prod_{i=1}^q GL(m_i;\CC) \qquad
\text{where}\quad \sum_{i=1}^q m_in_i=n \quad \text{and}\quad
n_id=d_in
\end{equation}
of the conjugacy classes of reductive subgroups of
$\mathcal{G}_c(n,d)$ which occur as stabilizer groups of
semistable bundles.
As in $\S$3 we have
\[ Z^{ss}_R \cong \prod_{i=1}^q \mathcal{C}(n_i,d_i)^{ss}\]
and if $N_0^R = N_0$ is the connected component of the normaliser $N^R=N$ of $R$
in $\mathcal{G}_c(n,d)$ then
\[ N_0^R=\prod_{i=1}^q (GL(m_i;\CC)\times
\mathcal{G}_c(n_i,d_i))/\CC^*,\]
\[ N_0^R/R=\prod_{i=1}^q
\mathcal{G}_c(n_i,d_i)/\CC^*=\prod_{i=1}^q
\overline{\mathcal{G}}_c(n_i,d_i),\]
\[ \pi_0N^R=\prod_{j\ge 0, k\ge 0}\Sym(\#\{i:m_i=j, n_i=k\})\]
and hence \[ H^*_{N^R_0/R}(Z^{ss}_R)\cong \bigotimes_{i=1}^q
H^*_{\overline{\mathcal{G}}(n_i,d_i)}(\mathcal{C}
(n_i,d_i)^{ss}).\] The obvious maps
\[ \overline{\mathcal{G}}_c\times_{N^R} Z^{ss}_R\to
\overline{\mathcal{G}}_cZ^{ss}_R
\hookrightarrow \mathcal{C} (n,d)^{ss}\] give rise to a map
\[H^*_{\overline{\mathcal{G}}(n,d) }(\mathcal{C}(n,d)^{ss})
\to H^*_{N^R}(Z^{ss}_R) {\cong} [H^*_{N^R_0/R}(Z^{ss}_R)\otimes
H^*_R]^{\pi_0N^R}\hookrightarrow \left(\bigotimes_{i=1}^q
H^*_{\overline{\mathcal{G}}(n_i,d_i)}(\mathcal{C}
(n_i,d_i)^{ss})\right)\otimes H^*_R  \]
where $H^*_R = H^*(BR)$ denotes the $R$-equivariant cohomology of
a point.
The subspace $V(n,d)$ is
defined to be the intersection of the kernels of the compositions
\[H^*_{\overline{\mathcal{G}}(n,d) }(\mathcal{C}(n,d)^{ss})
\to \left(\bigotimes_{i=1}^q
H^*_{\overline{\mathcal{G}}(n_i,d_i)}(\mathcal{C}
(n_i,d_i)^{ss})\right)\otimes H^*_R\to \left(\bigotimes_{i=1}^q
H^*_{\overline{\mathcal{G}}(n_i,d_i)}(\mathcal{C}
(n_i,d_i)^{ss})\right)\otimes H^{\ge n_R}_R\] for all $R$ where
$H^{\ge n_R}_R=\oplus_{j\ge n_R}H^j_R $ and $n_R$ is given by
\[ n_R =  \dim_\CC H^1(\Sigma, \mathrm{End}'_{\oplus} E) -\dim_\CC R\]
which is easily computable using Riemann--Roch. If we denote the
equivariant universal bundle over $\mathcal{C}(n,d)^{ss}\times
\Sigma$ by $\mathbb{U}(n,d)$ then the restriction of the
generators of the equivariant cohomology rings can also be easily
computed from the equation of Chern classes
\[c(\mathbb{U}(n,d))|_{Z^{ss}_R}=\prod_{i=1}^qc(\CC^{m_i}\otimes
\mathbb{U}(n_i,d_i))\]
(see \S7 for the computation in the rank 2 case). Therefore, given
an equivariant cohomology class in $H^*_K(\mu^{-1}(0))\cong
H^*_{\overline{\mathcal{G}}(n,d)}(\mathcal{C}(n,d)^{ss})$, it is
straightforward to check if it belongs to $V(n,d)$ or not. In
particular, when the rank $n$ is $2$, we can write down $V(2,d)$
explicitly as worked out in \cite{Kiem2} (Theorem 5.3) using the results of
\cite{Kiem3}.

Let $\alpha,\beta$ be two classes in $H^*_K(M^{ext}_{U(n)})$ with
$\deg \alpha+\deg \beta=\dim_{\RR}\mnd$ such that their
restrictions to $\mu^{-1}(0)$ lie in $V(n,d)\cong I\!H^*(\mnd)$.
Then the restriction of their cup product $\alpha\beta$ to
$\mu^{-1}(0)$ also lies in $V(n,d)$ by \cite{jkkw} Theorem 14 or
\cite{Kiem} Theorem 5.3. By the de Rham model for intersection
cohomology, any top degree class in $I\!H^*(\mnd)$ is represented by
a differential form $\eta$, \emph{compactly supported} on the
smooth part $\mnd^s$ of $\mnd$; the pairing
$\langle\kappa(\alpha),\kappa(\beta)\rangle$ in $I\!H^*(\mnd)$ of
the classes $\kappa(\alpha)$ and $\kappa(\beta)$ represented by
$\alpha$ and $\beta$ is the integral of any such differential form
representing $\alpha \beta$. By pulling back $\eta$ to
$\mu^{-1}(0)^s$, $\alpha\beta$ is represented by an equivariant
differential form compactly supported on $\mu^{-1}(0)^s$ where
$\mu^{-1}(0)^s$ denotes the smooth part in $\mu^{-1}(0)$ so that
$\mu^{-1}(0)^s/K=\mnd^s$. By Martin's argument (see (2.9) above)
using the fibration
\[
\begin{CD}
K/T @>>> \mu^{-1}(0)^s/T\\
@. @V{\pi}VV \\
@. \mu^{-1}(0)^s/K=\mnd^s
\end{CD}
\]
the pairing $\langle\kappa(\alpha),\kappa(\beta)\rangle$ is given
by
\[
\int_{\mnd} \kappa(\alpha \beta) =\frac{1}{n!}\int_{\mu^{-1}(0)/T}
\kappa( \alpha \beta \cald).
\]
Let $\varepsilon\in \mathbf{t}^*$ be a regular value sufficiently
close to $0$. Then there is a surjective map
\[
\mu^{-1}(\varepsilon)/T\to \mu^{-1}(0)/T
\]
induced by the gradient flow of minus the normsquare $-|\mu|^2$ of
the moment map, which is a diffeomorphism over the smooth part
$\mu^{-1}(0)^s/T$.\footnote{One way to see the diffeomorphism is
to view the $T$-quotients as the moduli spaces of parabolic
bundles (see \cite{J}). When the underlying vector bundle of a
parabolic bundle is stable and the parabolic weight $\varepsilon$
is sufficiently small, the stability of the parabolic bundle does not
change as we move $\epsilon$ around $0$. This gives us the
diffeomorphism.} Hence we have
\[
\frac1{n!}\int_{\mu^{-1}(0)/T}\kappa(\alpha \beta \cald)
=\frac1{n!}\int_{\mu^{-1}(\varepsilon)/T} \kappa^{(\varepsilon)}
(\alpha\beta\cald).
\]
Therefore we deduce that the pairing
$\langle\kappa(\alpha),\kappa(\beta)\rangle$
 in $I\!H^*(\mnd)$ is given by
\begin{equation}\label{yh-5.1}
\langle\kappa(\alpha),\kappa(\beta)\rangle
=\frac1{n!}\int_{\mu^{-1}(\varepsilon)/T} \kappa^{(\varepsilon)}
(\alpha\beta\cald)
\end{equation}
for any $\varepsilon\in \mathbf{t}^*$ sufficiently close to $0$.

To compute the right hand side of (\ref{yh-5.1}) by using the
formulas described in $\S$4 from \cite{JK2} $\S$8 and $\S$9 (where $M^{\rm ext}_{SU(n)}$
is used instead of $M^{\rm ext}_{U(n)}$, i.e. the determinant of
semistable bundles is fixed there), we consider the fibration
\[
\begin{CD}
\mu_{\Lambda}^{-1}(\varepsilon)/T @>>> \mu^{-1}(\varepsilon)/T\\
@. @V{det}VV\\
@. {\rm Jac}=(S^1)^{2g}
\end{CD}
\]
where $\mu_{\Lambda}:M^{\rm ext}_{SU(n)}\to \mathbf{k}^*$ is the
moment map for $M^{\rm ext}_{SU(n)}$. Integrating over the fiber
first, we get
\[
\frac1{n!}\int_{\mu^{-1}(\varepsilon)/T} \kappa^{(\varepsilon)} (\alpha\beta\cald)
=\frac1{n!}\int_{(S^1)^{2g}}\int_{\mu_{\Lambda}^{-1}(\varepsilon)/T}
\kappa^{(\varepsilon)} (\alpha\beta\cald ).
\]
The inner integral is given by Theorem \ref{t9.5} and so we have
proved the following.

\newcommand{\bk}{\mathbf{k} }

\begin{theorem} \label{t:t5.1}
Let $\alpha,\beta$ be two classes in
$H^*_{{\overline{\mathcal{G}}(n,d)}}({\mathcal{C}}(n,d))$ with
$\deg \alpha+\deg \beta=\dim_{\RR}\mnd$ such that their
restrictions to ${\mathcal{C}}(n,d)^{ss}$
 lie in $V(n,d)\cong I\!H^*(\mnd)$. For
$\mathbf{k}=(k_2,\cdots,k_n)\in \ZZ^{n-1}_{\ge 0}$, let
$f^\bk=\prod_{r=2}^nf_r^{k_r}$, let $\bk !=\prod_{r=2}^nk_r!$ and
let $\delta^\bk=\prod_{r=3}^n\delta^{k_r}_r$. Write
$$\alpha\beta=\sum_\bk  Q_\bk(a_2,\cdots,a_n,b_1^1,\cdots,b_n^{2g})\frac{f^\bk}{\bk !} $$
where $Q_\bk$ is a polynomial of the Atiyah--Bott
classes $a_r, b_r^j \in
H^*_{{\overline{\mathcal{G}}(n,d)}}({\mathcal{C}}(n,d))$ defined
in Remark \ref{abf lift to equivariant coh}. Here, and later, for simplicity of notation we think of this as representing the class 
$$
\sum_\bk  Q_\bk(\tilde a_2,\cdots,\tilde a_n,\tilde b_1^1,\cdots,\tilde b_n^{2g})\frac{\tilde f^\bk}{\bk !}
$$
in the equivariant cohomology $H^*_K(\mext)$ of the extended moduli space.  Then, using Theorem \ref{t9.5}, the pairing $\langle\kappa(\alpha),\kappa(\beta)\rangle$
of $\kappa(\alpha)$ and $\kappa(\beta)$ in $I\!H^*(\mnd)$ is given
by
\[
 \frac{(-1)^{n_+(g-1)}
}{n!}\sum_\bk
 \mathrm{Coeff}_{\delta^\bk}\left( \sum_{w \in W_{n-1} }
   \res_{Y_{1} = 0} \dots \res_{Y_{n-1} = 0}
\frac { \int_{(S^1)^{2g}}\int_{T^{2g}\times  \{ -\tildarg{w\tc}
\}  }\Bigl ( e^{ \tfq(X) }   Q_\bk(X) \Bigr )
}{\nusym(X)^{2g-2} \prod_{j = 1}^{n-1} \exp( -B(-X)_{j} - 1)
 }\right)\]
for any $\tc$ sufficiently close to $\tc_0$ where $\tilde{c}_0$ is
the unique element of $\liet_n$ which satisfies $e^{2\pi
i\tilde{c}_0}=\diag\, (e^{2\pi i d/n},\ldots, e^{2\pi i d/n})$ and
belongs to the fundamental domain defined by the simple roots for
the translation action on $\liet_n$ of the integer lattice $\L^I$.
Here
$$ \tfq =
\tf_2 +
 \sum_{r = 3}^n \delta_r \tf_r$$
as at (\ref{9.001}) where the  $\delta_r$ are formal
  nilpotent parameters, and we
expand
 $\exp \tfq$ as a formal power series in the
$\delta_r$; the coefficient
of $\delta^\bk$
is denoted by $\mathrm{Coeff}_{\delta^\bk}$ and $B(X)_j = -(dq)_X (\he{j}) $ as in Theorem
\ref{t9.5}.
\end{theorem}

In particular, if $\alpha\beta$ is a polynomial $Q$ in the classes
$a_2,\cdots, a_n, b_1^1,\cdots,b_n^{2g}$ then by Proposition
\ref{beg} we have
\[
\langle\kappa(\alpha),\kappa(\beta)\rangle=
\frac{(-1)^{n_+(g-1)}}{n!} \res_{Y_{1} = 0} \dots \res_{Y_{n-1} =
0} \Bigl ( \frac{\sum_{w \in W_{n-1}} e^{\inpr{ \tildarg{w \tc},X}
}
 \int_{(S^1)^{2g}}\int_{T^{2g}} Q
e^{ \omega} } { \nusym^{2g-2} \prod_{1\leq j \leq n-1} ( \exp
(Y_j)-1 ) } \Bigr ) \] as in the coprime case, while if
$\alpha\beta$ is the product of $(f_2)^k$ with a polynomial
$Q$ in the classes
$a_2,\cdots, a_n, b_1^1,\cdots,b_n^{2g}$  then
\[
\langle\kappa(\alpha),\kappa(\beta)\rangle=
\frac{(-1)^{n_+(g-1)}k!}{n!} \res_{Y_{1} = 0} \dots \res_{Y_{n-1}
= 0} \Bigl ( \frac{\sum_{w \in W_{n-1}} e^{\inpr{ \tildarg{w
\tc},X} }
 \int_{(S^1)^{2g}}\int_{T^{2g}} Q
e^{ \omega} } { \nusym^{2g-2} \prod_{1\leq j \leq n-1} ( \exp
(Y_j)-1 ) } \Bigr ). \]

\renorm

\renorm
\section{Pairings on the partial desingularisation $\tilde{\calm}(n,d)$}

In this section we will study pairings on the partial desingularisation
$\tilde{{\cal M}}(n,d)$ using the method described in $\S$8 of
\cite{jkkw}.
In the notation of $\S$2 above, \cite{jkkw} $\S$8 provides a
method for calculating pairings
$ \kappa_{{\tilde{M}}}(\a\b)[{\tilde{M}}/\!/G]$ in the cohomology
$H^*(\tilde M \git G)$ of classes $ \kappa_{{\tilde{M}}}(\a)$ and
$ \kappa_{{\tilde{M}}}(\b)$  in the image of the composition
$$H^*_K(M) \rightarrow H^*_K(\tilde M) \rightarrow H^*(\tilde M \git
G)$$
of the pullback from $M$ to $\tilde{M}$ and the map $\kappa_{\tilde M}$.
The first observation is that
\begin{equation} \label{new7}
\kappa_{{\tilde{M}}}(\a\b)[{\tilde{M}}/\!/G] = \frac{
1}{ |W|} \kappa_{{\tilde{M}}}^T(\a \b
{\mathcal{D}})[{\tilde{\mu}}^{-1}(0)/T]= \frac{
(-1)^{n_+} }{ |W|} \kappa_{{\tilde{M}}}^T(\a \b
{\mathcal{D}}^2)[{\tilde{M}}/\!/T_c] \end{equation}
by (2.8) and (2.9), and the next is that if  $\xi$
is any regular value of the $T$-moment map $\mu_T$ for $M$ and
we choose
the symplectic structure appropriately (as a sufficiently small perturbation
of the pullback of the symplectic structure on $M$; cf. \cite{jkkw}
\S4) then
\begin{equation} \label{new2}
\kappa_{{\tilde{M}},\xi}^T(\a \b
{\mathcal{D}})[{\tilde{\mu}}^{-1}(\xi)/ T] = \kappa_{{{M},\xi}}^T(\a \b
{\mathcal{D}})[\mu^{-1}(\xi)/ T],
\end{equation}
where the latter expression can in our case (for $M = M^{\rm ext}_{U(n)}$)
be calculated as in $\S$5.

\begin{rem} \label{fcknew}
Indeed the proof of Theorem \ref{t:t5.1} tells us that if
$\alpha$ and $\beta$ are two classes in
$H^*_{{\overline{\mathcal{G}}(n,d)}}({\mathcal{C}}(n,d))$ with
$\deg \alpha+\deg \beta=\dim_{\RR}\mnd$ and
$$\alpha\beta=\sum_\bk
Q_\bk(a_2,\cdots,a_n,b_1^1,\cdots,b_n^{2g})
\frac{f^\bk}{\bk !} $$
for polynomials $Q_\bk$ where the sum runs over
$\mathbf{k}=(k_2,\cdots,k_n)\in \ZZ^{n-1}_{\ge 0}$,  then the
pairing $\kappa_{{{M},\xi}}^T(\a \b {\mathcal{D}})[\mu^{-1}(\xi)/
T]$ is given by
\begin{equation}\label{pertpairing}
\frac{(-1)^{n_+(g-1)}
}{n!} \sum_\bk
 \mathrm{Coeff}_{\delta^\bk}\left(\sum_{w \in W_{n-1} }
   \res_{Y_{1} = 0} \dots \res_{Y_{n-1} = 0}
\frac { \int_{(S^1)^{2g}}\int_{T^{2g}\times  \{ -\tildarg{w\tc}
\}  }\Bigl ( e^{ \tfq(X) }   Q_\bk(X) \Bigr )
}{\nusym(X)^{2g-2} \prod_{j = 1}^{n-1} \exp (-B(-X)_{j} - 1)
 }\right)
 \end{equation}
where $\tc = \tc_0 + \xi$ and $\tilde{c}_0$ is
the unique element of $\liet_n$ which satisfies $e^{2\pi
i\tilde{c}_0}=\diag\, (e^{2\pi i d/n},\ldots, e^{2\pi i d/n})$ and
belongs to the fundamental domain defined by the simple roots for
the translation action on $\liet_n$ of the integer lattice $\L^I$.
Here $f^\bk=\prod_{r=2}^nf_r^{k_r}$, $\bk
!=\prod_{r=2}^nk_r!$ and $\delta^\bk=\prod_{r=3}^n\delta^{k_r}_r$
where the  $\delta_r$ are formal
  nilpotent parameters, and
$$ \tfq =
\tf_2 +
 \sum_{r = 3}^n \delta_r \tf_r$$
as at (\ref{9.001}). Finally $B(X)_j = -(dq)_X (\he{j}) $ as in
Theorem \ref{t9.5}.
\end{rem}

This means that to calculate the pairing of the cohomology classes
induced by $\alpha$ and $\beta $ in $H^*(\tilde{{\cal M}}(n,d))$
it suffices to calculate the difference between
\begin{equation} \label{anew}
\kappa_{{\tilde{M}},\xi}^T(\a \b
{\mathcal{D}})[{\tilde{\mu}}^{-1}(\xi)/ T]
\end{equation}
and
\begin{equation} \label{bnew}
\kappa_{{\tilde{M}}}^T(\a \b
{\mathcal{D}})[{\tilde{\mu}}^{-1}(0)/ T]
\end{equation}
when $M=M^{\rm ext}_{U(n)}$.

Since the partial desingularisation process takes place in stages, it is
easiest to consider a
single stage of the construction, when
a blow-up along the proper transform $\hat{Z}_R/\!/N$ of $Z_R/\!/N$ with
$$ R= GL(m_1;\CC) \times ... \times GL(m_q;\CC) $$
as at (\ref{labelR}) results in $\hat{M}/\!/G$ which in our case
can be written as $\hat{{\cal M}}(n,d) = \hat{{\cal C}}/\!/{\cal
G}_c$. (Note that in \cite{jkkw} the superscript $\hat{ }$ in
$\hat{Z}_R/\!/N$ is omitted since, for simplicity, it is assumed
there that the blow up along $\hat{Z}_R/\!/N$ is the first in the
partial desingularisation process and hence $\hat{Z}_R/\!/N =
{Z}_R/\!/N$.) Let $\hat{\mu}$ and $\hat{\mu}_T$ be the moment maps
for the actions of $K$ and $T$ on $\hat{M}$; then as at
(\ref{new2}) when $\xi$ is a regular value of $\mu_T$ we have
\begin{equation} \label{new4}
\kappa_{{\hat{M}},\xi}^T(\a \b
{\mathcal{D}})[{\hat{\mu}}^{-1}(\xi)/ T] = \kappa_{M,\xi}^T(\a \b
{\mathcal{D}})[{{\mu}}^{-1}(\xi)/T]
\end{equation}
provided that
we have made a sufficiently small perturbation of the pullback to $\hat{M}$
of the symplectic form  on $M$ (where \lq sufficiently small' depends
on $\xi$).
We choose $\xi \in \liets$  to lie in a connected
component $\Delta_{(i)}$ of the set of regular values of $\mu_T$
for which $0 \in \bar{\Delta}_{(i)}$, and choose $\hat{\xi} \in
\liets$ to lie in the intersection of $\Delta_{(i)}$ and a
connected component of the set of regular values of $\hat{\mu}_T$
which contains $0$ in its closure. We cannot necessarily choose
$\hat{\xi} = \xi$ because the choices of symplectic
structure on $\hat{M}$ and the moment maps
$\hat{\mu}$ and $\hat{\mu}_T$ depend on $\xi$, but it is enough to
calculate the difference
\begin{equation} \label{new5}
\kappa_{{\hat{M}},\xi}^T(\a \b
{\mathcal{D}})[{\hat{\mu}}^{-1}(\xi)/ T] -
\kappa_{{\hat{M},\hat{\xi}}}^T(\a \b
{\mathcal{D}})[{\hat{\mu}}^{-1}({\hat{\xi}})/T],
\end{equation}
since repeating this for each stage and using (\ref{new2}) and
(\ref{new4}) will give us the
difference between $\kappa_{{{M},\xi}}^T(\a \b
{\mathcal{D}})[{{\mu}}^{-1}({\xi})/ T]$
and
$\kappa_{{\tilde{M}},\tilde{\xi}}^T(\a \b
{\mathcal{D}})[{\tilde{\mu}}^{-1}({\tilde{\xi}})/ T]$
for any $\tilde{\xi}$ in a connected component of the set of regular
values of
$\tilde{\mu}_T$ which contains 0 in its closure. Since 0 is itself a
regular value
of $\tilde{\mu}_T$, we can choose $\tilde{\xi}$ to be 0 and thus
calculate the pairing (\ref{new7}).

Recall that the image $\hat{\mu}_T(\hat{M})$ of the moment map
$\hat{\mu}_T$ is a convex polytope which is divided by walls of
codimension one into subpolytopes whose interiors consist of
regular values of $\hat{\mu}_T$. The pairings we wish to calculate
are unchanged as $\xi$ varies within a connected component of the
set of regular values of $\hat{\mu}_T$, so it is enough to be able
to calculate the change  as $\xi$ crosses a wall of codimension
one. Any such wall is of the form
$$\hat{\mu}_T ( \hat{M}_1)$$
where $\hat{M}_1$ is a connected component of the fixed point set
in $\hat{M}$ of a circle subgroup $T_1$ of $T$, and it is shown in
\cite{jkkw} Lemma 23 that in order to calculate the difference
(\ref{new5}) the only wall crossing terms we need to consider
correspond to components $\hat{M_1}$ of fixed point sets of circle
subgroups $T_1$ satisfying
$$\emptyset \ne \pi(\hat{M_1}) \cap M^{ss} \subset G \hat{Z}_R^{ss}$$
so that $\hat{M_1}$ is contained in the exceptional divisor of $\pi:
\hat{M} \to M$.
Moreover if $N^{T_1}_0$ is the identity component of the normaliser
of $T_1$ then the subset
$$\{g \in G: T_1 \subseteq gRg^{-1} \}$$
of $G$ is the disjoint union
$$\bigsqcup_{1\le i\le m} N^{T_1}_0g_iN^R$$
of finitely many
double cosets for the left $N^{T_1}_0$ action and the right $N^R$
action, and
 the $T_1$-fixed point set in
$G\hat{Z}^{ss}_{R}$ is the disjoint union
$$\bigsqcup_{1\le i\le m}N^{T_1}_0 \hat{Z}^{ss}_{R_i}$$
where $R_i=g_iRg_i^{-1}\supset T_1$.
Let   the $T_1$-eigenbundles
 of the restriction  to $N_0^{T_1} \hat{Z}_{R_i}^{ss}$
of the normal bundle to  ${\cal G}_c \hat{Z}^{ss}_R$ be denoted by
\begin{equation} \label{defn8.6} {\cal W}_{i,j}\to N^{T_1}_0
\hat{Z}^{ss}_{R_i} \end{equation}
for $1 \le j \le l_i$. If $\hat{M_1}$ is a component of a fixed point
set of a circle
subgroup $T_1$ satisfying
$ \emptyset \ne \pi(\hat{M_1}) \cap M^{ss} \subset G \hat{Z}_R^{ss}$
then
$$ \hat{M_1} \cap \pi^{-1} (M^{ss}) = \PP {\cal W}_{i,j} $$
for some $i,j,$ and it is shown in \cite{jkkw} that the
corresponding wall crossing term
is
\begin{equation} \label{star}
\int_{\PP {\cal W}_{i,j}\quott_{\xi_1}  {(T/T_1)_c} }
 \kappa^{T/T_1}_{\PP {\cal W}_{i,j}, \xi_1}
\left ( {\rm res}_{X_1 = 0}  \frac{\alpha \beta \cald^2} {e_{\PP {\cal
W}_{i,j}}  }
 \right ) \end{equation}
where $e_{\PP {\cal W}_{i,j}}$ is the $T$-equivariant Euler class of the
normal bundle
to $\PP {\cal W}_{i,j}$, and the wall is crossed at $\xi_1 + \xi_2$
where
 $\xi_2$ is the constant $T_1$-component of $\hat{\mu}$ on
$\hat{M}_1$ and $\xi_1$ is orthogonal to the Lie algebra of $T_1$.

From \eqref{new7}, \eqref{new2} and \eqref{pertpairing}, we deduce
the following.
\begin{theorem}\label{mthmsect6}
Let $\alpha,\beta$ be two classes in
$H^*_{{\overline{\mathcal{G}}(n,d)}}({\mathcal{C}}(n,d))$ with
$\deg \alpha+\deg \beta=\dim_{\RR}\mnd$ and
$$\alpha\beta=\sum_\bk
Q_\bk(a_2,\cdots,a_n,b_1^1,\cdots,b_n^{2g})
\frac{f^\bk}{\bk !} $$ where $Q_\bk$ is a
polynomial.  Then using the notation of Remark
\ref{fcknew}, the intersection pairing on the partial
desingularisation $\tilde{\calm}(n,d)$ is given by
\[
\langle\kappa(\alpha),\kappa(\beta)\rangle= \sum_\bk
 \mathrm{Coeff}_{\delta^\bk}\left(\frac{(-1)^{n_+(g-1)}
}{n!}\right.
\]
\[
\left. \sum_{w \in W_{n-1} }
   \res_{Y_{1} = 0} \dots \res_{Y_{n-1} = 0}
\frac { \int_{(S^1)^{2g}}\int_{T^{2g}\times  \{ -\tildarg{w\tc} \}
}\Bigl ( e^{ \tfq(X) }   Q_\bk(X) \Bigr )
}{\nusym(X)^{2g-2} \prod_{j = 1}^{n-1} (\exp -B(-X)_{j} - 1)
 }\right)\]
\[
-(-1)^{n(n-1)/2}\sum_{R\in \mathcal{R}}\sum _{T_1}\sum_{i,j}
\int_{\PP {\cal W}_{i,j}\quott_{\xi_1}  {(T/T_1)_c} }
 \kappa^{T/T_1}_{\PP {\cal W}_{i,j}, \xi_1}
\left ( {\rm res}_{X_1 = 0}  \frac{\alpha \beta \cald^2} {e_{\PP
{\cal W}_{i,j}}  }
 \right )
\]
where the last sum runs over all reductive groups $R$ coming from
partitions of $n$ (see \eqref{labelR}), over all circle subgroups
$T_1$ of $T$ that appear as the stabiliser of a point in the
exceptional divisor of the blowup of $G\hat{Z}^{ss}_R$ and over
all $i,j$ such that $\hat{\mu}_T(\PP \cal W_{i,j})$ is a wall
between $0$ and $\hat{\xi}$.
\end{theorem}

The computation of the intersection pairing on the partial
desingularisation $\tilde{\calm}(n,d)$ is completed by computing
the last terms (the wall crossing terms) in the above theorem. In
the subsequent section, we will see how we can calculate the wall
crossing terms.

\renorm
\section{Wall crossing terms}

The purpose of this section is to compute the wall crossing term
\[
\int_{\PP {\cal W}_{i,j}\quott_{\xi_1}  {(T/T_1)_c} }
 \kappa^{T/T_1}_{\PP {\cal W}_{i,j}, \xi_1}
\left ( {\rm res}_{X_1 = 0}  \frac{\alpha \beta \cald^2} {e_{\PP
{\cal W}_{i,j}}  }
 \right )
 \]
which appears in Theorem \ref{mthmsect6} and hence to complete our calculation of the intersection pairing
on the partial desingularisation $\tilde{\calm}(n,d)$.

\begin{rem}\label{extra2} It is noted in \cite{jkkw} Remarks 26 and
27 that
this wall crossing term
 is an integral over a
quotient of the form $\hat{\mu}_T^{-1}(\xi_1 + \xi_2) \cap
\hat{M}_1 /T$ where $\xi_1$ and $\xi_2$ can be taken to be
arbitrarily close to 0. Since $\PP {\cal W}_{i,j}$
 is a projective
bundle over $N_0^{T_1} \hat{Z}_{R_i}^{ss}$,  it is helpful to use
the method of reduction to a maximal torus (cf. (2.8)) to relate
integrals over quotients of $\PP {\cal W}_{i,j}$ by $T$ to
integrals over quotients of $\PP {\cal W}_{i,j}$ by $N_0^{T_1}$,
{\em provided that} $\xi_1$ is centralised by $N_0^{T_1}$, so that
the quotient
$$\PP {\cal W}_{i,j} /\!/_{\xi_1} (N_0^{T_1}/(T_1)_c) =
\hat{\mu}^{-1}(\xi_1 + \xi_2) \cap \hat{M}_1/(N_0^{T_1} \cap K)$$
is well defined. Luckily in our situation we can use (\ref{stabsubset})
to allow us to assume that
$\xi_1$ is centralised by $N_0^{T_1}$, and this simplifies the
calculations described in \cite{jkkw} $\S$8
considerably. The wall crossing term (\ref{star}) becomes
\begin{equation} \label{*}
 \frac{(-1)^{n_+^{N_0^{T_1}}}}{|W_{N_0^{T_1}}|}
\int_{\PP {\cal W}_{i,j}\quott_{\xi_1}  {(N_0^{T_1}/(T_1)_c)} }
 \kappa^{T/T_1}_{\PP {\cal W}_{i,j}, \xi_1}
\left ( {\rm res}_{X_1 = 0}  \frac{\alpha \beta \cald^2} {({\cal
D}_{N_0^{T_1}})^2 e_{\PP {\cal W}_{i,j}}  } \right)
\end{equation}
where $n_+^{N_0^{T_1}}$ is the number of positive roots of
$N_0^{T_1}$, while ${\cal D}_{N_0^{T_1}}$ is the product of the
positive roots of $N_0^{T_1}$ and $W_{N_0^{T_1}}$ is the Weyl
group of $N_0^{T_1}$. Moreover by \cite{jkkw} (8.13) and Lemma 31
there is a fibration
\begin{equation} \label{commutdiag}
{\PP {\cal W}}_{i,j}  \quott_{\xi_1} (N_0^{T_1}/(T_1)_c) \cong
{\PP {\cal W}}_{i,j}  \mid_{\hat{Z}_{R_i}^{ss} } \quott N_0^{T_1} \cap
N^{R_i}
\stackrel{\Psi}{\longrightarrow} \hat{Z}_{R_i} \quott N_0^{T_1} \cap
N^{R_i}
\end{equation}
with fiber ${\PP ({\cal W}_{i,j})_x} \quott Stab
(x)\cap N_0^{T_1} \cap N^{R_i} . $ Thus we can calculate (\ref{*}) by
integrating
over the fibers of $\Psi$.
\end{rem}

Recall from (\ref{defn8.6}) that ${\cal W}_{i,j}\to N^{T_1}_0
\hat{Z}^{ss}_{R_i}$
is a $T_1$-eigenbundle
 of the restriction  to $N_0^{T_1} \hat{Z}_{R_i}^{ss}$
of the normal bundle to  ${\cal G}_c \hat{Z}^{ss}_R$, and
$R_i=g_iRg_i^{-1}\supset T_1$. Let us drop the indices $i$ and $j$
and consider one wall crossing term (\ref{*}) for
$$ R_i = R \cong GL(m_1;\CC) \times ... \times GL(m_q;\CC). $$

We have
\begin{equation} \label{stabx} \stab(x) = R \mbox{ for every } x \in
Z_R^s \end{equation}
and
by (\ref{N_0}) and (\ref{N/N_0}) the normaliser $N=N^R$ of $R$ in
${\cal G}_c$ has identity component
$$ N^R_0 \cong \prod_{1 \leq i \leq q} (GL(m_i;\CC) \times {\cal
G}_c(n_i,d_i))/\CC^*
$$
and
$$ \pi_0(N^R) = \prod_{j\geq 0, k \geq 0} \Sym({\#\{i:m_i=j \mbox{ and
}n_i=k\}})
.$$
Moreover
\begin{equation} \label{prodmod}  \hat{Z}_R \quott N^R \cong \left[
\prod_{i=1}^q \tilde{{\cal M}}(n_i,d_i)
\right] / \pi_0(N^R) \end{equation}
where $\left[ \prod_{i=1}^q \tilde{{\cal M}}(n_i,d_i)
\right] $ is the \lq blow up along diagonals' of $ \prod_{i=1}^q
\tilde{{\cal M}}(n_i,d_i) $,
defined as follows (cf. [38] Definition 3.9 and Lemma 3.11, but noting
that the definition
given in [38] of the \lq blow up along diagonals' is incorrect).

\begin{definition} \label{bdiag}
Let $Y_1,...,Y_s$ be quasi-projective varieties and let $\Pi = \{
I_1,...,I_k \}$ be a partition
of $\{1,...,s\}$ with associated equivalence relation $\sim$ such that
$Y_i = Y_j$ if and only if $i \sim j$. Given any partition $\Pi' = \{
I'_1,...,I'_\ell \}$
of $\{1,...,s\}$ which refines $\Pi$, there is a nonsingular closed
subvariety $A_{\Pi'}$ of the
product $\prod_{i=1}^s Y_i$ defined by
$$A_{\Pi'} = \{ (y_1,...,y_s) \in \prod_{i=1}^s Y_i : y_i = y_j \mbox{
if } i \sim' j \}$$
where $\sim'$ is the equivalence relation on $\{1,...,s\}$ induced by
the partition
$\Pi'$. Let
$$\left[ \prod_{i=1}^s Y_i \right]$$
denote the result of blowing up $\prod_{i=1}^s Y_i$ along the proper
transforms of the subvarieties
$A_{\Pi'}$ (where $\Pi'$ runs over all refinements of $\Pi$) in increasing order
of dimension. We call
$\left[ \prod_{i=1}^s Y_i \right]$ the product of $Y_1,...,Y_s$ blown up
along all diagonals.
\end{definition}

We saw in $\S$3 that the wall $\mu_T(\hat{M}_1)$ lies in an affine
hyperplane $\beta + \beta^{\perp}$ in $\liet_R$ where $\beta$ generates
$T_1$ and is
determined by a partition
$$\{ \Delta_{h,m}:(h,m) \in J \}$$
of $\{1,...,M \}$,  a nonempty subset $S$ of $\{ (i,j) \in \ZZ \times
\ZZ: 1 \leq i,j \leq M \}$,
and a directed graph $G(S)$ with vertices $1,...,M$ and directed
edges from $i$ to $j$ whenever $(i,j)\in S$, as in Proposition
\ref{5.1}.
Recall from Remark \ref{rem1}
that the graph $G(S)$ is connected and so we can
omit the index $h$ and take $J$ to be of the form
$$J = \{1, \ldots, t\}.$$
Our aim is to calculate the wall crossing term (\ref{*}) by integrating
over the
fibers
$${\PP ({\cal W})_x} \quott Stab
(x)\cap N_0^{T_1} \cap N^{R} = {\PP ({\cal W})_x} \quott
 N_0^{T_1} \cap N^{R}  $$
 of the fibration
$$ \Psi:
{\PP {\cal W}}  \quott_{\xi_1} (N_0^{T_1}/(T_1)_c) \cong
{\PP {\cal W}}  \mid_{\hat{Z}_{R}^{ss} } \quott N_0^{T_1} \cap
N^{R}
\longrightarrow \hat{Z}_{R} \quott N_0^{T_1} \cap N^{R}
$$
defined at (\ref{commutdiag}) (see Remark \ref{extra2} and (\ref{stabx})),
where
${\cal W} = {\cal W}_{i,j}$ and
$N_0^{T_1}$ is the identity component of the normaliser $N^{T_1}$ in
${\cal G}_c$ of the one-parameter subgroup $T_1$ generated by $\beta$,
so that
\begin{equation} \label{N=stabbeta}
N_0^{T_1} = N^{T_1} = \stab \beta .\end{equation}

From \cite{jkkw} Lemma 29 and Corollary 1 and from (\ref{prodmod})
above we know that the connected component $R (N_0^{T_1} \cap
N^{R} )_0$ has finite index in $N^R$ and that there are
isomorphisms
$$ \hat{Z}_{R} \quott (N_0^{T_1} \cap N^{R} )_0 \cong \hat{Z}_{R} \quott
N_0^{R} \cong
\left[  \prod_{i=1}^q \tilde{{\cal M}}(n_i,d_i) \right]
 $$
and
$$ \hat{Z}_{R} \quott N_0^{T_1} \cap N^{R}  \cong
\left[  \prod_{i=1}^q \tilde{{\cal M}}(n_i,d_i)
\right]/ \pi_0( N_0^{T_1} \cap N^{R} )$$
where
$$ \pi_0(N_0^{T_1} \cap N^{R} ) = \prod_{j_1\geq 0, j_2 \geq 0, j_3 \geq
0} \Sym({\#\{(i,k):m_i=j_1 \mbox{ and }
m^k_i = j_2 \mbox{ and }n_i=j_3\}})$$
(cf. (3.5)). Thus once we have reduced to integrals over
$ \hat{Z}_{R} \quott N_0^{T_1} \cap N^{R} $ by integrating over
the fibers of $\Psi$ we can complete the calculation by using induction
on $n$ to compute integrals over the product $\prod_{i=1}^q \tilde{{\cal M}}(n_i,d_i)$
and its blowup along diagonals $[\prod_{i=1}^q \tilde{{\cal M}}(n_i,d_i)]$.

We can integrate (\ref{*}) over the fibers
${\PP ({\cal W})_x} \quott
 N_0^{T_1} \cap N^{R}  $
of $\Psi$
by using the formula (\ref{res})
for pairings on the symplectic quotient of a projective space.
In this calculation
 the terms $\alpha$, $\beta$, ${\cal D}$ and ${\cal
D}_{N_0^{T_1}}$ all restrict to equivariant classes on the projective
space ${\PP ({\cal W})_x} $ which are pulled back from
the equivariant cohomology of a point
and are easy to calculate. The remaining term to consider is the equivariant
Euler class
$e_{\PP {\cal W}}$ of the normal bundle to ${\PP {\cal W}}$.

Let $ E \cong (\CC^{m_1} \otimes D_1) \oplus \cdots \oplus
(\CC^{m_q} \otimes D_q) $ represent an element of ${\cal G}_c
Z_R^s$ as at (3.2) above, with  $D_1, \ldots ,D_q$ all stable and
not isomorphic to one another, and $D_i$ of rank $n_i$ and degree
$d_i$. Recall from \S 3 that ${\cal C}$ is an infinite dimensional
affine space, and if we fix a $C^{\infty}$ identification of the
fixed $C^{\infty}$ hermitian bundle ${\cal E}$ with
$\bigoplus_{i=1}^q \CC^{m_i} \otimes D_i$ then we can identify
${\cal C}$ with the infinite dimensional vector space
$$\Omega^{0,1}(\mbox{End}(\bigoplus_{i=1}^q \CC^{m_i} \otimes D_i))$$
in such a way that the zero element of
$\Omega^{0,1}(\mbox{End}(\bigoplus_{i=1}^q \CC^{m_i} \otimes
D_i))$ corresponds to the given holomorphic structure on
$E=\bigoplus_{i=1}^q \CC^{m_i} \otimes D_i$. With respect to this
identification, the action of $R = \prod_{i=1}^q GL(m_i;\CC)$ on
${\cal C}$ is the action induced by the obvious action of $R$ on
$\bigoplus_{i=1}^q \CC^{m_i} \otimes D_i$.  The tangent space to
the ${\cal G}_c$-orbit through this holomorphic structure is the
image of the differential
$$\Omega^0 (\End E) \to \Omega^{0,1}(\End E),$$
and  the normal ${\cal N}_{R}$ to ${\cal G}_cZ^{ss}_{R}$ at
$E$  is naturally isomorphic to
the cokernel of the restriction of this differential to
$$\Omega^0 (\End'_\oplus E) \to \Omega^{0,1}(\End'_\oplus E),$$
which is
\[
H^1(\Sigma,\mbox{End}'_{\oplus} E) \cong \bigoplus_{i_1,i_2=1}^{q}
\CC^{m_{i_1}m_{i_2}-\delta_{i_1}^{ i_2}} \otimes H^1(\Sigma,D^*_{i_1}
\otimes D_{i_2}).
\]
We have $H^0(\Sigma,\mbox{End}'_{\oplus} E) =0$,
so the normal bundle
 to ${\cal G}_cZ^{s}_{R}$ in ${\cal C}$ is
$$-\pi_! \left( \bigoplus_{i_1,i_2=1}^{q}
\CC^{m_{i_1}m_{i_2}-\delta_{i_1}^{ i_2}} \otimes
{\mathbb{U}}^*_{i_1} \otimes {\mathbb{U}}_{i_2} \right)
$$
where ${\mathbb{U}}_{i_1}$ and ${\mathbb{U}}_{i_2}$ are the
appropriate universal bundles on ${\cal C}(n_{i_1},d_{i_1}) \times
\Sigma$ and ${\cal C}(n_{i_2},d_{i_2}) \times \Sigma$ pulled back
to ${\cal C}(n_{i_1},d_{i_1}) \times {\cal C}(n_{i_2},d_{i_2})
\times \Sigma$, and $\pi$ denotes the projection from ${\cal
C}(n_{i_1},d_{i_1}) \times {\cal C}(n_{i_2},d_{i_2}) \times
\Sigma$ to ${\cal C}(n_{i_1},d_{i_1}) \times {\cal
C}(n_{i_2},d_{i_2})$.

We need to extend this description to a description of the normal bundle
${\cal N}_R$
 to the proper transform ${\cal G}_c \hat{Z}^{ss}_{R}$ of ${\cal G}_c
{Z}^{ss}_{R}$. For simplicity we
consider the case when $ \hat{Z}^{ss}_{R}$ is obtained from
${Z}^{ss}_{R}$ via a single blow up
along its intersection with ${\cal G}_c {Z}^{ss}_{R'}$ for some $R'$
containing $R$; the general case
can then be obtained inductively (cf. \cite{K4} \S 8). From \cite{K4}
Corollary 8.11 we know that
${\cal G}_c {Z}^{ss}_{R'} \cap Z_R^{ss}$ is the disjoint union
\begin{equation} \label{disjoi} {\cal G}_c {Z}^{ss}_{R'} \cap Z_R^{ss} =
\bigsqcup_{r=1}^{m'} N^R_0 g'_r Z_{R'}^{ss} \end{equation} where
the subset $\{g \in {\cal G}_c: R \subseteq gR'g^{-1} \}$ of
${\cal G}_c$ is the disjoint union
$$\{g \in {\cal G}_c: R \subseteq gR'g^{-1} \}  =
\bigsqcup_{r=1}^{m'} N^R_0 g'_r N^{R'} $$ of finitely many double
$(N_0^R,N^{R'})$ cosets\footnote{Strictly speaking here we should
apply \cite{K4} Corollary 8.11 to a finite dimensional description
of $\mnd$ as a quotient (cf. $\S$3 above and [37]).}. Moreover if $E$ represents
an element of one of the components $N_0^R g Z_{R'}^{ss}$ of
${\cal G}_c {Z}^{ss}_{R'} \cap Z_R^{ss}$ where $R \subset gR'
g^{-1}$ then we have
\begin{equation} \label{canonical} E = (\CC^{m'_1} \otimes D'_1) \oplus
\cdots \oplus (\CC^{m'_Q} \otimes D'_Q)
\end{equation}
with $D'_1, \ldots, D'_Q$ all stable and not isomorphic to each other
and $D'_j$ of rank $n'_j$ and
degree $d'_j$ where
$$R' \cong \prod_{j=1}^Q GL(m'_j;\CC).$$
But also since $E \in Z_R^{ss}$ we have a decomposition
$$ E \cong (\CC^{m_1} \otimes D_1) \oplus \cdots \oplus (\CC^{m_q}
\otimes D_q)
$$
where $D_i$ is semistable of rank $n_i$ and degree $d_i$. As the
decomposition (\ref{canonical})
of $E$ is canonical, for $1 \leq i \leq q$ we must have
$$D_i \cong (\CC^{M_{i,1}} \otimes D'_1) \oplus \cdots \oplus
(\CC^{M_{i,Q}} \otimes D'_Q) $$
for some $M_{ij} \geq 0$ satisfying
\begin{equation} \label{mdash} m'_j = \sum_{i=1}^q m_i M_{ij}
\end{equation}
for $1 \leq j \leq Q$, and we can assume that $R \cong \prod_{i=1}^q
GL(m_i;\CC)$ is
embedded in $gR' g^{-1} \cong \prod_{j=1}^Q GL(m'_j;\CC)$ via
decompositions
$$\CC^{m'_j} \cong \bigoplus_{i=1}^q \CC^{m_i} \otimes \CC^{M_{ij}}$$
coming from (\ref{mdash}). Note that the normal to $N^R_0 g'_r
Z_{R'}^{ss}$ in $Z_R^{ss}$ is then
\begin{equation} \label{hort}
-\pi_! \left(\frac{ \bigoplus_{i=1}^q \bigoplus_{j_1,j_2=1}^{Q}
\CC^{M_{i,j_1}M_{i,j_2}} \otimes ({\mathbb{U}}'_{j_1})^* \otimes
{\mathbb{U}}'_{j_2} }{\bigoplus_{j=1}^Q ({\mathbb{U}}'_{j})^*
\otimes {\mathbb{U}}'_{j} }\right)
\end{equation}
where ${\mathbb{U}}'_j$ denotes the pullback of the appropriate
universal bundle on ${\cal C}(n'_j,d'_j) \times \Sigma$.

\begin{lemma} \label{Fulton15.4} There are short exact sequences of
sheaves over $\hat{Z}_R^{ss}$ as follows,
where ${\cal N}_{A|B}$ denotes the normal bundle to $A$ in $B$ when $A$
is a smooth submanifold of a manifold $B$:
$$ (i) \hspace{.5in} 0 \to {\cal N}_{\hat{Z}_R^{ss}|\hat{{\cal C}}^{ss}}
\to p^*( {\cal N}_{{Z}_R^{ss}|{{\cal C}}^{ss}})
\to p^* \iota_! \left(\frac{T{\cal C}}{T{\cal G}_c Z_{R'}^{ss} +
TZ_R^{ss}} \right) \to 0$$
where $p: \hat{Z}_R^{ss} \to Z_R^{ss}$ is the blow down map and $\iota:
Z_R^{ss} \cap {\cal G}_c Z_{R'}^{ss} \to Z_R^{ss}$
is the inclusion, and
$$\iota_! \left(\frac{T{\cal C}}{T{\cal G}_c Z_{R'}^{ss} + TZ_R^{ss}}
\right)$$
is shorthand for the extension by zero over $Z_R^{ss}$ of the vector
bundle
$$\frac{T{\cal C}|_{Z_R^{ss} \cap {\cal G}_c Z_{R'}^{ss}}
}{T{\cal G}_c Z_{R'}^{ss}|_{Z_R^{ss} \cap {\cal G}_c Z_{R'}^{ss}} +
TZ_R^{ss}|_{Z_R^{ss} \cap {\cal G}_c Z_{R'}^{ss}}} $$
which fits into an exact sequence
$$0 \to T(N_0^R g'_i Z^{ss}_{R'}) \to TZ_R^{ss}|_{ N_0^R g'_i
Z_{R'}^{ss}} \to {\cal N}_{R'}
\to \frac{T{\cal C}|_{Z_R^{ss} \cap {\cal G}_c Z_{R'}^{ss}}
}{T{\cal G}_c Z_{R'}^{ss}|_{Z_R^{ss} \cap {\cal G}_c Z_{R'}^{ss}} +
TZ_R^{ss}|_{Z_R^{ss} \cap {\cal G}_c Z_{R'}^{ss}}} \to 0 $$
on each component $N_0^R g'_i Z_{R'}^{ss}$ of $Z_R^{ss} \cap {\cal G}_c
Z_{R'}^{ss}$;

$$ (ii) \hspace{.5in} 0 \to \Lie({\cal G}_c)/\Lie(N^R) \to {\cal
N}_{\hat{Z}_R^{ss}|\hat{{\cal C}}^{ss}}
\to {\cal N}_{{\cal G}_c \hat{Z}_R^{ss}|\hat{{\cal
C}}^{ss}}|_{\hat{Z}_R^{ss}}
= {\cal N}_R|_{\hat{Z}_R^{ss}} \to 0$$
where $\Lie({\cal G}_c)/\Lie(N^R)$ is shorthand for the trivial vector
bundle on $\hat{Z}_R^{ss}$ whose fiber is
$$\Lie({\cal G}_c)/\Lie(N^R) \cong \Omega^0(\End'_\oplus E).$$
\end{lemma}

\noindent{{\bf Proof:}} (i) follows directly from \cite{F} Lemma 15.4
(i) and (iv), while (ii) is an immediate
consequence of the fact that
\begin{equation}\label{which} {\cal G}_c \hat{Z}_R^{ss} \cong {\cal G}_c
\times_{N^R} \hat{Z}_R^{ss} \end{equation} by \cite{K4} Corollary
5.6.\hfill $\Box$

\begin{corollary} \label{cor6.3}
The equivariant Chern polynomial of the normal bundle ${\cal N}_R$ to
${\cal G}_c \hat{Z}_R^{ss} $ in
$\hat{{\cal C}}^{ss}$ is given by
$$c({\cal N}_R)(t) = p^* \left( \frac{ c( -\pi_! \left(
\bigoplus_{i_1,i_2=1}^{q} \CC^{m_{i_1}m_{i_2}-\delta_{i_1}^{ i_2}}
\otimes {\mathbb{U}}^*_{i_1} \otimes {\mathbb{U}}_{i_2} \right)
)(t)} {c( \iota_! \left(\frac{T{\cal C}}{T{\cal G}_c Z_{R'}^{ss} +
TZ_R^{ss}} \right) )(t)} \right)
$$
with
$$c( \iota_! \left(\frac{T{\cal C}}{T{\cal G}_c Z_{R'}^{ss} + TZ_R^{ss}}
\right) )(t) = \prod_{r=1}^{m'} \left( \frac{ c\left( \iota_!
\left( -\pi_! \left( \bigoplus_{j_1,j_2 = 1}^Q
\CC^{m'_{j_1}m'_{j_2} } \otimes ({\mathbb{U}}'_{j_1})^{*} \otimes
{\mathbb{U}}'_{j_2} \right) \right) \right) (t) } {
 c\left( \iota_! \left( -\pi_! \left( \bigoplus_{i=1}^q
\bigoplus_{j_1,j_2 = 1}^Q \CC^{M_{i,j_1}M_{i,j_2} } \otimes
({\mathbb{U}}'_{j_1})^{*} \otimes {\mathbb{U}}'_{j_2} \right)
\right) \right) (t) } \right)$$ in the notation of Lemma
\ref{Fulton15.4} and the preceding paragraph.
\end{corollary}

\noindent{\bf Proof:} Since $H^*_{{\cal G}}({\cal G}_c \hat{Z}_R^{ss})
\cong H^*_{N^R}(\hat{Z}_R^{ss})$ by
(\ref{which}), it suffices to consider the restriction of ${\cal N}_R$
to $\hat{Z}_R^{ss}$.
 When $E \cong
  (\CC^{m_1} \otimes D_1) \oplus \cdots \oplus (\CC^{m_q} \otimes D_q)
$ represents an element of $Z_R^{ss}$ as above the
normal to $Z_R^{ss}$ in ${\cal C}^{ss}$ at $E$ is naturally isomorphic
to
$\Omega^{0,1}(\End'_\oplus E)$.
The result now follows from (\ref{hort}),
Lemma \ref{Fulton15.4} and the exact sequence
$$0 \to H^0(\Sigma,\End'_\oplus E) \to \Omega^0 (\End'_\oplus E) \to
\Omega^{0,1}(\End'_\oplus E) \to H^1(\Sigma, \End'_\oplus E) \to
0.$$\hfill $\Box$

\bigskip

The one-parameter subgroup $T_1$
of $R$
generated by $\beta$ acts diagonally on $\CC^{m_1} \oplus ... \oplus
\CC^{m_q}$ with
weights $\beta.e_j$ for $j \in \{1,...,M\}$ where $M = m_1 +...+m_q$,
and so it acts on
$$\Omega^{0,1}(\mbox{End}(\bigoplus_{i=1}^q \CC^{m_i} \otimes D_i))
= \bigoplus_{i_1,i_2=1}^q \Omega^{0,1} (\CC^{m_{i_1}} \otimes
(\CC^{m_{i_2}})^* \otimes D_{i_1}
\otimes D_{i_2}^*)$$
with weights $\beta.(e_i - e_j)$ for $i,j \in \{1,...,M\}$.
By Proposition \ref{5.1} we have a partition
$$\{ \Delta_{m}:m \in J = \{ 1, \ldots, t\}  \}$$
of $\{1,...,M \}$ such
that
$$\frac{\beta}{|\!| \beta |\!|^2} = \sum_{m =1}^t
\sum_{j \in \Delta_{m}} (\epsilon - m) \frac{e_j}{|\!| e_j |\!|^2},$$
and
$\beta. (e_i - e_j) = |\!| \beta |\!|^2 $ if and only if  $i \in
\Delta_{k}$ and
$j \in \Delta_{k+1}$ for some $k \in \{1, \ldots, t-1 \}$ by \cite{K8}
Lemma 5.3.
Recall that if $1 \leq i \leq q$ then
$e_{m_1 +...+m_{i-1} +1},...,e_{m_1 + ... + m_i}$ are the weights of the
standard
representation on $\CC^{m_i}$ of the component $GL(m_i;\CC)$ of
$R = \prod_{i=1}^q GL(m_i;\CC)$. If  $1 \leq i \leq q$ and $k \in J =
\{1, \ldots t \}$, then
let
\begin{equation} \label{deltaki}  \Delta^k_i = \Delta_{k} \cap \{m_1 +
... + m_{i-1} + 1, ..., m_1 + ... + m_i \}
\end{equation}
and let $m_i^k$ denote the size of $\Delta^k_i$, so that $\sum_{k =1}^t
m^k_i = m_i$.
If we make the induced identifications
$$ \CC^M = \bigoplus_{i=1}^q \CC^{m_i} = \bigoplus_{i=1}^q \bigoplus_{k
=1}^t \CC^{m_i^k}$$
 then the $T_1$-eigenbundle
${\cal W} = {\cal W}_{i,j}\to N^{T_1}_0 \hat{Z}^{ss}_{R}$
 of the restriction  to $N_0^{T_1} \hat{Z}_{R}^{ss}$
of the normal bundle to  ${\cal G}_c \hat{Z}^{ss}_R$ which corresponds
to the affine hyperplane $\beta + \beta^{\perp}$
is given by the image of
$$-\pi_! \left( \bigoplus_{i_1,i_2=1}^q \bigoplus_{k=1}^{t-1}
 (\CC^{m^{k}_{i_1}} \otimes (\CC^{m^{k+1}_{i_2}})^* \otimes {\mathbb{U}}_{i_1}
\otimes {\mathbb{U}}_{i_2}^*) \right)$$ and has equivariant Chern
polynomial equal to the pullback of
\begin{equation}\label{fillin}
c({\cal W})(t) = p^* \left( \frac{ c\left( -\pi_! \left(
\bigoplus_{i_1,i_2=1}^{q} \bigoplus_{k=1}^{t-1} \CC^{m_{i_1}^k
m_{i_2}^{k+1} } \otimes {\mathbb{U}}_{i_1} \otimes
{\mathbb{U}}_{i_2}^* \right) \right)(t)}{\prod_{r=1}^{m'}
 c\left( \iota_!
\left( -\pi_! \left( \bigoplus_{j_1,j_2 = 1}^Q
\bigoplus_{k=1}^{t-1} \CC^{m^{'k}_{j_1}m^{'k+1}_{j_2}} \otimes
{\mathbb{U}}'_{j_1} \otimes ({\mathbb{U}}'_{j_2})^* \right)
\right) \right) (t)} \right)
\end{equation}
where $m^{'k}_j=\sum_{i=1}^qm_i^kM_{i,j}$.

\begin{rem}
Recall that the final ingredient of the wall crossing term (\ref{*}) which is needed
is the equivariant Euler
class $e_{\PP {\cal W}}$ of the normal bundle
in  $\hat{{\cal C}}^{ss}$ to the projectivisation $\PP {\cal W}$ of
${\cal W} $. This can of course be obtained from the
equivariant Chern polynomial
which we are now in a position to calculate as follows.

The projective bundle
$\PP {\cal W} \to N^{T_1}_0 \hat{Z}^{ss}_{R}$ is a subbundle of the
restriction to $N^{T_1}_0 \hat{Z}^{ss}_{R}$
of the exceptional
divisor $\PP {\cal N}_R$ for the blow-up along ${\cal G}_c
\hat{Z}^{ss}_{R}$. The normal to
$\PP {\cal N}_R$ is ${\cal O}(-1)$, and by \cite{F} Lemma 15.4(ii) there
is an exact
sequence
$$0 \to {\cal O} \to p^*{\cal N}_R \otimes {\cal O}(1) \to T_{\PP {\cal
N}_R} \to p^* T_{{\cal G}_c \hat{Z}_R^{ss}}
\to 0$$
where $T_X$ denotes the tangent bundle to $X$ and $p:\PP {\cal N}_R \to
{\cal G}_c \hat{Z}^{ss}_{R}$
is the natural projection. Similarly we have
$$0 \to {\cal O} \to p^*{\cal W} \otimes {\cal O}(1) \to T_{\PP {\cal
W}} \to (p|_{\PP {\cal W}})^* T_{N^{T_1}_0 \hat{Z}_R^{ss}}
\to 0$$
so the equivariant Chern polynomial of the normal to $\PP {\cal W}$ in
$\PP {\cal N}_R$ is
\begin{equation} \label{e:eqchernnormal}
c( p^*({\cal N}_{N^{T_1}_0 \hat{Z}_R^{ss}|{\cal G}_c
\hat{Z}^{ss}_R})) (t) c(p^*({\mathcal{N}_R}/{\cal W}) \otimes
{\cal O}(1))(t).
\end{equation}
where
${\cal N}_{N^{T_1}_0 \hat{Z}_R^{ss}|{\cal G}_c
\hat{Z}^{ss}_R}$ is the
normal bundle to
$N^{T_1}_0 \hat{Z}_R^{ss}$ in ${\cal G}_c \hat{Z}_R^{ss}$.
We have
$${\cal G}_c \hat{Z}_R^{ss} \cong {\cal G}_c \times_{N^R}
\hat{Z}_R^{ss}$$
and
$$N_0^{T_1} \hat{Z}_R^{ss} \cong N_0^{T_1} \times_{N_0^{T_1} \cap N^R}
\hat{Z}_R^{ss},$$
so the normal  to
$N^{T_1}_0 \hat{Z}_R^{ss}$ in ${\cal G}_c \hat{Z}_R^{ss}$ is isomorphic
to  ${\cal G}_c / N_0^{T_1} N^R$,
and therefore
its equivariant Chern roots are
the weights of the natural action on ${\cal G}_c / N_0^{T_1} N^R$.

Let $\lambda_1, \ldots, \lambda_l$ be the equivariant Chern roots of
${\cal N}_R/{\cal W}$, so that
its equivariant Chern polynomial is given by
$$c\left( {\cal N}_R/{\cal W}\right)(t) = \prod_{j=1}^l (1 + \lambda_j t
).$$
Then the equivariant Chern polynomial of $p^*({\cal N}_R/{\cal
W})\otimes {\cal O}(1)$
is
\begin{equation} \label{e:eqchern}
 \prod_{j=1}^l (1 + \lambda_j t + \zeta t) = (1 + \zeta t)^l c\left(
{\cal N}_R/{\cal W}\right)\Bigl ( \frac{t}{1 + \zeta t}\Bigr ) \end{equation}
where $\zeta$ is the standard generator of the cohomology of the
projective bundle over the
cohomology of the base.
But by Corollary \ref{cor6.3} we have
$$c({\cal N}_R)(t) = p^* \left( \frac{ c( -\pi_! \left(
\bigoplus_{i_1,i_2=1}^{q} \CC^{m_{i_1}m_{i_2}-\delta_{i_1}^{ i_2}}
\otimes {\mathbb{U}}^*_{i_1} \otimes {\mathbb{U}}_{i_2} \right)
)(t)} {c( \iota_! \left(\frac{T{\cal C}}{T{\cal G}_c Z_{R'}^{ss} +
TZ_R^{ss}} \right) )(t)} \right)
$$
with
$$c( \iota_! \left(\frac{T{\cal C}}{T{\cal G}_c Z_{R'}^{ss} + TZ_R^{ss}}
\right) )(t)
= $$
$$ \prod_{r=1}^{m'} \left( \frac{
c\left( \iota_! \left( -\pi_! \left( \bigoplus_{j_1,j_2 = 1}^Q
\CC^{m'_{j_1}m'_{j_2}} \otimes ({\mathbb{U}}'_{j_1})^{*} \otimes
{\mathbb{U}}'_{j_2} \right) \right) \right) (t) } {
 c\left( \iota_! \left( -\pi_! \left( \bigoplus_{i=1}^q
\bigoplus_{j_1,j_2 = 1}^Q \CC^{M_{i,j_1}M_{i,j_2}} \otimes
({\mathbb{U}}'_{j_1})^{*} \otimes {\mathbb{U}}'_{j_2} \right)
\right) \right) (t) } \right) .$$ Moreover by (\ref{fillin}) we
have
\begin{equation}
c({\cal W})(t) = p^* \left( \frac{ c\left( -\pi_! \left(
\bigoplus_{i_1,i_2=1}^{q} \bigoplus_{k=1}^{t-1} \CC^{m_{i_1}^k
m_{i_2}^{k+1} } \otimes {\mathbb{U}}_{i_1} \otimes
{\mathbb{U}}_{i_2}^* \right) \right)(t)}{\prod_{r=1}^{m'}
 c\left( \iota_!
\left( -\pi_! \left( \bigoplus_{j_1,j_2 = 1}^Q
\bigoplus_{k=1}^{t-1} \CC^{m^{'k}_{j_1}m^{'k+1}_{j_2}} \otimes
{\mathbb{U}}'_{j_1} \otimes ({\mathbb{U}}'_{j_2})^* \right)
\right) \right) (t)} \right)
\end{equation}
so
\begin{equation} \label{e:final}  c\left( {\cal N}_R/{\cal W}\right)(t) =
\frac{c\left( {\cal N}_R\right )(t)}{c\left (\cal W\right)(t)} =
\frac{p^* \Bigl (c(- \pi_! V_0)(t)c(\iota_!(- \pi_!)
V'_0)(t)\Bigr)} {p^*\Bigl ( c(-\pi_! V_1) (t)c(\iota_!(- \pi_!)
V'_1)(t) \Bigr)}
\end{equation}
where
\begin{equation} V_0 =
\bigoplus_{i_1,i_2=1}^{q} \CC^{m_{i_1}m_{i_2}-\delta_{i_1}^{ i_2}}
\otimes {\mathbb{U}}^*_{i_1} \otimes {\mathbb{U}}_{i_2}
\end{equation}
and
\begin{equation} V_1 =
\bigoplus_{i_1,i_2=1}^{q} \bigoplus_{k=1}^{t-1} \CC^{m_{i_1}^k
m_{i_2}^k } \otimes {\mathbb{U}}_{i_1} \otimes
{\mathbb{U}}^*_{i_2}
\end{equation}
while
$$ V'_0 =  \bigoplus_{r=1}^{m'}
\left( \bigoplus_{i=1}^q \bigoplus_{j_1,j_2 = 1}^Q
\CC^{M_{i,j_1}M_{i,j_2}} \otimes ({\mathbb{U}}'_{j_1})^{*} \otimes
{\mathbb{U}}'_{j_2} \oplus \bigoplus_{j_1,j_2 = 1}^Q
\bigoplus_{k=1}^{t-1} \CC^{m^{'k}_{j_1}m^{'k+1}_{j_2}} \otimes
{\mathbb{U}}'_{j_1} \otimes ({\mathbb{U}}'_{j_2})^* \right)
$$
and
$$ V'_1 = \bigoplus_{r=1}^{m'} \bigoplus_{j_1,j_2 = 1}^Q
\CC^{m'_{j_1}m'_{j_2}} \otimes ({\mathbb{U}}'_{j_1})^{*} \otimes
{\mathbb{U}}'_{j_2}.$$

By the Grothendieck--Riemann--Roch theorem
(\cite{F}, Theorem 15.2) we have
\begin{equation} \label{e:grr} {\rm ch} (f_! \alpha) {\rm td} (T_{Y})
= f_* \Bigl ( {\rm ch (\alpha)} {\rm td} (T_{X}) \Bigr )
\end{equation}
when
   $f: X \to Y$ is proper.
When $f=\pi$ is a fibration with fiber $\Sigma$ this gives us
\begin{equation} \label{e:grr'}
{\rm ch} \Bigl (\pi_! \alpha\Bigr )
= \pi_* \Bigl( {\rm ch (\alpha)} {\rm td} (\Sigma) \Bigr  )
= \pi_* \Bigl( {\rm ch (\alpha)} (1-(g-1)\omega) \Bigr  )
 \end{equation}
where $\omega$ is the standard generator of $H^2(\Sigma)$.
When $f=\iota:
Z_R^{ss} \cap {\cal G}_c Z_{R'}^{ss} \to Z_R^{ss}$
is a closed embedding we get
\begin{equation} \label{e:grr''}
{\rm ch} \Bigl (\iota_! \alpha\Bigr ) = \iota_* \Bigl( {\rm ch
(\alpha)} \left({\rm td} \left(\mathcal{N}_{Z^{ss}_R\cap
\mathcal{G}_cZ^{ss}_{R'} | Z^{ss}_R} \right)\right)^{-1} \Bigr  )
 \end{equation}
(see \cite{F} $\S$15.2; in particular the formula immediately
before Corollary 15.2.1). Then we use (\ref{disjoi}) and
(\ref{hort}) to compute the Chern polynomials.

By putting all this together we can (at least in principle) calculate the equivariant
Chern polynomial of the normal bundle to $\PP {\cal W}$,
and hence  calculate the equivariant Euler class $e_{\PP
{\cal W}}$ (cf. similar
calculations in \cite{EK}, in particular \cite{EK} Proposition 9.2).

\end{rem}

We now have all the ingredients
needed to calculate the wall crossing terms
(\ref{*})
by integrating over the fibers  ${\PP ({\cal W}_{i,j})_x} \quott Stab
(x)\cap N_0^{T_1} \cap N^{R_i}$ of the map
$$ \Psi:
{\PP {\cal W}}_{i,j}  \quott_{\xi_1} (N_0^{T_1}/(T_1)_c) \cong
{\PP {\cal W}}_{i,j}  \mid_{\hat{Z}_{R_i}^{ss} } \quott N_0^{T_1} \cap
N^{R_i}
\longrightarrow \hat{Z}_{R_i} \quott N_0^{T_1} \cap N^{R_i}
$$
defined at (\ref{commutdiag}). The results can be expressed (as in
 \cite{jkkw}) in terms of integrals
over projective subbundles of ${\PP ({\cal W}_{i,j})}$ which can be
calculated with
the following standard lemma.

\begin{lemma} Let $E$ be a rank $r$ complex vector bundle over
a manifold $M$ and let $\eta \in H^*(\PP E) \cong
H^*(M)[y]/(p(y))$ where $p(y) =  y^r + c_1(E) y^{r-1} + \dots +
c_r(E) $. Then
$$ \int_{{\PP} E} \eta = {\rm res}_{y = 0 } \int_M \frac{\eta}{p(y)}.
$$
\end{lemma}

This lemma reduces the wall crossing terms (\ref{*}) to integrals over
spaces of the form
$$  \hat{Z}_R \quott N^R \cong \left[ \prod_{i=1}^q \tilde{{\cal
M}}(n_i,d_i)
\right] / \pi_0(N^R) $$
(see (\ref{prodmod})) which can be calculated using induction on $n$.

In the next
section we will carry out the details of all these calculations in the
case when $n=2$.

\renorm
\section{The case when the rank $n$ is two}

In this section we explicitly compute intersection numbers in the
partial desingularization $\widetilde{\calm} (2,0)$. Let
$M^{\rm ext}=M^{\rm ext}_{U(2)}$ be the extended moduli space with the
group $U(2)$ on which $K=SU(2)$ acts by conjugation, i.e.
$M^{\rm ext}$ is the fiber product
 \begin{equation}
   \begin{CD}
   M^{\rm ext}    @>{\mu}>> {\mathbf{su} (2)}  \\
   @VVV             @VVV \exp   \\
   U(2)^{2g}         @>{\Phi}>>  SU(2)
   \end{CD}
\end{equation}

Let $T$ denote the maximal torus of $K=SU(2)$ consisting of
diagonal matrices. Fix a sufficiently small positive number
$\varepsilon$. From \eqref{tor} we have

\begin{equation} \label{e:partint}
\begin{array}{lll}
\kappa_{\widetilde{M}^{\rm ext}}\big(\eta e^{\bar{\omega}}\big)
[\widetilde{\calm}(2,0)] &=& \frac12 \left (
-\int_{\tilde{\mu}_T^{-1}(0)/T }\kappa_{\widetilde{M}^{\rm
ext}}^T\big( \eta e^{\bar{\omega}} \cald^2\big) +
\int_{\tilde{\mu}_T^{-1}(\varepsilon)/T
}\kappa_{\widetilde{M}^{\rm ext}}^{T,\varepsilon}\big( \eta
e^{\bar{\omega}} \cald^2\big) \right )\\ &-& \frac12
\int_{\mu_T^{-1}(\varepsilon)/T } \kappa_{M^{\rm
ext}}^{T,\varepsilon}\big(\eta  e^{\bar{\omega} } \cald^2\big)
\end{array}\end{equation}
where $\mu_T$ and $\tilde{\mu}_T$ are moment maps for the
$T$-action on $M^{\rm ext}$ and $\widetilde{M}^{\rm ext}$ respectively.
Notice that the second term in (\ref{e:partint}) may be computed
using periodicity as in \cite{JK2}, as explained in $\S$4 of
the present paper.

To compute the first term, we need to examine the walls (images
under $\tilde{\mu}_T$ of components of the fixed point set of $T$)
crossed in passing from $0$ to $\varepsilon$ in $\liet^*=\RR$. The
components of the fixed point set that are relevant to us are
those that meet the exceptional divisors
 in $\widetilde{M}^{\rm ext}$
(see \cite{jkkw}, Lemma 23).

\subsection{Wall crossing term from the first blow-up}
To form $\widetilde{M}^{\rm ext}$ from $M^{\rm ext}$, we first blow up
along the set $\Delta$ of the points with stabilizer $K = SU(2)$
which is $(S^1)^{2g}$ where $S^1$ represents the center of $U(2)$.
With the natural identification of $(S^1)^{2g}$ with the Jacobian
${\rm Jac}$ of $\Sigma$, the normal bundle of the $K$-fixed point locus
$\Delta$ is
\[
  R^1\pi_* \End (\mathcal{L}\oplus \mathcal{L})_0\cong
  R^1\pi_*\mathcal{O}\otimes {\mathbf{sl} (2)}
\]
where $\mathcal{L}\to \Sigma\times {\rm Jac}$ is the Poincar\'e
bundle and $\pi:\Sigma\times {\rm Jac}\to {\rm Jac}$ is the
projection. Hence, the exceptional divisor of the first blow-up in
the partial desingularisation process is $\PP
(R^1\pi_*\mathcal{O}\otimes {\mathbf{sl} (2)})$. The $T$-fixed
point component with positive moment map value in the exceptional
divisor is thus the projectivization
\[
\PP\left[R^1\pi_*\mathcal{O}\otimes
\left(\begin{array}{lll}0&1\\
0&0\end{array}\right)\right]\cong \PP T_{\rm Jac}\cong {\rm Jac}\times
\PP^{g-1}
\]
of the tangent bundle of the Jacobian which is trivial. The normal
bundle to this fixed point component is thus
\[
\mathcal{O}_{\PP^{g-1}}(-1)\oplus \left[R^1\pi_*\mathcal{O}\otimes
\left(\begin{array}{lll}1&0\\
0&-1\end{array}\right)\otimes \mathcal{O}_{\PP^{g-1}}(1) \right]
\oplus \left[R^1\pi_*\mathcal{O}\otimes
\left(\begin{array}{lll}0&0\\
1&0\end{array}\right)\otimes \mathcal{O}_{\PP^{g-1}}(1) \right].
\]
The first summand is normal to the exceptional divisor and the
last two are normal in the exceptional divisor. The torus $T$ acts
on the three summands with weights $2,-2, -4$ respectively. As
$R^1\pi_*\mathcal{O}\cong T_{\rm Jac}$ is trivial, the equivariant
Euler class of the normal bundle is
\[
(-y+2Y)(y-2Y)^g(y-4Y)^g
\]
where $y=c_1(\mathcal{O}_{\PP^{g-1}}(1))$ is the generator of
$H^*(\PP^{g-1})$ and $Y$ is the generator of $H^*_T(pt)$. Hence
the wall crossing term from the first blow-up in the partial
desingularization process is $$ -\frac12\int_{\PP^{g-1}\times {\rm Jac}}
\res_{Y=0} \frac{\eta e^{\bar{\omega}}|_\Delta
\cald^2}{(-y+2Y)(y-2Y)^g(y-4Y)^g} $$
\begin{equation}\label{yh-8.1}=\frac12\int_{\Delta}\int_{\PP^{g-1}}\res_{Y=0}\frac{\eta
e^{\bar{\omega}}|_\Delta (4Y^2)
}{(-2)^{g+1}(-4)^gY^{2g+1}}(1-\frac{y}{2Y})^{-g-1}(1-\frac{y}{4Y})^{-g}\end{equation}
$$=-\frac1{2^{3g}}\int_{\Delta}\res_{y=0}\res_{Y=0}\frac{\eta
e^{\bar{\omega}}|_{\Delta}}{y^gY^{2g-1}}
(1-\frac{y}{2Y})^{-g-1}(1-\frac{y}{4Y})^{-g}.$$

A universal bundle $\mathcal{U}\to \Sigma\times \mathcal{C}^{ss}$
restricted to $\Sigma\times Z^{ss}_{SU(2)}$ where $Z^{ss}_{SU(2)}$
is the locus of $SU(2)$-fixed points,  decomposes as
$\mathcal{L}\oplus \mathcal{L}$ for a universal line bundle
$\mathcal{L}$ over $\Sigma\times Z^{ss}_{SU(2)}$ by abuse of
notation. Since we are considering the moduli space of degree zero
semistable bundles, the equivariant first Chern class is of the
form
\[
c_1(\mathcal{L})=Y\otimes 1+\sum d^j\otimes \sigma_j
\]
where $d^j\in H^1(\Delta)$  and $\sigma_j$ is a symplectic basis
of $H^1(\Sigma)$ so that $\sigma_j\sigma_{j+g}$ is the fundamental
class $\rho\in H^2(\Sigma)$. Hence $c_2(\mathcal{U})$ restricts as
\[
c_2(\mathcal{L}\oplus \mathcal{L})=c_1(\mathcal{L})^2=Y^2\otimes
1+\sum Yd^j\otimes \sigma_j-2\gamma\otimes \rho
\]
where $\gamma=\sum_{j=1}^g d^jd^{j+g}$. This implies that the
classes $a_2, b_2^j, f_2$ defined in $\S$4 restrict to $Y^2$,
$Yd^j$ and $-2\gamma$ respectively. Therefore, we can explicitly
compute the restriction of $\eta e^{\bar{\omega}}$ to $\Delta$ in
terms of $Y$ and $d^j$.

The last expression in (\ref{yh-8.1}) is nonzero only when $\eta
e^{\bar{\omega}}|_\Delta\in H^*_T(\Delta)=H^*_T({\rm Jac})$ is a
constant multiple of the product of the fundamental class
$\frac{\gamma^g}{g!}$ of $\Delta$ and $Y^{3g-3}$, in which case
the wall crossing term is computed as follows:
$$
-\frac1{2^{3g}}\int_{\Delta}\res_{y=0}\res_{Y=0}\frac{\frac{\gamma^g}{g!}Y^{3g-3}}{y^gY^{2g-1}}
(1-\frac{y}{2Y})^{-g-1}(1-\frac{y}{4Y})^{-g}$$
\begin{equation}\label{yh-8.2} =
-\frac1{2^{3g}}\res_{y=0}\res_{Y=0}\frac{Y^{g-2}}{y^g}
(1-\frac{y}{2Y})^{-g-1}(1-\frac{y}{4Y})^{-g}\end{equation} $$=
-\frac1{2^{3g}}\mathrm{Coeff}_{y^{g-1}Y^{-g+1}}
(1-\frac{y}{2Y})^{-g-1}(1-\frac{y}{4Y})^{-g} $$
$$= -\frac1{2^{3g}} \mathrm{Coeff}_{t^{g-1}}
(1-\frac{t}{2})^{-g-1}(1-\frac{t}{4})^{-g}.$$

\subsection{Wall crossing term from the second blow-up}
Let $\Gamma$ denote the set of points in $M^{\rm ext}$ fixed by the
action of $T$ and $\widehat{\Gamma}$ be the proper transform after
the first blow-up.
 To get the partial desingularization we blow up
along $G\widehat{\Gamma}$ for $G=K^{\CC}$.

Note that $\Gamma=(S^1\times S^1)^{2g}$ where $(S^1\times
S^1)\subset U(2)$ is the set of diagonal matrices in $U(2)$.
$\Gamma$ can be naturally identified with the product ${\rm Jac}\times
{\rm Jac}$ of the Jacobian of $\Sigma$ and $\widehat{\Gamma}$ is the
blow-up of ${\rm Jac}\times {\rm Jac}$ along the diagonal $\Delta\cong {\rm Jac}$.
The normal bundle $\mathcal{N}$ to $G\widehat{\Gamma}$ splits as
the direct sum $W_+\oplus W_-$ on which $T$ acts with weights $2$
and $-2$ respectively. Hence the wall to be crossed is exactly
$\PP W_+$.

Let us compute the equivariant Euler class of the normal bundle of
the wall $\PP W_+$. Let $\mathcal{L}_1$ (resp. $\mathcal{L}_2$) be
the pull-back of the Poincar\'e bundle $\mathcal{L}\to \Sigma
\times {\rm Jac}$ via $\pi_{12}$ (resp. $\pi_{13}$) where
$\pi_{12}:\Sigma\times {\rm Jac}\times {\rm Jac}\to \Sigma\times {\rm Jac}$ (resp.
$\pi_{23}:\Sigma\times {\rm Jac}\times {\rm Jac}\to \Sigma\times {\rm Jac}$) is the
projection onto the first and second (resp. third) components. Let
$\pi_{23}:\Sigma\times {\rm Jac}\times {\rm Jac}\to {\rm Jac}\times {\rm Jac}$ and
$\pi:\Sigma\times {\rm Jac}\to {\rm Jac}$ denote the obvious projections.

As observed in $\S$6, the Chern class of the normal bundle
$\mathcal{N}$ of $G\widehat{\Gamma}$ restricted to
$\widehat{\Gamma}$ is
\[
c(\mathcal{N}|_{\widehat{\Gamma}})=\frac{c\left(
-(\pi_{23})_!(\mathcal{L}_1^\vee\otimes \mathcal{L}_2\oplus
\mathcal{L}_2^\vee\otimes
\mathcal{L}_1)\right)}{c\left(\imath_!(R^1\pi_*(\mathcal{O})\oplus
R^1\pi_*(\mathcal{O})) \right)}
\]
where $\imath:\Delta\hookrightarrow \Gamma$ is the diagonal
embedding. Similarly, we have
\[
c(W_-|_{\widehat{\Gamma}})=\frac{c\left(
-(\pi_{23})_!(\mathcal{L}_1^\vee\otimes
\mathcal{L}_2)\right)}{c\left(\imath_!(R^1\pi_*(\mathcal{O}))
\right)}
\]
\[
c(W_+|_{\widehat{\Gamma}})=\frac{c\left(
-(\pi_{23})_!(\mathcal{L}_2^\vee\otimes
\mathcal{L}_1)\right)}{c\left(\imath_!(R^1\pi_*(\mathcal{O}))
\right)}.
\]

After normalization, the ordinary first Chern classes of
$\mathcal{L}_1$ and $\mathcal{L}_2$ can be written as
$$\sum_{i=1}^{2g} d_1^i\otimes \sigma_i,\ \ \ \ \ \ \sum_{i=1}^{2g} d_2^i\otimes
\sigma_i$$ where $\sigma_i$ is a symplectic basis of $H^1(\Sigma)$
as before. By Grothendieck--Riemann--Roch, the Chern characters are
\[
ch\left(-(\pi_{23})_!(\mathcal{L}_1^\vee\otimes
\mathcal{L}_2)\right)
=ch\left(-(\pi_{23})_!(\mathcal{L}_2^\vee\otimes
\mathcal{L}_1)\right) = (g-1)+\hat{\gamma}
\]
where
\[
\hat{\gamma}=\sum_{i=1}^gd_1^id_1^{i+g}+d_2^id_2^{i+g}+d_1^id_2^{i+g}+d_2^id_1^{i+g}
.\] From a well-known combinatorial argument relating the Chern
characters with the Chern classes, we have
\[
c\left(-(\pi_{23})_!(\mathcal{L}_1^\vee\otimes
\mathcal{L}_2)\right)
=c\left(-(\pi_{23})_!(\mathcal{L}_2^\vee\otimes
\mathcal{L}_1)\right) = \exp(\hat{\gamma})
\]
Since $R^1\pi_*\mathcal{O}\cong \mathcal{O}^{\oplus g}$, we have
\[
c\left(\imath_!(R^1\pi_*(\mathcal{O}))
\right)=c\left(\imath_!(\mathcal{O}) \right)^g=(1+h)^{-g}
\]
by \cite{F} Example 15.3.5 where
$h=c_1(\mathcal{O}_{\hat{\Gamma}}(-E))$ for the exceptional
divisor $E$ of the blow-up $\hat{\Gamma}\to\Gamma$. Therefore we
have
\begin{equation}\label{yh-8.3}
c(W_{\pm}|_{\hat{\Gamma}})=(1+h)^g\exp(\hat{\gamma}).
\end{equation}
It is an illuminating exercise to check that the right hand side
indeed lies in $H^{\le 2g-2}(\hat{\Gamma})$. Let
$y=c_1(\mathcal{O}_{\PP W_+}(1))$. Since the normal bundle to $\PP
W_+$ in $\PP \mathcal{N}$ is the pull-back of $W_-$ tensored with
$\mathcal{O}_{\PP W_+}(1)$ and $T$ acts with weights $-4$, the
equivariant Euler class of the normal bundle to $\PP W_+$ in $\PP
\mathcal{N}$ is
\[
(y-4Y)^{g-1}\left(1+\frac{h}{y-4Y}\right)^g\exp\left(\frac{\hat{\gamma}}{y-4Y}\right).
\]
The normal bundle of $\PP \mathcal{N}$ restricted to $\PP W_+$ is
$\mathcal{O}_{\PP W_+}(-1)$ on which $T$ acts with weight $2$ so
that the equivariant Euler class is $(-y+2Y)$. Therefore the
equivariant Euler class of the normal bundle to $\PP W_+$ in
$\widetilde{M}^{\rm ext}$ is
\begin{equation}\label{yh-8.4}
e_{\PP W_+}=(-y+2Y)(y-4Y)^{g-1} \left(1+\frac{h}{y-4Y}\right)^g
\exp\left(\frac{\hat{\gamma}}{y-4Y}\right).
\end{equation}

As observed in the previous subsection, we can easily compute the
restriction to $G\hat{\Gamma}$ of $\eta e^{\bar{\omega}}$ and
express as a linear combination of classes of the form $\xi Y^n$
for some nonnegative integer $n$ and $\xi\in H^*(\Gamma)\subset
H^*(\hat{\Gamma})$. Since $\dim_\RR \calm (2,0)=8g-6$, the
intersection pairing is nonzero only for classes of degree $8g-6$.
So it suffices to consider the case when $\xi\in
H^{8g-6-2n}(\Gamma)$. The wall crossing term from the second
blow-up is then
\begin{equation}\label{yh-8.5}
\begin{array}{lll}
\int_{\PP W_+} \res_{Y=0}\frac{\xi Y^n (4Y)^2}{e_{\PP W_+}} &=&
-\int_{\PP W_+|_{\hat{\Gamma}}}\res_{Y=0} \frac{\xi
Y^n}{(-y+2Y)(y-4Y)^{g-1} \left(1+\frac{h}{y-4Y}\right)^g
\exp\left(\frac{\hat{\gamma}}{y-4Y}\right)}\\
&=& -\int_{\PP W_+|_{\hat{\Gamma}}}\xi \res_{Y=0}
\frac{Y^{n-g}}{2\cdot (-4)^{g-1}}\sum_{r,s,l\ge
0}A_{r,s,l}\hat{\gamma}^rh^sy^l(-4Y)^{-r-s-l}\\
&=& -\frac1{2\cdot (-4)^n} \int_{\PP W_+|_{\hat{\Gamma}}}\xi
\sum_{r+s+l=n-g+1}A_{r,s,l}\hat{\gamma}^rh^sy^l\\
&=& -\frac1{2\cdot (-4)^n}\sum_{r+s+l=n-g+1}A_{r,s,l} \int_{\PP
W_+|_{\hat{\Gamma}}}\xi \hat{\gamma}^rh^sy^l\end{array}
\end{equation}
where $A_{r,s,l}$ are defined by the power series expansion in $t$
\[
\left((1+2t)\sum_{k=0}^{g-1}(1+t)^{g-1-k}\sum_{r+s=k}\frac1{r!}
\frac{g!}{s!(g-s)!} z^rx^s
\right)^{-1}=\sum_{r,s,l}A_{r,s,l}z^rx^st^l.
\]
So it suffices to compute
\[
\int_{\PP W_+|_{\hat{\Gamma}}}\xi \hat{\gamma}^rh^sy^l
\]
for $r+s+l=n-g+1$. By the residue formula applied to
$W_+|_{\hat{\Gamma}}$, whose Chern class is
$(1+h)^g\exp(\hat{\gamma})$, with respect to the action of $U(1)$
by usual complex multiplication on fibers, we have
\begin{equation}\label{yh-8.6}
\begin{array}{lll}
\int_{\PP W_+|_{\hat{\Gamma}}}\xi \hat{\gamma}^rh^sy^l &=&
\int_{\hat{\Gamma}} \res_{Y=0} \frac{\xi
\hat{\gamma}^rh^sy^l}{y^{g-1}(1+\frac{h}{y})^g
\exp(\frac{\hat{\gamma}}{y})}\\
&=& \int_{\hat{\Gamma}} \res_{Y=0} \frac{\xi
\hat{\gamma}^rh^sy^l}{y^{g-1}}\sum_{a,b\ge
0}B_{a,b}(\frac{\hat{\gamma}}{y})^a(\frac{h}{y})^b\\
&=&\sum_{a,b}\int_{\hat{\Gamma}} \res_{Y=0} B_{a,b}
\frac{\xi\hat{\gamma}^{r+a}h^{s+b}}{y^{g-1-l+a+b}}\\
&=&
\sum_{a+b=l-g+2}B_{a,b}\int_{\hat{\Gamma}}\xi\hat{\gamma}^{r+a}
h^{s+b}\end{array}
\end{equation}
where $B_{a,b}$ are defined by the power series expansion
\[
\frac1{(1+x)^g\exp z}= \sum B_{a,b}z^ax^b.
\]
Since $\hat{\gamma}, \xi\in H^*(\Gamma)$ and $-h$ is the
Poincar\'e dual of the exceptional divisor $E$ in $\hat{\Gamma}$,
it is easy to see that the integral
\[
\int_{\hat{\Gamma}}\xi\hat{\gamma}^{r+a} h^{s+b}
\]
is nonzero only if $s+b=0$ or $s+b=g$.

First suppose $b=-s$ so that $a=n-2g+3-r$ since $r+s+l=n-g+1$ and
$a+b=l-g+2$. Then we have
\[
\int_{\hat{\Gamma}}\xi\hat{\gamma}^{r+a} h^{s+b}=\int_{\Gamma} \xi
\hat{\gamma}^{n-2g+3}.
\]
Since $\xi\in H^{8g-6-2n}(\Gamma)$, $\xi \hat{\gamma}^{n-2g+3}$ is
a constant multiple of $\prod_{j=1}^gd_1^jd_1^{j+g}d_2^jd_2^{j+g}$
in which case we can complete the computation from
\[
\int_{\Gamma}\prod_{j=1}^gd_1^jd_1^{j+g}d_2^jd_2^{j+g}=1
\]

Next suppose $b=g-s$ so that $a=n-3g+3-r$. In this case,
\[
\int_{\hat{\Gamma}}\xi\hat{\gamma}^{r+a}h^{s+b}=
\int_{\hat{\Gamma}}\xi\hat{\gamma}^{n-3g+3}h^{g}= -\int_E\xi
(4\gamma)^{n-3g+3}h^{g-1}=-\int_\Delta \xi (4\gamma)^{n-3g+3}
\]
where $E$ is the exceptional divisor in $\hat{\Gamma}$ which is a
projective bundle over $\Delta$ and $\gamma\in H^2(\Delta)$ is the
class introduced in the previous subsection.

Since $\xi\in H^{8g-6-2n}(\Gamma)$, the integral
\[
\int_{\Delta}\xi (4\gamma)^{n-3g+3}=4^{n-3g+3} \int_{\Delta}\xi
\gamma^{n-3g+3}
\]
is nonzero only when $\xi|_\Delta$ is a constant multiple of
$\gamma^{4g-3-n}$. Hence the computation is now complete from
\[
\int_\Delta \gamma^g=g!
\]

To be quite explicit, let us consider the classes $a_2^mf_2^n$
with $2m+n=4g-3$. Recall that $a_2|_\Delta=Y^2$ and
$f_2|_{\Delta}=-2\gamma$. Similarly, from
$$c_2(\mathcal U)|_\Gamma=c_2(\mathcal L_1\oplus\mathcal
L_2)=c_1(\mathcal L_1)c_1(\mathcal L_2)=Y^2\otimes
1+\sum_{i=1}^{2g} Y(d^i_1+d^i_2)\otimes \sigma_i-(\sum_{j=1}^g
d_1^jd_2^{j+g}+d_2^jd_1^{j+g})\otimes \rho$$ we see that
$a_2|_\Gamma=Y^2$, $f_2|_\Gamma=-\gamma_{12}$ where $\gamma_{12}=
\sum_{j=1}^g d_1^jd_2^{j+g}+d_2^jd_1^{j+g}$. Therefore
\[
a_2^mf_2^n|_\Delta=Y^{2m}(-2\gamma)^n,\ \ \ \ \
a_2^mf_2^n|_\Gamma=Y^{2m}(-\gamma_{12})^n.
\]
Notice that $\gamma_{12}|_{\Delta}=2\gamma$. Let $\gamma_1=
\sum_{j=1}^g d_1^jd_1^{j+g}$ and
$\gamma_2=\sum_{j=1}^gd_2^jd_2^{j+g}$ so that
$\hat{\gamma}=\gamma_1+\gamma_2+\gamma_{12}$.

The wall crossing term from the first blow-up is zero if $g$ is
even so that $3g-3$ is odd or if $n\ne g$. If $g$ is odd and
$n=g$,
$$a_2^mf_2^n|_\Delta=Y^{3g-3}(-2\gamma)^g=(-2)^gg!\frac{\gamma^g}{g!}Y^{3g-3}
=-2^gg!\frac{\gamma^g}{g!}Y^{3g-3}$$ and thus from (\ref{yh-8.2})
the wall crossing term is
\[
\frac{g!}{2^{2g}} \mathrm{Coeff}_{t^{g-1}}
(1-\frac{t}{2})^{-g-1}(1-\frac{t}{4})^{-g} .
\]

From (\ref{yh-8.5}) and (\ref{yh-8.6}), the wall crossing term
from the second blow-up is
\[
\begin{array}{lll}
&-&\frac1{2\cdot (-4)^{2m}}\sum_{r+s+l=2m-g+1}A_{r,s,l}
\sum_{a+b=l-g+2}B_{a,b}\int_{\hat{\Gamma}}(-\gamma_{12})^n\hat{\gamma}^{r+a}
h^{s+b}
\\
&=& - \frac1{2^{4m+1}}\sum_{r+s+l=2m-g+1}A_{r,s,l}
\sum_{a+b=l-g+2}\\
&& \big(B_{2m-2g+3-r,-s} \int_{\Gamma} (-\gamma_{12})^n
\hat{\gamma}^{2m-2g+3} +B_{2m-3g+3-r,g-s}
\int_{\Delta}(-2\gamma)^n (4\gamma)^{2m-3g+3} \big) .
\end{array}
\]
The last integral is just
\[
\int_{\Delta}(-2\gamma)^n
(4\gamma)^{2m-3g+3}=(-1)^n2^{2m-2g+3}\int_\Delta \gamma^g=
(-1)^n2^{2m-2g+3}g!
\]
Also, by combinatorial computation, we have
\[
\int_{\Gamma} (-\gamma_{12})^n \hat{\gamma}^{2m-2g+3}
=(-1)^n\int_{\Gamma}\gamma_{12}^n(\gamma_{12}+\gamma_1+\gamma_2)^{2g-n}
\]
\[
=(-1)^n
\sum_{k=0}^{[\frac{2g-n}{2}]}\frac{(-1)^k(2g-2k)!(2g-n)!g!}{(2g-2k-n)!k!(g-k)!}.
\]
So the computation is complete and we have proved the following.

\begin{theorem}\label{mthmsect8}
For a pair $(m,n)$ of nonnegative integers satisfying $2m+n=4g-3$,
\[
\kappa(a_2^mf_2^n)
[\widetilde{\calm}(2,0)]\quad =\quad
\frac{(-1)^{g-1-m}n!}{2^{2m-g+1}}\res_{Y=0} \frac{1}{Y^{2g-2-2m}
(e^Y-1)}
\]
\[
+ \left(\frac{g!}{2^{2g}}\mathrm{Coeff}_{t^{g-1}}
(1-\frac{t}2)^{-g-1}(1-\frac{t}4)^{-g} \right) \delta_{g,n} -
\frac1{2^{4m+1}}\sum_{r+s+l=2m-g+1}A_{r,s,l} \sum_{a+b=l-g+2}
\]
\[
\left[ (-1)^n
\sum_{k=0}^{[\frac{2g-n}{2}]}\frac{(-1)^k(2g-2k)!(2g-n)!g!}{(2g-2k-n)!k!(g-k)!}
B_{2m-2g+3-r,-s} +(-1)^n2^{2m-2g+3}g!B_{2m-3g+3-r,g-s}  \right]
\]
where $\delta_{g,n}$ is Kronecker's delta and the constants $A_{r,s,l}$ and $B_{a,b}$ are defined as after (\ref{yh-8.5}) and (\ref{yh-8.6}) respectively.
\end{theorem}

\renorm\section{Witten's integrals}

\newcommand{\baromega}{{\bar{\omega}}}

In the case when $0$ is a regular value of the moment map $\mu$
Witten \cite{tdgr} relates the intersection pairings of two
classes $\k_M(\a)$, $\k_M(\b)$ of complementary degrees in
$H^*(M/\!/G)$ coming from $\a$, $\b $ $\in H^*_K(M)$ to the
asymptotic behaviour of the integral $\cali^\epsilon (\a \b
e^{i\baromega})   $ given by \beq \label{iepstwo} \cali^\epsilon
(\a\b  e^{i\baromega}) = \frac{1}{(2 \pi )^s  \vol(K) }  \int_{X
\in \liek} [d X] e^{-\e<X,X>/2} \int_M \a(X) \b(X) e^{i \omega}
e^{i \mu(X)}  \end{equation} \beq \label{iepstwo'} = \frac{1}{(2
\pi)^l |W| \vol( T) } \int_{X \in \liet} [d X] e^{-\e<X,X>/2}
\int_M \a(X) \b(X)\cald^2(X)  e^{i \omega} e^{i \mu(X)},
\end{equation} where as before $\baromega = \omega + \mu$ and
the notations $s$ and $l$ are as after (\ref{res}). $\S$2. He
expresses the integral as a sum of contributions, one of which is
localized near $\mu^{-1}(0)$ and reduces to the intersection
pairing required, while the rest tends to 0 exponentially fast as
$\e$ tends to 0. These results were described in  \cite{jkkw},
$\S$9 and extended there to the case when 0 is not a regular value
of $\mu$.

\begin{rem}
Note that in the  conventions of \S \ref{ss:eqpd}, $ \bar{\omega}
= \omega + \mu $ satisfies $d_K \bar{\omega} = 0$. On the other
hand, Witten uses the convention that \beq \label{e:eq4} d_K = d -
i \iota_{\nu_\xi} \eeq for $\xi \in \liek$ so that in his
convention $\omega + i \mu$ is equivariantly closed. For this
reason we have chosen to retain $e^{i \omega} $ in our integral
(\ref{iepstwo}), whereas Witten writes the integral
(\ref{iepstwo}) without the factor $i$ multiplying $\omega$.
\end{rem}

Even when $0$ is not a regular value of $\mu$, Witten's integral
$\cali^\epsilon (\a\b  e^{i\baromega}) $ decomposes into the sum
of a term $\cali^\epsilon_0(\a\b  e^{i\baromega})$ determined by
the action of $K$ on an arbitrarily small neighbourhood of
$\mu^{-1}(0)$, and other terms which tend to zero exponentially
fast as $\epsilon \to 0 $. Moreover there is a residue formula for
$\cali^\epsilon_0(\a\b  e^{i\baromega})$ which is a sum over
components of the fixed point set of $T$ on $M$ and reduces to the
residue formula (\ref{res}) when 0 is a regular value of $\mu$
(see \cite{jkkw} $\S$9). When 0 is not a regular value of $\mu$
then this residue formula is related to, but not quite the same
as, the formulas for intersection pairings given in previous
sections; it is not in general a polynomial in $\epsilon$ but
instead it is a polynomial in $\sqrt{\e}$ (as was proved by
Paradan in \cite{paradan} Cor.5.2).

\begin{rem} \label{jkkw1rem9.3}
Let $\Theta$ denote the equivariant class given by the invariant
polynomial function $\Theta(X) = \langle X,X \rangle$ on $\liet$.
It is shown in \cite{jkkw} $\S$9 that if a component $F \in \calf$
of the fixed point set $M^T$ is such that $\mu_T(F)$ does not lie
on a wall through 0 (or a wall such that the affine hyperplane
spanned by the wall passes through 0), then the contribution of
$F$ to the residue formula for $\cali^\epsilon_0(\a\b
e^{i\baromega})$ is the same as the contribution of $F$ to the
pairing
$$\kappa_{{{M},\xi}}^T(\a \b
{\mathcal{D}}e^{i\baromega - \epsilon \Theta/2})[\mu^{-1}(\xi)/
T]$$ for $\xi$ sufficiently close to 0, which was calculated in
Remark \ref{fcknew}: it is zero when $F \notin \calf_+$ and
$$\frac{
 (-1)^{s+n_+}}{|W| \vol(T)} \res \bigl( \cald(\xvec)^2
 \int_{{F }}
\frac{i_{{F}}^*(\alpha \beta
e^{i(\bom-\delta}))(\xvec)e^{-\epsilon \langle X,X \rangle
/2}}{e_{F}(\xvec)} [d\xvec] \bigr )$$ for any sufficiently small
$\delta \in \liets$ which is a regular value of the moment map
$\mu_T$ when $F \in \calf_+$. When the degrees of $\alpha$ and
$\beta$ sum to the real dimension of $\mnd$ and $\alpha$ and
$\beta$ both restrict to elements of the subspace $V(n,d)$ of
$H^*_{\overline{\mathcal{G}}(n,d)}({\mathcal{C}}(n,d)^{ss})$ which
is isomorphic to $I\!H^*(\mnd)$ (see $\S$5 above) then this
contribution is also the same as the contribution of $F$ to the
residue formula for the pairing $\langle \k_{{}}(\a),\k_{{}}(\b
)\rangle$ in the intersection cohomology  of $\mnd$. If however
$\mu_T(F)$ does lie on a wall through 0 then the contribution of
$F$ to the residue formula for $\cali^\epsilon_0(\a\b
e^{i\baromega})$ involves Gaussian integrals over cones $C$ in
$\liets$  of the form (cf. \cite{jkkw}, (9.5))
\begin{equation} \label{e:gaussianfixed}
 \frac{i^l(2 \pi)^{-l/2} }{|W|\vol(T) \epsilon^{s/2}}
\int_{y \in C  } [dy]  \cald(y) e^{- <y, y>/{2\epsilon} } \res^{}
\Bigl (     \cald(X) \int_F\frac{i_F^* (\a(X) \b(X) e^{i \omega})
}{e_F(X) } e^{i <\mu(F)-y,X>} \Bigr )
\end{equation}
which reduces to the expression above when $C$ is the whole of
$\liets$.
\end{rem}

We now study $\cali^\epsilon (\alpha \beta e^{i \bar{\omega}})$
when $M = M^{\rm ext}_K$; for simplicity we will only consider the
case when $n=2$ and $d=0$. As before $T$ denotes the maximal torus
$U(1)$  of $SU(2)$, and we identify the Lie algebra of $U(1)$ with
$\RR$ by equating $X \in \RR$ with $iX \gamma$ where $\gamma =
(1,-1)$ is a chosen root of $SU(2)$. Note that $<\gamma, \gamma>
= 2$ in the  Euclidean inner product.

$M^{\rm ext}_K$ is of course not compact and there are infinitely
many components $F$ of the fixed point set of the action of $T$ on
$M^{\rm ext}_K$, which are indexed by the value of the moment map
(or equivalently by the integers), each diffeomorphic to $T^{2g}$.
But by Theorem 34 and Remark 36 of \cite{jkkw} (cf. Remark
\ref{jkkw1rem9.3} above), the difference between
$$\cali^\epsilon_0(\a\b  e^{i\baromega})$$
and the pairing $\kappa_{{{M},\xi}}^T(\a \b
{\mathcal{D}}e^{i\baromega - \epsilon \Theta/2})[\mu^{-1}(\xi)/
T]$ for any $\xi$ sufficiently close to 0 (which can be calculated
as at Remark \ref{fcknew})
 is a sum of contributions
corresponding to those components $F$ of the fixed point set of
$T$ for which $\mu_T(F)$ lies on a wall through 0 in $\liets \cong
\RR$, i.e. for which $\mu_T(F)=0$. In fact there is only one
such component $F_0 = T^{2g} \times \{0\} $ in the extended moduli
space $M^{\rm ext}_K $ (see $\S$4.2). The Euler class
$e_{F_0}(X)$ is given by $(2X)^{2g}$  as it is the product of $2g$
copies of the root $\gamma(X) = 2X$, and we have $\cald^2 (X) =
(2X)^2$ since $\cald(X) = \gamma(X).$ Hence we can conclude that
\begin{equation} \label{e:ieps}
\cali^\epsilon(\a\b  e^{i\baromega}) =  E_0  + E_1 + \frac{
(-1)^{g-1}   } {2  }  {\rm res}_{X = 0 } \Bigl ( \sum_{j \geq 0}
\frac{(-\epsilon X^2)^j}{j!}
 \frac{\int_{T^{2g}}\alpha \beta e^{i\omega} }  {(2X)^{2g-2} (e^{2X} - 1) }
\Bigr )
 \end{equation}
where $E_0 $ is the contribution of $F_0$ to
(\ref{e:gaussianfixed}), while
 $E_1$ is a sum of terms vanishing exponentially in $\epsilon$
as $\epsilon \to 0$ (see \cite{JK1} Theorem 4.1 for a more precise
definition of ``vanishing exponentially'') and the last term is
given by Remark \ref{fcknew}.

For simplicity we take $\a=\b  = 1$ from now on. Then the
contribution $E_0$ to the integral $\cali^\epsilon(
e^{i\baromega})$
 from  the component $F_0$
is given by
\begin{equation} \label{e:gauss2} E_0 =
C_0 \int_{X \in \RR - i \delta} \frac{e^{- \epsilon X^2}}{X^{2g-2} }
dX \end{equation} where
$$ C_0 = \frac{1}{4 \pi}  \frac{\int_{F_0}e^{i\omega} }{2^{2g-2}} $$
and $ \int_{F_0} e^{i \omega} = (2i)^g $ (see \cite{JK2}, Lemma
10.10). We have treated overall normalization constants using the
conventions of Corollary 8.2 of \cite{JK1}, with the correction
made in Footnote 9 of \cite{JK2}. Here $\delta > 0 $ is a small
real parameter which  ensures convergence of the integral
(\ref{e:gauss2}) despite the pole at $X = 0 $.

\begin{rem}One readily sees (by changing variables from $X$ to $Y  = \epsilon X^2$) that  the integral in
(\ref{e:gauss2}) is homogeneous in $\epsilon$. In fact
$$\int_{X \in \RR - i \delta}\frac{e^{- \epsilon X^2} } {X^{2g-2}} dX
=  \frac{\epsilon^{g-3/2}}{2} \int_{\epsilon Y^2 \in \RR - i
\delta} e^{-Y} Y^{1/2 - g} dY .$$
\end{rem}

In the situation studied in \cite{tdgr}, Witten found that
 the integral analogous to $\cali^\epsilon(e^{i\baromega})$ (which Witten denotes by $Z(\epsilon)$) is computed
when $n=2$ and $d=0$ as
\begin{equation}  \label{e:witans}
Z(\epsilon) = \frac{1}{ (2 \pi^2)^{g-1} } \sum_{n = 1}^\infty
\frac{\exp (-  \epsilon \pi^2 n^2) }{ n^{2g-2} } \end{equation}
(see \cite{tdgr}, (4.43)). It follows (using the Poisson summation
formula as in (4.53) of \cite{tdgr}) that
\begin{equation} \label{e:wit2}
\left ( \frac{\partial}{\partial \epsilon} \right )^{g-1}
Z(\epsilon) = C_1 \epsilon^{-1/2}  + E'_1 \end{equation} where
$C_1 = \frac{(-1)^{g-1}}{2^g \sqrt{\pi}} $ and $E'_1$  vanishes
exponentially as $\epsilon \to 0 $. Thus we have (integrating with
respect to $\epsilon$)  that (see \cite{tdgr}, (4.54))
\begin{equation} \label{e:wit3}
Z(\epsilon) = C_1 C_2 \epsilon^{g - 3/2} + E_1 + E_2
\end{equation}
where
\begin{equation} \label{e:c2}
C_2 = 2 \cdot (2/3) \cdot (2/5) \dots (2/(2g-3))  \end{equation}
and $E_1$ vanishes exponentially as $\epsilon \to 0$ while $E_2$
is a polynomial in $\epsilon$ of degree $g-2$.

We can now compare our expression for
$\cali^\epsilon(e^{i\baromega})$ to Witten's $Z(\epsilon)$. It is
easy to see that if $j \geq g-1$ then
$$
 {\rm res}_{X = 0
} \Bigl (
 \frac{(- X^2)^j\int_{T^{2g}} e^{i\omega}}  {j! (2X)^{2g-2} (e^{2X} - 1) }
\Bigr ) = 0.$$ Hence it follows from (\ref{e:ieps}) and
(\ref{e:gauss2}) that
\begin{equation} \label{e:gauss3}
 \left ( \frac{\partial}{\partial \epsilon} \right )^{g-1}\cali^\epsilon(e^{i\baromega}) =
(-1)^{g-1} C_0 \int_{X \in \RR - i \delta} dX e^{- \epsilon X^2} =
(-1)^{g-1} C_0 {\sqrt{\pi} }\epsilon^{-1/2} + E'_1
\end{equation}
where $E'_1$  vanishes exponentially as $\epsilon \to 0$ (cf.
\cite{JK1} Theorem 4.1). From this and Remark \ref{jkkw1rem9.3} we
obtain
\begin{equation} \label{e:7.8}
 \cali^\epsilon(e^{i\baromega}) = \frac{ (-1)^{g-1}
 i^g }{2^{g} \sqrt{\pi}  }  C_2\epsilon^{g-3/2}+ E_1 +
\sum_{j = 0}^{g-2}
 {\rm res}_{X = 0
} \Bigl (
 \frac{(-\epsilon X^2)^j\int_{T^{2g}} e^{i\omega}}  {j! (2X)^{2g-2} (e^{2X} - 1) }
\Bigr )
 \end{equation}
where $E_1$ is a sum of terms vanishing exponentially in
$1/\epsilon$ as $\epsilon \to 0$ and $C_2$ was introduced in
(\ref{e:c2}) above.

Up to multiplication by the constant $  i^g$, the coefficient of
$\epsilon^{g-3/2}$ in (\ref{e:7.8}) agrees with the coefficient of
$\epsilon^{g-3/2}$ in
 Witten's expression  (\ref{e:wit3})
for $Z(\epsilon)$. The factor $i^g$ is accounted for because
Witten computes
$$Z(\epsilon) = \frac{1}{(2 \pi)^s \vol(K)}
\int_{X \in \liek}  e^{- \epsilon X^2} e^{\omega + i \mu X} $$
which can be transformed into our $\cali^\epsilon$ by the
substitution $\omega \mapsto i \omega$.

We can also investigate the relation between the polynomial part
of $Z(\epsilon)$ and our expression
\begin{equation} \label{jkkwexp} \sum_{j = 0}^{g-2}
 {\rm res}_{X = 0
} \Bigl (
 \frac{(-\epsilon X^2)^j\int_{T^{2g}} e^{i\omega}}  {j! (2X)^{2g-2} (e^{2X} - 1) }
\Bigr )\end{equation} for the polynomial part of
$\mathcal{I}^\epsilon(e^{i\baromega})$. The formula (8.7) for
$Z(\epsilon)$ may be expanded to find the coefficient of
$\epsilon^k$ for $0 \leq k \leq g-2$ (cf. (4.49) of \cite{tdgr})
 as
\begin{equation}
\Bigl ( \frac{1}{2 \pi^2} \Bigr )^{g-1} \sum_{n = 1}^\infty
 \frac{ (- \pi^2)^k}
{ n^{2g-2-2k} } \end{equation} which equals
\begin{equation}  \label{e:epsexp}
\Bigl ( \frac{1}{2 \pi^2} \Bigr )^{g-1} {(- \pi^2 )^k} \zeta
(2g-2-2k)
\end{equation}
where $\zeta$ denotes the Riemann zeta function.

By an elementary contour integral argument (cf. \cite{rank2} Lemma
5.12) we see that
\begin{lemma} \label{l:l5.12} For any positive integer $m$, we have
$$\zeta(2m) = \sum_{n > 0 } \frac{1}{n^{2m}} = \pi i {\rm res}_{Y= 0}
\frac{1}{Y^{2m} (e^{2 \pi i Y} - 1) }$$
$$ = {(-1)^m (\pi)^{2m}} {\rm res}_{X = 0 } \frac{1}{X^{2m} (e^{2X} - 1)}. $$
\end{lemma}

Using Lemma \ref{l:l5.12} to express the  residue in
(\ref{jkkwexp}) as a multiple of a zeta function, together with
the fact from \cite{JK2}, Lemma 10.10 that $
\int_{T^{2g}}e^{i\omega} = (2i)^g$, we see that (\ref{jkkwexp})
agrees with (\ref{e:epsexp}).

\end{document}